\documentclass[11pt]{amsart}
\usepackage{eucal}
\usepackage[all]{xy}
\usepackage{amsfonts,amssymb}
\usepackage{epsfig}
\usepackage[usenames]{color}
\usepackage{float} 
\xyoption{arc}

\newtheorem{theorem}{Theorem}[section]
\newtheorem{lemma}[theorem]{Lemma}

\newtheorem{corollary}[theorem]{Corollary}
\newtheorem{defn}[theorem]{Definition}

\newcommand{\Q}{\mathbb Q}
\newcommand{\C}{\mathbb C}
\newcommand{\R}{\mathbb R}
\newcommand{\K}{\mathbb K}
\renewcommand{\P}{\mathbb P}
\renewcommand{\O}{\EuScript O}

\newcommand{\D}{\mathbb D}

\newcommand{\Z}{\mathbb Z}
\newcommand{\N}{\mathbb N}

\newcommand{\bracket}[1]{\{#1\}}

\title{The growth rate of symplectic homology and affine varieties}
\author{Mark McLean}

\begin{document}

\begin{abstract}
We will show that the cotangent bundle of a manifold whose free loopspace
homology grows exponentially is not symplectomorphic to any smooth affine variety.
We will also show that the unit cotangent bundle of such
a manifold is not Stein fillable by a Stein domain
whose completion is symplectomorphic to a smooth affine variety.
For instance, these results hold for end connect sums of
simply connected manifolds whose cohomology
with coefficients in some field has at least two generators.
We use an invariant called the growth rate of symplectic homology
to prove this result.
\end{abstract}

\maketitle

\bibliographystyle{halpha}


\tableofcontents

\section{Introduction}

The aim of this paper is to use an invariant
called the growth rate of symplectic homology
to construct many examples of cotangent bundles
which are not symplectomorphic
to any smooth affine variety.
We also have a contact analogue of this result which says
that if we take the unit cotangent bundle of one of these
manifolds then it 
is not Stein fillable by a Stein domain whose completion
is symplectomorphic smooth affine variety.
All of our manifolds are assumed to be oriented
unless explicitly stated otherwise.

If we have some Liouville domain $M$ (defined in Section
\ref{section:liouvilldomaindefinition}) and a class $b \in H^2(M,\Z / 2\Z)$,
then the growth rate of symplectic homology is
an invariant $\Gamma(M,b) \in \{-\infty\} \cup [0,\infty)$.
This invariant was originally defined in
\cite[Section 4]{Seidel:biasedview}.
Every Liouville domain also has a completion $\widehat{M}$,
and the growth rate is in fact an invariant of the
completion $\widehat{M}$ up to symplectomorphisms preserving the class $b$
and so we will often write $\Gamma(\widehat{M},b)$.

A smooth affine variety $A$ has a symplectic structure
obtained by embedding it into $\C^N$ algebraically
and then pulling back the standard symplectic form onto $A$.
It turns out that this is symplectomorphic to the completion
of some Liouville domain $\overline{A}$ obtained from
taking a large closed ball in $\C^N$ and intersecting it with $A$.
Because the symplectic form is a biholomorphic invariant of $A$,
we can assign the invariant $\Gamma(A,b) = \Gamma(\overline{A},b)$
(see \cite{EliahbergGromov:convexsymplecticmanifolds}).

We say that a contact manifold $C$ is {\it algebraically Stein fillable}
if it is Stein fillable by a Stein domain $D$ whose completion
$\widehat{D}$ is symplectomorphic to a smooth affine variety $A$.
For a smooth affine variety $A$ we also have a finite
invariant $m_A \in \N_{>0}$ defined as follows:
choose a compactification $X$ of $A$ by a smooth normal crossing
divisor (i.e. $X \setminus A$ is a union of
smooth transversely intersecting connected complex hypersurfaces $S_i$,
$i=1,\cdots,k$).
For $I \subset \{1,\cdots,k\}$ we write $S_I := \cap_{i \in I} S_i$.
Let $d := \text{max} \{ n - \text{dim}_\C(S_I) | S_I \neq \emptyset\}$
where $n = \text{dim}_\C A$. 
We define $m_A$ to be the minimum of
$d$ over all compactifications described above.

The main result of this paper is:
\begin{theorem} \label{thm:maingrowthratebound}
$\Gamma(A,b) \leq m_A$ for any class $b \in H^2(A,\Z/2\Z)$.
In particular $\Gamma(A,b) < \infty$.
\end{theorem}

In fact we will prove a more general theorem:
\begin{theorem} \label{thm:mainfillingobstruction}
Suppose that the boundary of a Liouville domain $M$
is algebraically Stein fillable by $A$, then $\Gamma(M,b) \leq m_A$
for any $b \in H^2(M,\Z/2\Z)$.
\end{theorem}
This theorem will be proven in Section
\ref{subsection:smoothaffinevarieities}.
These theorems have been proven in
\cite[Section 4]{Seidel:biasedview} when $A$ has complex
dimension $2$.
Let $Q$ be a compact manifold.
Choose some metric on $Q$ and let ${\mathcal L}^q Q$
be the space of all free loops on $Q$ of length less than or
equal to $q$, and let ${\mathcal L} Q$ be the space of
all free loops on $Q$.
Let $\K$ be any field and let $a^{\K}(q)$ be the rank of the image of the
inclusion map
\[H_*({\mathcal L}^q(Q),\K) \hookrightarrow H_*({\mathcal L}(Q),\K).\]
The manifold $Q$ is said to have {\it exponential growth} if for some $\K$,
the function $a^\K(q)$ grows faster than any polynomial
(we mean here that for any polynomial $p(q)$,
$a^\K(q) > p(q)$ for $q$ large enough).
If the manifold $Q$ has finite fundamental group,
then this definition of exponential growth is the
same as saying that $\sum_{i=1}^k b_i({\mathcal L}(Q),\K)$
grows faster than any polynomial in $k$ for some $\K$
by \cite{Gromov:homotopical}.
There are many examples of manifolds with exponential growth. Here are some:
\begin{enumerate}
\item Any surface of genus $2$ or higher.
\item The end connect sum of two simply connected manifolds
$M_1 \# M_2$ where $H^*(M_1,\K)$ has at least two generators
and $H^*(M_2,\K)$ has rank at least $3$ (See \cite{Lambrechts:Betti}).
\item Any manifold $Q$ whose fundamental group
is the free product of at least 3 non-trivial groups
(See Lemma \ref{lemma:growthratefundamentalgroup},
and Lemma \ref{lemma:triplefreeproduct}).
\end{enumerate}

Conjecturally (see \cite{ViguePoirrier:homotopie}) there should be many
more manifolds with exponential growth such as simply connected manifolds whose
Betti numbers are greater than that of the torus.
Given any Riemannian manifold $Q$, we can construct a symplectic manifold
$T^*Q$ where locally the symplectic form is $\sum_i dq_i \wedge dp_i$
where $q_i$ are the position coordinates and $p_i$ are the momentum
coordinates.
We also have a Liouville domain $D^* Q$
which is the manifold with boundary consisting of
covectors of length less than or equal to $1$.
Its boundary is the contact manifold denoted by $S^*Q$
consisting of covectors of length $1$ and the contact
form here is $\sum_i p_i dq_i$.

By using work from 
\cite[Corollary 1.2]{salamonweber:loop}
(although we could have used ideas from
\cite{AbbondandoloSchwartz:cotangentloop}
or \cite{Viterbo:functorsandcomputations2}),
it can be shown that
$\Gamma(T^*Q,\omega_2) = \infty$ when $Q$ has exponential growth
and $\omega_2$ is the pullback of the second Stiefel-Whitney class of $Q$.
Hence we have the following corollary of Theorem
\ref{thm:maingrowthratebound}:
\begin{corollary} \label{corollary:cotangentaffine}
If $Q$ has exponential growth then $T^*Q$ is not symplectomorphic
to any smooth affine variety.
\end{corollary}
This is because $\Gamma(T^*Q,\omega_2) = \infty$ and $\Gamma(A,b)<\infty$
for any smooth affine variety $A$ and any $b \in H^2(A,\Z / 2\Z)$.
We also have the following corollary of Theorem
\ref{thm:mainfillingobstruction}:
\begin{corollary} \label{corollary:unitcotangent}
The unit cotangent bundle of $Q$ is not algebraically
Stein fillable.
\end{corollary}

If we wish to consider unoriented manifolds $Q$ then we need to do the following:
A contact manifold $C$ is said to be {\it covered algebraically Stein fillable}
if there is a Stein filling $M$ of $C$ whose completion $\widehat{M}$
is symplectomorphic to a finite cover of a smooth affine variety.

\begin{theorem} \label{thm:coverfillingobstruction}
Suppose the boundary of a Liouville domain $M$
is covered algebraically Stein fillable by a smooth affine variety $A$,  then $\Gamma(M,b) \leq m_A$
for any $b \in H^2(M,\Z/2\Z)$.
\end{theorem}

We will prove this theorem in Section
\ref{section:coversofsmoothaffinevarieties}.
We then immediately get the following corollaries:
\begin{corollary} \label{corollary:covercotangentaffine}
If $Q$ is a compact manifold of exponential growth which is possibly unoriented
such that some finite cover is oriented and is has exponential growth then
$T^*Q$ is not symplectomorphic to any smooth affine variety.
\end{corollary}

This corollary is true for the following reason:
If $T^* Q$ is symplectomorphic to some smooth affine variety $A$
then for every oriented finite cover $\widetilde{Q}$ of $Q$,
$T^* \widetilde{Q}$ is symplectomorphic to some finite cover of $A$.
Hence by the above theorem we have that $\Gamma(T^* \widetilde{Q},\omega_2) < \infty$
which means that $\widetilde{Q}$ is cannot have exponential growth.

Similar reasoning also gives us the following corollary:
\begin{corollary} \label{corollary:coverunitcotangent}
The unit cotangent bundle of any (possibly unoriented) $Q$ of exponential growth with dimension $>2$ is not
covered algebraically Stein fillable.
In particular the unit cotangent bundle of any such $Q$ is not
algebraically Stein fillable.
\end{corollary}
The reason why we need dimension $>2$ is that if we have some algebraic Stein filling $C$
of $S^* Q$ then we require that $\pi_1(Q) = \pi_1(S^*Q) = \pi_1(C)$ so that we can take
an appropriate cover of $C$.
In future work the author hopes to prove a based loopspace
version of the above theorems using wrapped Floer homology.
This will enlarge the class of manifolds whose cotangent bundle
is not symplectomorphic to a smooth affine variety.
For instance this would be true
if their Betti numbers are greater than that of the torus.

There is a similar conjecture relating smooth affine varieties
with cotangent bundles which we will now describe.
Suppose that $U$ is a real affine algebraic variety
over $\R$ which is diffeomorphic to our manifold $Q$.
We will suppose that $Q$ is simply connected.
We say that $Q$ has a {\it good complexification} if there exists
a $U$ as described above such that the natural map $U(\R) \rightarrow U(\C)$
is a homotopy equivalence.
There is a question in \cite{Totaro:complexifications}
which asks if $Q$ has a Riemannian metric of non-negative sectional curvature
when $Q$ has a good complexification.
There is also a conjecture attributed to Bott which conjectures
that any simply connected manifold with non-negative sectional curvature is
rationally elliptic. Rationally elliptic means that
$\text{dim}\pi_*(Q) \otimes \Q < \infty$
(see \cite{FelixHalperinThoamas:homotopylie}).
This conjecture is mentioned in
\cite[Question 12, page 519]{FelixHalperinThomas:rationalhomotopytheory}.
Felix and Thomas in \cite{FelixThomasRadius} proved that if a simply connected manifold is rationally elliptic
then the growth rate of the Betti numbers of its based loopspace grows sub exponentially.
Vigu{\'e}-Poirrier in \cite[Page 415]{ViguePoirrier:homotopie}
conjectured that the growth rate of the Betti numbers of the based loopspace grow exponentially
if and only if the same thing is true for the free loopspace.
If we look at all of these statements then we get the following question:
Suppose $Q$ is simply connected of exponential growth, then is it true that
$Q$ does not admit a good complexification?
This question has some similarities with Corollary
\ref{corollary:cotangentaffine}
because we can choose a symplectic form on $U(\C)$
such that a neighbourhood of $U(\R)$ is symplectomorphic to
a neighbourhood of the zero section of $T^*Q$ and also
because both $T^*Q$ and $U(\C)$ deformation retract onto $Q$.

We will now give a brief sketch of the proof of
Theorem \ref{thm:mainfillingobstruction}.
Symplectic homology of $\widehat{M}$ is the homology of a chain complex having the following
generators:
\begin{enumerate}
\item Critical points of some Morse function on $M$.
\item Two copies of each closed Reeb orbit of $\partial M$.
\end{enumerate}
This chain complex has a natural filtration by $\R_{\geq 0}$ called the {\it action}.
Critical points have action $0$ and the action of each Reeb orbit is its length.
Hence we can define $SH_*^{\leq \lambda}(M)$ for each $\lambda$ as the
homology of the subcomplex of orbits of action $\leq \lambda$.
The growth rate is the rate at which
\[\text{rank im}(SH_*^{\leq \lambda}(M) \rightarrow \varinjlim_l SH_*^{\leq l}(M)) \]
grows with respect to $\lambda$.
For instance if the rank is bounded above by some polynomial of degree $k$
then the growth rate is less than or equal to $k$.
It turns out that our smooth affine variety $A$ is symplectomorphic to $\widehat{M}$ for some $M$.
We need to provide an upper bound for the growth rate of $M$.
One way of doing this is to show that the number of Reeb orbits
of length $\leq \lambda$ is bounded above by some polynomial of degree
at most $m_A$.
This is the method used in \cite[Section 4]{Seidel:biasedview}
when Theorem \ref{thm:mainfillingobstruction} is proven in complex dimension $2$.
This upper bound is achieved by compactifying $A$ by a smooth normal crossing
divisor and carefully constructing a Liouville domain using this divisor.
The problem is that this upper bound is much harder to achieve in higher dimensions.
In this paper we define growth rate in a slightly more flexible
way so that it is easier to provide such an upper bound.

We have some Hamiltonian $H \geq 0$ and an almost complex structure $J$
called a {\it growth rate admissible pair} which satisfies certain technical properties.
From this pair $(H,J)$ we can define a series of homology groups
$SH_*^{\#}(\lambda H,J)$ generated by $1$-periodic orbits of $\lambda H$.
For $\lambda_1 \leq \lambda_2$ there is a natural map
\[SH_*^{\#}(\lambda_1 H,J) \rightarrow SH_*^{\#}(\lambda_2 H,J).\]
We prove that the growth rate of
\[\text{rank im}(SH_*^{\#}(\lambda H,J) \rightarrow \varinjlim_l SH_*^{\#}(l H,J))\]
is equal to the original definition of growth rate (Section \ref{subsection:alternategrowthratedefinition}).
Finally we carefully construct a Hamiltonian $H$ so that
there is a degree $m_A$ polynomial bound on the number of $1$-periodic orbits
of a very small perturbation of $\lambda H$ and this gives us an upper
bound on the growth rate (Section \ref{subsection:smoothaffinevarieities}).
We construct this Hamiltonian by first compactifying our smooth affine variety $A$
with a smooth normal crossing divisor $D$.
We deform $D$ without changing the symplectic structure on its complement $A$
so that there is a nice symplectic structure near $D$ (Section \ref{secion:compactificationsofalgebraicvarieties}).
We can construct $H$ so that around some point $p \in D$ we have that
$D$ looks like $\{x_1=0,\cdots,x_k=0,y_1=0,\cdots,y_k=0\}$
where $x_1,\cdots,x_n,y_1,\cdots,y_n$ is some nice coordinate chart around $p$
and $H$ is a product Hamiltonian.
The symplectic structure is not quite the standard one on this chart
but the Hamiltonian flow of $H$ is exactly the same as the Hamiltonian
flow with respect to the standard symplectic structure.
This enables us to find all the $1$-periodic orbits of $H$ and give
a bound proportional to $\lambda^k \leq \lambda^{m_A}$ in each chart.
Hence we get an upper bound of $m_A$.

{\bf Acknowledgements:}
I would like to thank Mohammed Abouzaid,
Dietmar Salamon, Burt Totaro, Paul Seidel and Ivan Smith for extremely
useful comments.
In addition, the author is grateful to the referee for many valuable comments and suggestions.
The author was partially supported by
NSF grant DMS-1005365.

\section{Main definitions and properties}
\subsection{Liouville domains}
\label{section:liouvilldomaindefinition}

A {\it Liouville domain} is a compact manifold $N$
with boundary and a $1$-form $\theta_N$ satisfying:
\begin{enumerate}
\item $\omega_N := d\theta_N$ is a symplectic form.
\item The $\omega_N$-dual of $\theta_N$ is transverse
to $\partial N$ and pointing outwards.
\end{enumerate}
We will write $X_{\theta_N}$ for the $\omega_N$-dual of $\theta_N$.
Sometimes we have manifolds with corners with $1$-forms
$\theta$ satisfying the same properties as above.
We view these as Liouville domains by smoothing the corners slightly.
By flowing $\partial N$ backwards along $X_{\theta_N}$ we have a small collar
neighbourhood of $\partial N$ diffeomorphic to $(0,1] \times \partial N$
with $\theta_N = r_N \alpha_N$.
Here $r_N$ parameterizes $(0,1]$ and $\alpha_N$ is the contact form on the boundary
given by $\theta_N |_{\partial N}$.
The completion $\widehat{N}$ is obtained
by gluing $[1,\infty) \times \partial N$ to this collar neighbourhood
and extending $\theta_N$ by $r_N \alpha_N$.

The Liouville domains that we will be studying are called
Stein domains. A {\it Stein manifold}
is a complex manifold that can be properly embedded in $\C^N$.
An equivalent definition of a Stein manifold is a complex manifold
admitting an exhausting plurisubharmonic function.
A function is {\it exhausting} if it is bounded
from below and the preimage of every compact set is compact.
A function $f : S \rightarrow \R$ is {\it plurisubharmonic}
if $(-d d^c f)(X,JX) > 0$ for all non-zero vectors $X$
where $d^c := d \circ J$ and $J$ is the complex structure on the complex manifold $S$.
This implies that $\omega_f := -dd^cf$ is a symplectic form.
If $c$ is a regular value of an exhausting plurisubharmonic function $f$ then $f^{-1}(-\infty,c]$
is a Liouville domain with Liouville form $-d^c f$.
We call such a domain a {\it Stein domain}.
A Stein manifold is of {\it finite type}
if it admits an exhausting plurisubharmonic function with only finitely
many singularities.
A plurisubharmonic function $f$ is said to be {\it complete} if the $\omega_f$-dual
of $-d^cf$ is integrable. Every finite type Stein manifold
admits an exhausting plurisubharmonic function which is complete.
Also any two such functions on $S$ give symplectomorphic Stein manifolds.
Hence the symplectic form is a biholomorphic invariant.
We will write $\omega_S$ for this symplectic form.
If $f$ is finite type then it turns out for $C \gg 0$
that $(S,\omega_S)$ is symplectomorphic to the completion of $f^{-1}(-\infty,C]$.
One way of obtaining such a symplectic form is embedding $S$
as a closed proper holomorphic embedding into $\C^N$ and then pulling
back the standard symplectic form $\omega_{\text{std}}$ on $\C^N$.
Important examples of finite type Stein manifolds $A$ are
smooth affine varieties (see \cite[Section 4b]{Seidel:biasedview}).
In fact from \cite[Section 4b]{Seidel:biasedview} we have
for any smooth affine variety $A$, a Liouville domain $\overline{A}$
which unique up to isotopy through Liouville domains such that
$\widehat{\overline{A}}$ is symplectomorphic to $A$.
Uniqueness here means that if $B$ is isomorphic as a variety
to $A$ then $\overline{A}$ is isotopic to $\overline{B}$.
From now on throughout this paper if we have a smooth affine variety $A$ we will
write $\overline{A}$ for such a Liouville domain.
We will call any such Liouville domain an {\it algebraic Stein domain}.
Here is a direct way of constructing this Liouville domain:
Choose any algebraic embedding $\iota$ of $A$ into $\C^N$
(so it is a closed subvariety).
Let $(r_i,\vartheta_i)$ be polar coordinates for the $i$'th
factor in $\C^N$.
Let $R := \sum_i \frac{r_i^2}{4}$. We have $\theta_A := -d^c R = \sum_i\frac{r_i^2}{2} d\vartheta_i$.
We have that $d\theta_A$ is equal to the standard symplectic structure on $\C^N$.
By abuse of notation we write $\theta_A$ for $\iota^* \theta_A$,
and $\omega_A := d\theta_A$.
\begin{lemma} \label{lemma:affineLiouvilledomain}
There is a $C > 0$ such that
for all $c \geq C$,$((R|_A)^{-1}(-\infty,c],\theta_A)$
is a Liouville domain whose completion is symplectomorphic to $(A,\omega_A)$.
\end{lemma}
\proof of Lemma \ref{lemma:affineLiouvilledomain}.
We will first show that for some $C \geq 0$, $R|_A$ has no singularities
when $R|_A \geq C$.
Compactify $\C^N$ to $\P^N$ and let $X$ be the projective subvariety
which is the closure of $A$ in $\P^N$.
Let $H$ be the hypersurface $X \setminus A$.
We have a line bundle $O(1)$ described as follows:
If we think of a point on $\P^N$ as a line $\Lambda$ in $\C^{N+1}$
then the fiber at this point is the quotient $(\C^{N+1})^* / V$
where $V$ is the vector subspace of linear $1$-forms whose kernel contains $\Lambda$.
The hyperplane $H$ is represented by some linear hyperplane $\widetilde{H}$ in $\C^{N+1}$
and so we have a section $s$ given by a non-zero linear $1$-form
whose kernel is $\widetilde{H}$. This section vanishes on $H$
with order $1$.

Another way of describing the fiber of this bundle over the point represented by the line $\Lambda$
is as the space of linear $1$-forms $W_{\Lambda}$ which vanish in the hyperplane orthogonal to $\Lambda$.
We can put a metric $\|\cdot\|$ on $W_{\Lambda}$ induced by the standard
metric on $(\C^{N+1})^*$ hence we view $\|\cdot\|$ as a metric on $O(1)$.
The curvature form is a positive $(1,1)$ form
and this gives us a symplectic form on $\P^N$.
We have an action of $U(N+1)$ on $\P^N$ induced by the standard
action on $\C^{N+1}$. This action also naturally lifts to an action
on the total space of $O(1)$.
The metric $\|\cdot\|$ is invariant under the action of $U(N+1)$.
We can ensure that $U(N) \subset U(N+1)$ is the subgroup which sends $H$ to $H$.
This subgroup also fixes our section $s$ so
$u^*\|s\| = \|s\|$ for each $u \in U(N)$.
Hence $\|s\|$ must be a function of $R$. If it wasn't then there exists
two elements $a,b \in \C^N$ of the same modulus such that $\|s\|(a) \neq \|s\|(b)$
but this is impossible as there is a $u \in U(N)$ such that $u(a) = b$ which implies
that $\|s\|(a) = u^*(\|s\|)(b) = \|s\|(b)$.
The function $-\log{\|s\|} |_A$ is a plurisubharmonic function
on $A$.
It is also equal to $f(R|_A)$ for some function $f$.
We will first show that $-\log{\|s\|} |_A$ has no singularities near infinity.
If $D$ was a smooth normal crossing divisor then \cite[Section 4b]{Seidel:biasedview}
tells us that $d\log{\|s\|} |_A$ is non-zero near infinity.
But the problem is that $D$ could be anything.
So by the Hironaka resolution of singularities
\cite{hironaka:resolution} we blow up $X$ away from $A$
to $\tilde{X}$ where $\tilde{X}$ is smooth and the
pullback $\tilde{D}$ is a smooth normal crossing divisor.
Let $L$ be the pullback of $\O(1)$ and $\tilde{s}$ the pullback of $s$.
We can also pull back the metric $\|.\|$ to $L$.
To show that $d\log{\|s\|}$ is non-zero we just show
that $d\log{\|\tilde{s}\|}$ is non-zero in exactly the same
way as in Lemma \cite[Section 4b]{Seidel:biasedview}:
Let $p \in \tilde{D}$ and choose local holomorphic coordinates
$z_1,\cdots,z_n$ and a trivialization of $L$ around $p$ so that
$s = z_1^{w_1} \cdots z_n^{w_n}$ ($w_i \geq 0$).
The metric $\|.\|$ on $L$ is equal to $e^{\psi}|.|$
for some function $\psi$ with respect to this trivialization
where $|.|$ is the standard metric on $\C$.
So \[d\log{\|\tilde{s}\|} = -\psi - (\sum_i w_i d\log{|z_i|}).\]
If we take the vector field
$Y := -z_1 \partial_{z_1} \cdots - z_n \partial_{z_n}$,
then $d\log{(|z_j|)}(Y) = -1$ and $Y.\psi$ tends to zero.
Hence $d\log{\|\tilde{s}\|}$ is non-zero near infinity which implies that
$f(R) = -\log{\|s\|}$ has no singularities near infinity.
Because $f(R)$ and $R$ are both exhausting, we get that
$f'>0$ near infinity which implies that $R|_A$ has no singularities
in the region $R \geq C$  for some $C>0$.

Because $A$ is a holomorphic submanifold
of $\C^N$, we have that $\omega_A$ is a symplectic form.
Let $c \geq C$ and write  $A_c := (R|_A)^{-1}(-\infty,c]$.
We have that $A_c$ is a Liouville domain because if $X_{\theta_A}$
is the $\omega_A$-dual of $\theta_A$ then it is equal to the
gradient of $R|_A$ and hence is transverse to $\partial A_c$ and pointing outwards
because $c$ is a regular value of $R|_A$ and $f'(c)>0$.
Also we have that the $\omega_A$-dual $X_{\theta_A}$ of $\theta_A$
is an integrable vector field because
$d(R|_A)(X_{\theta_A}) \leq dR(\sum \frac{r_i}{2} \frac{\partial}{\partial r_i})|_A = R|_A$
and hence if $\Phi_t^{X_{\theta_A}}$ is the flow of this vector field then
$(R|_A)(\Phi_t^{X_{\theta_A}}(x))$ increases at a rate of at most $e^R|_A$.
This first inequality is true because $X_{\theta_A}$ is the
orthogonal projection of $\sum_i \frac{r_i}{2} \frac{\partial}{\partial r_i}$
onto $TA$ so its length $l$ decreases, but
$R(\sum_i \frac{r_i}{2} \frac{\partial}{\partial r_i}) = \text{sup}_{|V|\leq l}
R(V)$.
By flowing $\partial A_c$ along $X_{\theta_A}$
we get that $\overline{A \setminus A_c}$
is diffeomorphic to $[1,\infty) \times \partial A_c$ and $\theta_A = r_A \alpha_{A_c}$
where $r_A$ parameterizes $[1,\infty)$ (this is because $d(R|_A)(X_{\theta_A}) > 0$ for $R|_A \geq C$).
Hence $A$ is symplectomorphic to $\widehat{A_c}$.
\qed

A contact manifold $C$ is said to be {\it Stein fillable}
if it is contactomorphic to the boundary of some
Stein domain.
It is said to be {\it algebraically Stein fillable}
if the completion of this Stein domain is symplectomorphic to $A$.

\subsection{Symplectic homology} \label{section:symplectichomology}

Let $N$ be a Liouville domain with $c_1 = 0$.
We make some additional choices $\eta := (\tau,b)$ for $N$.
The element $\tau$ is a choice of trivialization of the canonical
bundle of $N$ up to homotopy and $b$ is an element of $H^2(N,\Z / 2\Z)$.
We will assume that $\partial N$ has discrete period spectrum ${\mathcal P}_N$
(the set of lengths of Reeb
orbits of the $(\partial N',\alpha_{N'})$).
For each pair of numbers $a<b$ where $a,b \in [-\infty,\infty]$
we will define a symplectic homology group $SH_*^{(a,b]}(N,\eta)$
which is an invariant of $(N, \theta_N)$ up to exact symplectomorphism (although we suppress $\theta_N$ in
the notation unless the context is unclear).
Here an exact symplectomorphism between two Liouville domains
$N$ and $N'$ is a smooth diffeomorphism $\Phi$ from $N$ to $N'$
such that $\Phi^* \theta_{N'} = \theta_N + df$ for some smooth function $f: N \rightarrow \R$.

A family of Hamiltonians $H : S^1 \times \widehat{N} \rightarrow \R$
is said to be {\it admissible} if $H(t,x) = \lambda r_N(x)$
near infinity where $r_N$ is the cylindrical coordinate of $\widehat{N}$ and $\lambda$
is some positive constant. Here $\lambda$ is some positive constant
which is not in the period spectrum of $\partial N$.
We also require that $H_t|_N < 0$.
We sometimes view $H$ as a time dependent Hamiltonians $H_t : \widehat{N} \rightarrow \R$ where $t \in S^1$.
We have an $S^1$ family of vector fields $X_{H_t}$ and it
has an associated flow $\Phi^t_{X_{H_t}}$
(a family of symplectomorphisms parameterized by $t \in \R$ satisfying
$\frac{\partial}{\partial t} \Phi^t_{X_{H_t}} = X_{H_t}$
where we identify $S^1 = \R / \Z$).
A $1$-periodic orbit $o : S^1 \rightarrow \widehat{N}$
is a map which satisfies $o(t) = \Phi^t_{X_{H_t}}(x)$ for some $x \in \widehat{N}$.
We say that $o$ is non-degenerate if
$D\Phi^1_{X_{H_t}} : T_x\widehat{N} \rightarrow T_x\widehat{N}$
has no eigenvalue equal to $1$.
By \cite{DostoglouSalamon:instantonscurves},
we can perturb $H$ by a $C^\infty$ small amount so that
all of its $1$-periodic orbits are non-degenerate.
The problem here is that this Hamiltonian may not be admissible
after perturbing it, so we need a lemma:
\begin{lemma} \label{lemma:perturbinghamiltonian}
Let $H : S^1 \times N \rightarrow \R$ where $N$ is a symplectic manifold.
Let $U$ be a small neighbourhood of some of the $1$-periodic
orbits of $H$ such that no $1$-periodic orbits intersect
the boundary $\partial U$.
Then we can perturb $H$ by a $C^\infty$ small amount to $\tilde{H}$
such that all the $1$-periodic orbits of $\tilde{H}$ are non-degenerate
inside $U$ and $\tilde{H} = H$ outside $U$.
\end{lemma}
\proof of Lemma \ref{lemma:perturbinghamiltonian}.
We can choose $U$ small enough so that the time
$1$ flow of $X_{H_t}$ is well defined.
Let $U'$ be a smaller open neighbourhood of the $1$-periodic
orbits whose closure is contained in $U$.
By slightly extending the work of \cite{DostoglouSalamon:instantonscurves},
we can find for any positive function $f : N \rightarrow \R$,
a Hamiltonian $H_f$ such that
$\|H_f - H\|_{C^\infty} < f$ and such that
all the $1$-periodic orbits of $H_f$ are non-degenerate.
%
Let $R$ be an open neighbourhood of $\partial U'$ inside $U$
where $H$ has no $1$-periodic orbits which intersect
the closure of $R$.
Choose a bump function function $\rho : N \rightarrow \R$
which is $0$ on a neighbourhood of $U' \setminus R$
and $1$ on a neighbourhood of $N \setminus U'$.
Let $\tilde{H}_f := \rho H + (1-\rho) H_f$.
Suppose for a contradiction that for every $f$ sufficiently small,
there exists a $1$-periodic orbit of $\tilde{H}_f$
inside $U$ which is degenerate.
This orbit must intersect $R$
because all orbits of $H_f$ are non-degenerate away from $R$.
Then by a compactness argument, we have
a sequence of such orbits converging
to a $1$-periodic orbit of $H$ which intersects the closure of $R$.
This is impossible. Hence we have
for any $f$ arbitrarily small,
there is a Hamiltonian $H_f$ satisfying
$\|H_f - H\|_{C^\infty} < f$ such that
all of its orbits inside $U$ are non-degenerate
and $H_f = H$ outside $U$.
\qed

Because we have a trivialization $\tau$ of the canonical bundle
of $N$,
this gives us a canonical
trivialization of the symplectic bundle $TN$
restricted to an orbit $o$. Using this trivialization, we can define an index of $o$ called the
{\bf Robbin-Salamon} index (This is equal to the Conley-Zehnder
index taken with negative sign).
We will write $i(o)$ for the index of this orbit $o$.
For a 1-periodic orbit $o$ we define the {\bf action} $A_H(o)$:
\[A_H(o) := -\int_0^1 H(t,\gamma(t))dt -\int_o \theta_N.\]

Choose a coefficient field $\K$.
Let
\[CF_k^d(H,J,\eta) := \bigoplus_{o} \K \langle o \rangle\]
where we sum over $1$-periodic orbits $o$ of $H$ satisfying
$A_H(o) \leq d$ whose Robbin-Salamon index is $k$.
We write 
\[CF_k^{(c,d]}(H,J,\eta) := CF_k^d(H,J,\eta) / CF_k^c(H,J,\eta).\]
We need to define a differential for the chain
complex $CF_k^d(H,J,\eta)$ such that the natural inclusion
maps $CF_k^c(H,J,\eta) \hookrightarrow CF_k^d(H,J,\eta)$ for $c<d$ are chain maps.
This makes $CF_k^{(c,d]}(H,J,\eta)$ into a chain complex as well.
In order to define this we choose an ${\mathbb S}^1$ family of
almost complex structures $J_t$ compatible with
the symplectic form.
We assume that
$J_t$ is convex with respect to this cylindrical end outside some large compact set
(i.e. $\theta \circ J_t = dr$). We also say that $J_t$ is {\bf admissible}.

We will now describe the differential \[\partial : CF_k^d(H,J,\eta) \rightarrow CF_{k-1}^d(H,J,\eta).\]
We consider curves
$u : \R \times {\mathbb S}^1 \longrightarrow \widehat{N}$ satisfying
the perturbed Cauchy-Riemann equations:
\[ \partial_s u + J_t(u(s,t)) \partial_t u = \nabla^{g_t} H\]
where $\nabla^{g_t}$ is the gradient associated to the ${\mathbb S}^1$ family of metrics
$g_t := \omega(\cdot,J_t(\cdot))$.
For two periodic orbits $o_{-},o_{+}$ let
$\overline{U}(o_{-},o_{+})$ denote the set of all curves $u$ satisfying
the Cauchy-Riemann equations such that $u(s,\cdot)$ converges
to $x_{\pm}$ as $s \rightarrow \pm \infty$. This has a natural
$\R$ action given by translation in the $s$ coordinate.
Let $U(o_{-},o_{+})$ be equal to $\overline{U}(o_{-},o_{+}) / \R$.
For a $C^{\infty}$ generic admissible complex structure
we have that $U(o_{-},o_{+})$ is an $i^C(o_{-},o_{+}) -1$  dimensional
oriented manifold (see \cite{FHS:transversalitysymplectic}).
There is a maximum principle which ensures
that all elements of $U(o_{-},o_{+})$ stay inside
a compact set $K$ (see \cite[Lemma 1.5]{Oancea:survey},
\cite[Lemma 7.2]{SeidelAbouzaid:viterbo} or
Corollary \ref{corollary:maximumprinciplerescaling}).
Hence we can use a compactness theorem (see for instance
\cite{BEHWZ:compactnessfieldtheory}) which ensures that
if $i(o_{-}) - i(o_{+}) = 1$, then
$U(o_{-},o_{+})$ is a compact zero dimensional manifold.
The class $b \in H^2(N,\Z / 2\Z)$ enables us to orient this manifold
(see \cite[Section 3.1]{Abouzaid:contangentgenerate}).
Let $\# U(x_{-},x_{+})$ denote the number of
positively oriented points of $U(x_{-},x_{+})$
minus the number of negatively oriented points. Then we have a differential:
\[\partial : CF_k^d(H,J,\eta) \longrightarrow CF_{k-1}^d(H,J,\eta) \]
\[\partial \langle o_{-} \rangle := \displaystyle \sum_{i(o_{-}) - i(o_{+})=1 } \# U(o_{-},o_{+}) \langle o_{+} \rangle\]
By analyzing the structure of 1-dimensional moduli spaces, one shows
$\partial^2=0$ and defines
$SH_*(H,J,\eta)$ as the homology of the above chain complex.
As a $\K$ vector space $CF_k^d(H,J,\eta)$ is independent of $J$ and $b$,
but its boundary operator does depend on $J$. The homology
group $SH_*^d(H,J,\eta)$ depends on $H$ and $\eta$ but is independent of $J$
up to canonical isomorphism.
We define $SH_*^{(c,d]}(H,J,\eta)$ as
the homology of the chain complex \[CF_*^d(H,J,\eta) / CF_*^c(H,J,\eta).\]

If we have two non-degenerate admissible
Hamiltonians $H_1 < H_2$ and two admissible
almost complex structures $J_1,J_2$, then there is a natural
map:
\[SH_*^{(a,b]}(H_1,J_1,\eta) \longrightarrow SH_*^{(a,b]}(H_2,J_2,\eta)\]
This map is called a {\bf continuation map}.
This map is defined from a map $C$ on the chain level as follows:
\[C : CF_k^d(H_1,J_1,\eta) \longrightarrow CF_k^d(H_2,J_2,\eta) \]
\[\partial \langle o_{-} \rangle := \displaystyle
\sum_{i(o_{-}) = i(o_{-})) } \# P(x_{-},x_{+}) \langle o_{+} \rangle\]
where $P(o_{-},o_{+})$ is a compact oriented zero dimensional manifold
of solutions of the following equations:
Let $K_s$ be a smooth non-decreasing family of admissible Hamiltonians 
equal to $H_1$ for $s \ll 0$ and $H_2$ for $s \gg 0$
and $J_s$ a smooth family of admissible almost complex structures joining $J_1$ and $J_2$.
The set  $P(o_{-},o_{+})$ is the
set of solutions to the parameterized Floer equations
\[ \partial_s u + J_{s,t}(u(s,t)) \partial_t u = \nabla^{g_t} K_{s,t}\]
such that $u(s,\cdot)$ converges
to $x_{\pm}$ as $s \rightarrow \pm \infty$.
For a $C^{\infty}$ generic family $(K_s,J_s)$ this is a compact zero dimensional
manifold (if $o_-,o_+$ have the same relative index with respect to the cylinder
$C$ joining them).
Again the class $b \in H^2(N,\Z / 2\Z)$ enables us to orient this manifold.
If we have another such non-decreasing family admissible Hamiltonians joining $H_1$ and $H_2$
and another smooth family of admissible almost complex structures joining $J_1$ and $J_2$,
then the continuation map induced by this second family is the same as
the map induced by $(K_s,J_s)$.
The composition of two continuation maps is a continuation map.
If we take the direct limit of all these maps with respect
to admissible Hamiltonians $H$ ordered by $<$
such that $H|_N < 0$, then
we get our symplectic homology groups $SH_*^{(a,b]}(N,\eta)$.
We will write $SH_*^{\#}(N,\eta)$
or $SH_*^{\#}(H,J,\eta)$ for $SH_*^{[0,\infty)}(N,\eta)$
or $SH_*^{[0,\infty)}(H,J,\eta)$.

Also we will write $SH_*$ instead of $SH_*^{(-\infty,\infty)}$.
If we wish to stress which coefficient field we are using,
we will write $SH_*^{\#}(M,\eta,\K)$ or $SH_*^{\#}(H,J,\eta,\K)$ if the field is $\K$
for instance. We will also define $SH_*^{\leq b}(M,\eta,\K)$
to be the group $SH_*^{(-\infty,b]}(M,\eta,\K)$.

We will dealing with other pairs $(H,J)$ that are not necessarily
admissible. The definition of symplectic homology $SH_*^{(a,b]}(H,J,\eta)$ is still the same,
although $(H,J)$ has to satisfy some conditions to ensure
that we have a well defined symplectic homology group.
This will be discussed later in Section \ref{subsection:alternategrowthratedefinition}.
From now on instead of writing $SH_*^{(a,b]}(H,J,\eta)$ we will suppress
the term $\eta$ and just write $SH_*^{(a,b]}(H,J)$ instead when the context is clear.
Also from now on whenever we have a Liouville domain or symplectic manifold then
we will assume that we have chosen such a pair $\eta = (\tau,b)$.

\subsection{Growth rates} \label{section:growthrates}

In order to define growth rates, we will need some linear algebra first.
Let $(V_x)_{x \in [1,\infty)}$ be a family of vector spaces indexed by $[1,\infty)$.
For each $x_1 \leq x_2$ we will assume that there is a homomorphism
$\phi_{x_1,x_2}$ from $V_{x_1}$ to $V_{x_2}$ with the property that
for all $x_1 \leq x_2 \leq x_3$, $\phi_{x_2,x_3} \circ \phi_{x_1,x_2} = \phi_{x_1,x_3}$
and $\phi_{x_1,x_1} = \text{id}$.
We call such a family of vector spaces a {\it filtered directed system}.
Because these vector spaces form a directed system, we can take the direct
limit $V := \varinjlim_x V_x$.
From now on we will assume that each $V_x$ is finite dimensional.
For each $x \in [1,\infty)$ there is a natural map:
\[q_x : V_x \rightarrow \varinjlim_x V_x. \]
Let $a : [1,\infty) \rightarrow [0,\infty)$ be a function such that
$a(x)$ is the rank of the image of the above map $q_x$.
We define the growth rate as:
\[\Gamma( (V_x) ) : =\varlimsup_x \frac{\log{a(x)}}{\log{x}} \in \{-\infty\} \cup [0,\infty].\]
If $a(x)$ is $0$ then we just define $\log{a(x)}$ as $-\infty$.
If $a(x)$ was some polynomial of degree $n$ with positive leading coefficient, then the growth rate
would be equal to $n$.
If $a(x)$ was an exponential function with positive exponent, then the growth rate is $\infty$.
The good thing about growth rate is that
if we had some additional vector spaces $(V'_x)_{x \in [1,\infty)}$
such that the associated function $a'(x) := \text{rank}(V'_x \rightarrow \varinjlim_x V'_x)$
satisfies $a'(x) = A a(B x)$ for some constants $A,B>0$ then
$\Gamma(V'_x)  = \Gamma(V_x)$.
The notation we use for filtered directed systems
is usually of the form $(V_x)$ or $(V_*)$, and we will usually write
$V_x$ without brackets if we mean the vector space indexed by $x$.

In the previous section we defined for a Liouville domain $N$
(whose boundary had discrete period spectrum), $SH_*^{\leq \lambda}(N)$.
For $\lambda_1 \leq \lambda_2$, there is a natural map
$SH_*^{\leq \lambda_1}(N) \rightarrow SH_*^{\leq \lambda_2}(N)$
given by inclusion of the respective chain complexes.
This is a filtered directed system $(SH_*^{\leq \lambda}(N))$
whose direct limit is $SH_*(N)$.
\begin{defn} \label{defn:growthrate}
We define the growth rate $\Gamma(N,b)$ as:
\[\Gamma(N,b) := \Gamma(SH_*^{\leq \lambda}(N,b))\]
\end{defn}
We also have the following Theorem:
\begin{theorem} \label{thm:growthrateinvariance}
Let $N_1,N_2$ be two Liouville domains such that $\widehat{N_1}$
is symplectomorphic to $\widehat{N_2}$ where the symplectomorphism
pulls back $b_2 \in H^2(N_2,\Z / 2\Z)$ to $b_1 \in H^2(N_1,\Z / 2\Z)$
and $\tau_2$ to $\tau_1$ where $\tau_2$ and $\tau_1$ are trivializations of the canonical bundle.
Then
$\Gamma(N_1,(\tau_1,b_1)) = \Gamma(N_2,(\tau_2,b_2))$.
\end{theorem}
This theorem will be proven in Section \ref{subsection:alternategrowthratedefinition}.
Hence we will just write $\Gamma(\widehat{N},d\theta_N,(\tau,b))$ for the growth rate of $(N,\theta_N)$.
We will sometimes just write $\Gamma(\widehat{N})$ if the context makes it clear
that $d\theta_N$ is our symplectic form and $(\tau,b)$ is our associated
trivialization and homology class.

\section{Growth rate linear algebra}

Recall that a filtered directed system is a family of vector spaces $(V_x)$
parameterized by $[1,\infty)$ forming a category where for $x_1 \leq x_2$
there is a unique homomorphism from $V_{x_1}$ to $V_{x_2}$
and no other morphisms anywhere else.
For technical reasons we will define $V_x$ to be zero for $x < 1$ and
so all the morphisms starting with one of these vector spaces is also $0$.
A morphism of filtered directed systems $\phi : (V_x) \rightarrow (V'_x)$
consists of some constant $C_\phi$ and a sequence
of maps \[a_x : V_x \rightarrow V'_{C_\phi x}\]
so that we have the following commutative diagram:
\[
\xy
(0,0)*{}="A"; (40,0)*{}="B";
(0,-20)*{}="C"; (40,-20)*{}="D";
(0,-40)*{}="E"; (40,-40)*{}="F";
"A" *{V_{x_1}};
"B" *{V'_{C_\phi x_1}};
"C" *{V_{x_2}};
"D" *{V'_{C_\phi x_2}};
"E" *{V_{x_3}};
"F" *{V'_{C_\phi x_3}};
%
{\ar@{->} "A"+(12,0)*{};"B"-(15,0)*{}};
{\ar@{->} "C"+(12,0)*{};"D"-(15,0)*{}};
{\ar@{->} "E"+(12,0)*{};"F"-(15,0)*{}};
{\ar@{->} "A"+(0,-4)*{};"C"+(0,4)*{}};
{\ar@{->} "B"+(0,-4)*{};"D"+(0,4)*{}};
{\ar@{->} "C"+(0,-4)*{};"E"+(0,4)*{}};
{\ar@{->} "D"+(0,-4)*{};"F"+(0,4)*{}};
"A"+(20,3) *{a_{x_1}};
"C"+(20,3) *{a_{x_2}};
"E"+(20,3) *{a_{x_3}};
\endxy
\]
for all $x_1 \leq x_2 \leq x_3$
where the vertical arrows come from the directed system.

Let $\psi_{x_1,x_2}$ be the natural map from $V_{x_1}$ to $V_{x_2}$
in this filtered directed system for $x_1 \leq x_2$.
For each constant $C \geq 0$, we have an morphism $C_V$ from
$(V_x)$ to $(V_x)$ given by the map $\psi_{x,Cx}$.
We say that $(V_x)$ and $(V'_x)$ are {\it isomorphic}
if there is a morphism $\phi$ from $(V_x)$ to $(V'_x)$
and another morphism $\phi'$ from $(V'_x)$ to $(V_x)$
such that $\phi' \circ \phi = C_V$ and $\phi \circ \phi' = C'_{V'}$
where $C,C'\geq 0$ are constants and $C_V : (V_x) \rightarrow (V_x)$,
$C'_{V'} : (V'_x) \rightarrow (V'_x)$ are the morphisms described above.
One of the aims of this paper is to assign for each completion $\widehat{M}$
of a Liouville domain $M$, a filtered directed system which unique up to isomorphism.
From this we can define growth rate.
In order for growth rate to be well defined, we need to show:
\begin{lemma} \label{lemma:growthrateuptoisomorphism}
Let $(V_x),(V'_x)$ be two isomorphic filtered directed systems,
then $\Gamma(V_x) = \Gamma(V'_x)$.
\end{lemma}
\proof
Let $\phi_x : V_x \rightarrow V'_{C_\phi x}$
and $\phi'_x : V'_x \rightarrow V_{C_{\phi'} x}$
be our isomorphisms.
We have a morphism from $\varinjlim_x V_x$
to  $\varinjlim_x V'_x$ induced from $\phi$
and an inverse induced from $\phi'$.
This is because the morphism induced by $C_V$
is the identity map on $\varinjlim_x V_x$ and similarly
$C'_{V'}$ induces the identity map.
We will write $\phi$ and $\phi'$ for such maps by abuse of notation.
Let $a_x : V_x \rightarrow \varinjlim_x V_x$
and $a'_x : V'_x \rightarrow \varinjlim_x V'_x$.
Because $\phi$ is an isomorphism on $\varinjlim_x V_x$, we have that
the rank of the image of $\phi \circ a_x$ is equal to the rank
of the image of $a_x$.
We have that $a'_{C_\phi x} \circ \phi_x = \phi \circ a_x$
which implies that $\text{rank im}(a'_{C_\phi} x) \geq \text{rank im}(a_x)$.
Similarly $\text{rank im}(a_{C_{\phi'} x}) \geq \text{rank im}(a'_x)$
for all $x \in [1,\infty)$.
Hence:
\[\varlimsup_x \frac{\log{\text{rank im}(a_x)}}{\log{x}} \leq
\varlimsup_x \frac{\log{\text{rank im}(a'_{C_\phi x})}}{\log{x}}
= \varlimsup_x \frac{\log{\text{rank im}(a'_{C_\phi x})}}{\log{C_\phi x}} = \]
\[\varlimsup_x \frac{\log{\text{rank im}(a'_{x})}}{\log{x}}
\leq \varlimsup_x \frac{\log{\text{rank im}(a_{C_{\phi'} x})}}{\log{x}}
= \varlimsup_x \frac{\log{\text{rank im}(a_x)}}{\log{x}}.\]
This implies that
$\Gamma(V_x) = \Gamma(V'_x)$ as the first term in the above set of inequalities
is $\Gamma(V_x)$ and the fourth term is $\Gamma(V'_x)$.
\qed

We now need a Lemma giving us a sufficient condition
telling us when two filtered directed systems are equivalent.
\begin{lemma} \label{lemma:equivalentfiltereddirectedsystems}
Let $(V_x^j)$ ($j=1,2,3,4$) be filtered directed systems.
Let $u_j : V_x^j \rightarrow V_{C_j x}^{j+1}$ ($j=1,2,3$)
be morphisms of directed systems so that
composing any two of them gives us an isomorphism.
Then $V_x^2$ is isomorphic to $V_x^3$.
\end{lemma}
\proof of Lemma \ref{lemma:equivalentfiltereddirectedsystems}.
Let $\psi^j_{x_1,x_2}$ be the directed system
map $V^j_{x_1} \rightarrow V^j_{x_2}$ for $x_1 \leq x_2$.
By the definition of an isomorphism of filtered directed systems,
there exist morphisms
$b_1 : V_x^3 \rightarrow V_{D_1 x}^1$
and $b_2 : V_x^4 \rightarrow V_{D_2 x}^2$
so that
$u_2 \circ u_1 \circ b_1$
is the directed system map \[\psi^3_{x,D_1 C_1 C_2 x} : V_x^3 \rightarrow V_{D_1 C_1 C_2 x}^3.\]
Also
\[b_1 \circ u_2 \circ u_1 = \psi^1_{x, D_1 C_1 C_2 x},\]
\[u_3 \circ u_2 \circ b_2 = \psi^4_{x, D_2 C_2 C_3 x}\]
and
\[b_2 \circ u_3 \circ u_2 = \psi^2_{x, D_2 C_2 C_3 x}.\]

We define $\phi : V_x^2 \rightarrow V_{C_2 x}^3$
by $\phi = u_2$.
We define $\phi' : V_x^3 \rightarrow V_{D_1 D_2 C_1 C_2 C_3 x}^2$
by $\phi' = b_2 \circ u_3 \circ \psi^3_{x,D_1 C_1 C_2 x}$.

Let $x \in V^3_x$.
%
By abuse of notation for any $x$ we will just write
$\psi^j_K$ for $\psi^j_{x,Kx}$ for $j=1,2,3$.
Because
$\psi^3_{D_1 C_1 C_2} = u_2 \circ u_1 \circ b_1$,
\[\phi \circ \phi' = u_2 \circ b_2 \circ u_3 \circ u_2 \circ u_1 \circ b_1\]
\[ = u_2 \circ \psi^2_{D_2 C_2 C_3} \circ u_1 \circ b_1
= \psi^3_{D_2 C_2 C_3} \circ u_2 \circ u_1 \circ b_1\]
\[ =  \psi^3_{x, D_2 C_2 C_3 x} \circ \psi^3_{D_1 C_1 C_2}.\]
Also:
\[\phi' \circ \phi =  b_2 \circ u_3 \circ \psi^3_{D_1 C_1 C_2} \circ u_2 \]
\[= \psi^2_{D_1 C_1 C_2} \circ b_2 \circ u_3 \circ u_2\]
\[= \psi^2_{D_1 C_1 C_2} \circ \psi^2_{D_2 C_2 C_3}.\]

Hence $\phi$ and $\phi'$ give us our isomorphism and we have proven the Lemma.
\qed

\section{Growth rate geometry}

\subsection{Some alternate definitions of growth rate} \label{subsection:alternategrowthratedefinition}

We will define growth rate for a slightly larger class of manifolds
called finite type convex symplectic manifolds and also using a broader class of Hamiltonians.
There are three reasons for doing this.
The first reason is that we wish to prove that growth rate is an invariant
up to symplectomorphism and so we need a definition of growth rate
which is invariant under symplectomorphism.
The second reason is that the author wishes to use this larger class
of Hamiltonians in a future paper to prove that
growth rate behaves well under products and also with respect to Lefschetz fibrations.
A third reason is that this way of thinking might be useful for
answering various dynamical questions.

A {\it convex symplectic manifold} is a manifold $M$ with a $1$-form $\theta_M$
such that
\begin{enumerate}
\item $\omega_M := d\theta_M$ is a symplectic form.
\item There is an exhausting function $f_M : M \rightarrow \R$
and a sequence $c_1 < c_2 < \cdots$ tending to infinity such that
the $\omega_M$-dual $X_{\theta_M}$ of $\theta_M$ satisfies
$df_M(X_{\theta_M}) > 0$ along $f_M^{-1}(c_i)$ for each $i$.
\end{enumerate}
Some basic facts about convex symplectic manifolds are proven in the appendix.
We say that $M$ is of {\it finite type} if there is a $C \in \R$
such that $(f_M^{-1}(-\infty,c],\theta_M)$ is a Liouville domain
for all $c \geq C$.

Let $(M,\theta^t_M)$ be a smooth family of convex symplectic manifolds
parameterized by $t \in [0,1]$.
This is said to be a {\it convex deformation} if
for every $t \in [0,1]$ there is a $\delta_t>0$
and an exhausting function $f_M^t$ and a sequence
of constants $c_1^t < c_2^t < \cdots$ tending to infinity
such that $( (f_M^t)^{-1}(-\infty,c_i^t],\theta^s_M)$ is a Liouville domain
for each $s \in [t- \delta_t,t+\delta_t]$ and each $i$.
We do not require that $f_M^t$,$c_i^t$,$\delta_t$
smoothly varies with $t$. In fact it can vary in
a discontinuous way with $t$.

Let $M$ be a finite type convex symplectic manifold.
In order to define growth rate, we need a slightly larger class of Hamiltonians.
We will first describe Hamiltonians on $M$ that give us filtered directed systems.
Let $(S,j)$ be a complex surface possibly with boundary.
Let $\gamma$ be a $1$-form on $S$  so that $d\gamma \geq 0$.
Let $H : M \rightarrow \R$ be a Hamiltonian and $J$ an almost
complex structure compatible with the symplectic form $\omega$.
We have a Hamiltonian vector field $X_{H}$ defined by $\omega(X_H,\cdot) = dH$.
A map $u : S \rightarrow M$
satisfies the {\it perturbed Cauchy-Riemann equations} if
$(du - X_H \otimes \gamma)^{0,1} = 0$.
Here
$du - X_H \otimes \gamma$ is a $1$-form on $S$ with values in
the complex vector bundle $\text{Hom}(TS,u^* TM)$ where the complex
structure at a point $s \in S$ is induced from $j$ and $J$.
The equation $(du - X_H \otimes \gamma)^{0,1} = 0$ is written explicitly
as 
\begin{equation} \label{equation:cauchyriemann}
du - X_H \otimes \gamma + J \circ (du - X_H \otimes \gamma) \circ j = 0.
\end{equation}
Sometimes we will write
$(du - X_H \otimes \gamma)^{0,1}_J = 0$ if we wish to emphasise
the fact that we are using the almost complex structure $J$.
Here is a particular example.
Let $S = \R \times S^1 = \C / \Z$. We let $\gamma = dt$
where $t$ parameterizes $S^1 = \R / \Z$.
Then the perturbed Cauchy-Riemann equations become
\[\partial_s u + J \partial_t u = J X_H\]
which is just the Floer equation.

The pair $(H,J)$ on $M$
is said to satisfy a {\it maximum principle with respect to an open set} $U^H$
if there is a compact set $K' \subset M$ containing $U^H$
such that for any compact complex surface $(S,j)$ with $1$-form $\gamma$ 
($d\gamma \geq 0$) and map $u : S \rightarrow M$,
satisfying:
\begin{enumerate}
\item $u$ satisfies the perturbed Cauchy-Riemann equations.
\item $u(\partial S) \subset U^H$
\end{enumerate}
we have that $u(S) \subset K'$.
We also require that $U^H$ contains all the $1$-periodic
orbits of $(H,J)$ of action greater than some small negative constant.

A pair $(H,J)$ is said to be $SH_*$ admissible if
there is a discrete subset $A_H \subset (0,\infty)$
such that $(\lambda H,J)$ satisfies the
maximum principle for $\lambda \in (0,\infty) \setminus A_H$
with respect to some relatively compact open set $U^H_\lambda$.
We require that $U^H_{\lambda_1} \subset U^H_{\lambda_2}$ for $\lambda_1 \leq \lambda_2$.
Note that if $(H,J)$ is $SH_*$ admissible then so is
any positive multiple of $H$.

For an $SH_*$ admissible pair $(H,J)$, we define $SH_*^{\#}(\lambda H,J)$
($\lambda \in (0,\infty) \setminus A_H$) as follows:
By Lemma \ref{lemma:perturbinghamiltonian},
we can perturb $\lambda H$ by an arbitrarily small amount
to a non-degenerate time dependent Hamiltonian $H'_t : S^1 \times M \rightarrow \R$
so that it is equal to $\lambda H$ outside some closed subset of $U^H_\lambda$.
After subtracting a small constant from $H'_t$,
we can assume that $H'_t$ is equal to $\lambda H - \epsilon$
on a closed subset of $U^H_\lambda$ and $H'_t < \lambda H$ where $\epsilon>0$ is some constant.
For a generic $S^1$ family $J'_t$ of almost complex structures
such that $J_t = J$ outside a closed subset of $U^H_\lambda$,
we have that $SH_*^{\#}(H'_t,J'_t)$ is well defined.
This is because all the $1$-periodic orbits of $H'_t$
are contained in a compact set and the maximum principle
ensures that all the Floer trajectories also stay inside a compact set.
The pair $(H_t,J_t)$ constructed above is called an {\it approximating pair}
for $(H,J)$.
A similar argument ensures that 
if we have two approximating pairs $(H'_t,J'_t)$ and $(H''_t,J''_t)$ for
$(H,J)$
with $H'_t < H''_t$,
there is a well defined
continuation map $SH_*^{\#}(H'_t,J'_t) \rightarrow SH_*^{\#}(H''_t,J''_t)$
induced from a generic family
of approximating pairs joining $(H'_t,J'_t)$ and $(H''_t,J''_t)$
where the Hamiltonians are non-decreasing.
Also because the continuation map is induced from a non-decreasing 
family of Hamiltonians, we have that the continuation map
respects the filtration by action. So any orbit of action $\leq f$
is sent to another orbit of action $\leq f$ under the continuation map.
The set of approximating pairs $(H'_t,J'_t)$ induces a directed system
where $(H'_t,J'_t) \leq (H''_t,J''_t)$ if and only if $H'_t \leq H''_t$ for all $t \in S^1$.
We define $SH_*^{\#}(\lambda H,J)$ as $\varinjlim_{(H'_t,J'_t)} SH_*^{\#} (H'_t,J'_t)$
where the direct limit is taken over this directed system.

We also have a continuation map for $\lambda_1 \leq \lambda_2$
where $\lambda_i \in (0,\infty) \setminus A_H$
from $SH_*^{\#}(\lambda_1 H,J)$ to $SH_*^{\#}(\lambda_2 H,J)$.
This continuation map is induced by a family of pairs
$(H^s_t,J^s_t)$ such that
\begin{enumerate}
\item  they are equal to $(\lambda'_s H - \epsilon_s,J)$ outside a closed subset of
$U^H_{\lambda_2}$ for some
smooth non-decreasing family of constants $\lambda'_s$
and some smooth non-increasing family of constants $\epsilon_s$.
\item $(H^s_t,J^s_t)$ is a fixed approximating pair $(H^-_t,J^-_t)$ for $(\lambda_1 H,J)$
for $s \ll 0$.
\item $(H^s_t,J^s_t)$ is a fixed approximating pair $(H^+_t,J^+,t)$ for $(\lambda_2 H,J)$
for $s \gg 0$.
\item  $\frac{\partial}{\partial s} H^s_t \geq 0$.
\end{enumerate}
We need to show that the continuation map trajectories stay inside a compact set.
Suppose we have some map $u : \R \times S^1 \rightarrow M$
satisfying the continuation map equations which join orbits inside $U^H_{\lambda_2}$:
\[ \partial_s u + J^s_t \partial_t u = J^s_t X_{H^s_t}.\]
These can be rewritten in the following way:
\[ (du - X_{H^s_t} \otimes dt)^{0,1}_{J^s_t} = 0.\]
We have $(H^s_t,J^s_t) = (\lambda'_s H - \epsilon_s,J)$
outside some compact subset $R$ of $U^H_{\lambda_2}$.
Let $S \subset \R \times S^1$ be a compact submanifold so that
$u(S)$ is disjoint from $R$ and $u(\partial S) \subset U^H_{\lambda_2}$.
Hence we have that $u$ restricted to $S$ satisfies:
\[ (du - X_{\lambda'_s H - \epsilon_s} \otimes dt)^{0,1}_{J} = 0.\]
The constants $\epsilon_s$ do not matter so our equation becomes:
\[ (du - X_{H} \otimes \gamma)^{0,1}_{J} = 0\]
where $\gamma = \lambda'_s dt$. Here $d\gamma \geq 0$ because
$\lambda'_s$ is non-decreasing.
Hence by the maximum principle we have that $S$ must be contained entirely
inside $U^H_{\lambda_2}$.
This implies that the image of each continuation map $u$ must be contained inside $U^H_{\lambda_2}$
and hence inside a fixed compact set.
Hence we have a well defined continuation map
between $SH_*^{\#} (H^-_t,J^-_t)$ 
and $SH_*^{\#} (H^+_t,J^+_t)$.
This induces a continuation map between
$SH_*^{\#} (\lambda_1 H,J)$ and
$SH_*^{\#} (\lambda_2 H,J)$.
If $\lambda \in A_H$ then we define $SH_*^{\#}(\lambda H, J)$
as the direct limit of $SH_*^{\#}(\lambda' H,J)$ where
$\lambda' \notin A$ tends to $\lambda$ from below.
Hence we have a filtered directed system
$(SH_*^{\#}(\lambda H_t,J_t))$.

We wish to put some additional conditions on the pair
$(H,J)$ so that the associated filtered directed
system is an invariant related to symplectic homology.
Here are the additional conditions:

\begin{enumerate}
\item (bounded below property)

The Hamiltonian $H$ is greater than or equal to zero,
and there exists a compact set $K$ and a constant $\delta_H > 0$
such that:
$H > \delta_H$ outside $K$.

\item (Liouville vector field property)

There exists an exhausting function $f_H$, and $1$-form $\theta_H$
such that:
\begin{enumerate}
\item $\theta_M - \theta_H$ is exact where $\theta_M$ is the
Liouville form on $M$.
\item There exists a small $\epsilon_H > 0$ such that
$dH(X_{\theta_H}) > 0$ in the region $H^{-1}(0,\epsilon_H]$
where $X_{\theta_H}$ is the $\omega_M$-dual of $\theta_H$.
We also require that $d\left(dH(X_{\theta_H})\right)(X_{\theta_H}) > dH(X_{\theta_H})$ in this region.
\item
There is a constant $C$ such that
$df_H(X_{\theta_H}) > 0$
in the region $f_H^{-1}[C,\infty)$ and
$f_H^{-1}(-\infty,C]$ is non-empty and is contained in the interior of $H^{-1}(0)$.
\end{enumerate}

\item (action bound property)

There is a constant $C_H$ and $1$-form $\theta$ such that
the function
$-\theta(X_H) - H$ must be bounded above by $C_H$
where $X_H$ is the Hamiltonian vector field associated to $H_t$.
We also require that $\theta - \theta_M$ is exact.

\end{enumerate}
If an $SH_*$ admissible pair $(H,J)$ has these additional conditions,
then it is called {\it growth rate admissible}.

Here are some important examples:

{\bf First example :}
Let $X_{\theta_M}$ be the vector field given by the $\omega_M$-dual of $\theta_M$.
Because $M$ is a finite type convex symplectic manifold,
we have a function $f_M : M \rightarrow \R$ such that,
$df_M(X_{\theta_M}) > 0$ in the region where $f_M > c$ for some $c \gg 0$ and such that $f_M$ is exhausting.
Let $X_{f_M}$ be the Hamiltonian vector field of $f_M$.
Choose some $c < c_1 < c_2 < \cdots $ tending to infinity
so that $df_M(X_{\theta_M})>0$ on $f_M^{-1}(c_i)$ for each $i$.
We perturb $f_M$ by a $C^0$ small amount to $g_M$
so that it has the following property:
$g_M^{-1}(c_i) = f_M^{-1}(c_i)$
and $dg_M(X_{\theta_M})$ is constant on a small neighbourhood
of $g_M^{-1}(c_i)$ for each $i$ and $dg_M(X_{\theta_M})> 0$ for
$g_M > c$.
The level sets $A_y := g_M^{-1}(y)$ for $y > c$ are contact manifolds
with contact form $\alpha_y := \theta_M|_{A_y}$.
Because
\[\alpha_y(X_{g_M}) = \theta_M(X_{g_M}) =
\omega_M(X_{\theta_M},X_{g_M}) = -dg_M(X_{\theta_M}) < 0,\]
we have that $X_{g_M}$ is non-trivial in the region $\{g_M > c\}$.
Because the vector field $X_{g_M}$ is contained in $TA_y$ ($y>c$),
and because it is non-zero, it has a shortest orbit
inside $A_y$ (as $A_y$ is compact).
The reason for this is that for each $p \in A_y$
we can choose coordinates $x_1,\cdots,x_{2n}$ around $p$ so that
$X_{g_M} = \frac{\partial}{\partial x_1}$, which means that
any flowline going through a smaller neighbourhood $U$ of
$p$ must take some time $\delta_U >0$ to pass through $U$.
Because $A_y$ is compact, we can then cover it with finitely many such neighborhoods
$U_i$ which implies that any orbit must flow for time at least
$\text{min}_i \delta_{U_i}$.
Let $\delta(y)>0$ be smaller than $\text{min}_i \delta_{U_i}$.
We can assume that $\delta(y)$ is a smooth function of $y$.

Let $h : \R \rightarrow \R$ be a function with the property
that $h(x)=0$ for $x \leq c+1$ and $h(x)>0$ elsewhere.
For $x \geq c+1$ we let $h'(x)$ be smaller than
$\delta(x) / g_M$ and $h'>0$ whenever $h>0$.
We also assume that $h'$ is small enough so that
$h'(x) \alpha(x)(X_{g_M})$ is bounded above by a constant.
The Hamiltonian flow of $H := h(g_M)$ in the region
$g_M \leq c$ is $0$ and for $x > c+1$
in the region $A_x$, it is equal to some very small multiple of $X_{g_M}|_{A_x}$.
We also need an additional technical condition which ensures that
our Hamiltonian $H$ will satisfy the Liouville vector field property. 
We require that there exists some $\epsilon_h > 0$ with $h''(x) > 0$
for $x \in (c+1,c+1+\epsilon_h)$.
Fix $\lambda > c$. In the region $\{g_M > \lambda\}$,
there are no $1$-periodic orbits of $\lambda H$.
This is because all the orbits of $X_H|_{A_x} = h'(x) X_{g_M}|_{A_x}$ have
length greater than $\lambda$ because $h' < \delta_x / g_M \leq \delta_x / \lambda$.
Also because $H \geq 0$ and $\alpha_x(X_H)$ is bounded from
above for $x > c$, we have that $-\theta_M(X_H) - H$ is bounded from
above so $H$ satisfies the action bound property.
The Hamiltonian $H$ satisfies the bounded below property because
its zero set is compact, it is greater than or equal to zero
and is greater than some constant near infinity.
Note that $\log(h'(x))$ tends to minus infinity
as $x$ tends to $c+1$ from above. Also the derivative
of this function is positive.
This implies that the derivative $h''(x)/h'(x)$
tends to infinity. Hence (after shrinking $\epsilon_h$)
we can assume that
$h''(x) \geq \nu h'(x)$ for $x \in (c+1,c+1+\epsilon_h)$ where
\[ \nu > \left(1-d\left(g_M(X_{\theta_M})\right)(X_{\theta_M})\right) / g_M(X_{\theta_M})^2 \]
in the region $x \in (c+1,c+1+\epsilon_h)$.
Because $dH(X_{\theta_M}) = h'(g_M) dg_M(X_{\theta_M}) \geq 0$
and
\[ d\left(dH(X_{\theta_M})\right)(X_{\theta_M}) =
h''(g_M)dg_M(X_{\theta_M})^2 + h'(g_M) d\left(g_M(X_{\theta_M})\right)(X_{\theta_M}) > \]
\[ h'(g_M) dg_M(X_{\theta_M})\]
in the region $H^{-1}(0,\epsilon_H)$ for some small $\epsilon_H>0$,
$H$ satisfies the Liouville vector field property.
We need to find a $J$ so that
$(H,J)$ satisfies a maximum principle.
Because $dg_M(X_{\theta_M})$ is constant on a small neighbourhood of $c_i$
for each $i$, we have (by flowing $G_i$ along $X_{\theta_M}$) a small neighbourhood of
$G_i := g_M^{-1}(c_i)$ symplectomorphic to $(1-\epsilon_i,1+\epsilon_i) \times G_i$
with contact form $\kappa_i g_M \alpha_i$ where
$\alpha_i := \theta_M|_{G_i}$ is a contact form and $\kappa_i$ is a constant.
This is a slice of a positive cylindrical end so we can choose $J$ so that it looks cylindrical
in these regions.
We define $U^H_\lambda$ to be any family open sets of the form
$g_M^{-1}(-\infty,c_i)$ containing all the $1$-periodic
orbits of $\lambda H$ and such that $U^H_{\lambda_1} \subset U^H_{\lambda_2}$
for $\lambda_1 \leq \lambda_2$.
Hence by Lemma \ref{lemma:maximumprinciple},
there exists a $J$ such that
$(\lambda H,J)$ satisfies the maximum principle with respect to $U^H_\lambda$.
Hence $(H,J)$ is growth rate admissible.

{\bf Second example :}
Let $f_M$ and $X_{\theta_M}$ be as above.
Let $Q := f_M^{-1}(-\infty,C]$ be a manifold
with boundary such that $X_{\theta_M}$ is transverse to the
boundary of $Q$ and pointing outwards
and $df_M(X_{\theta_M}) > 0$ outside $Q$.
Even though the Liouville vector field $X_{\theta_M}$ may not be complete,
we can still flow $\partial Q$ along $X_{\theta_M}$ so that we have that
the set $\overline{(M \setminus Q)}$ is diffeomorphic
to some open subset $U$ of $[1,\infty) \times \partial Q$ containing $\{1\}
\times \partial Q$.
We have that $\theta_M|_U = r\alpha$ where $r$ parameterizes $[1,\infty)$
and $\alpha = \theta_M|_{\partial Q}$.
We call this a {\it partial cylindrical end} of $M$.
We can ensure that the period spectrum of the contact manifold $\partial Q$
is discrete and injective (after perturbing $f_M$ very slightly).
Let $J$ be an almost complex structure which is cylindrical inside $U$
(i.e. $dr \circ J = -\theta_M$).
Let $\epsilon>0$ be a constant small enough so that
$r^{-1}(1+\epsilon)$ is still a compact manifold transverse to $X_{\theta_M}$.
Let $H$ be a Hamiltonian such that $H = h(r)$ in $U$ and $0$ inside $Q$.
We require that $h(r)=0$ near $r=1$ and $h(r) = r-\frac{\epsilon}{4}$ for $r \geq 1+\epsilon/2$,
and $h',h'' \geq 0$ and $h' > 0$ when $h>0$.
We also require that $h''(x)>0$ for $x \in (1,1+\epsilon_H)$ where
$\epsilon_H>0$ is small.

\begin{figure}[H]
\centerline{
 \scalebox{1.0}{
 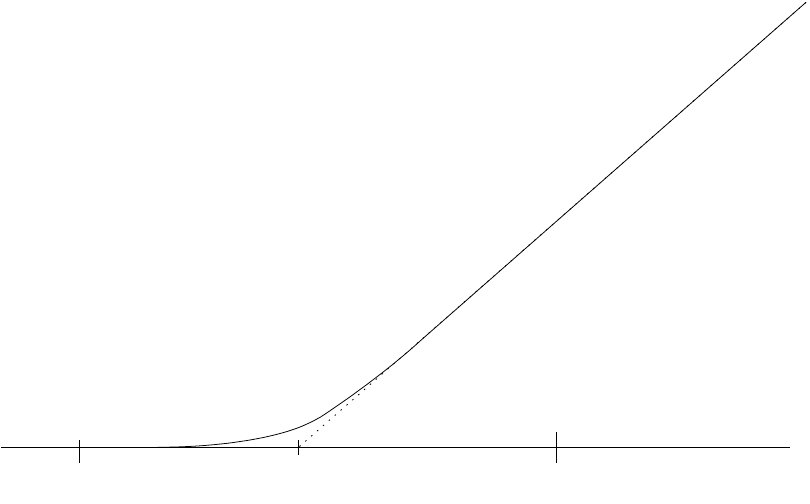
}
     }
\end{figure}

For $\lambda$ not in the period spectrum of $\partial Q$,
we have that all of the $1$-periodic orbits of $\lambda H$ lie
inside the compact set $\{r \leq \epsilon/2\}$.
We define $U^H_\lambda$ to be the open set
$Q \cup r^{-1}([1,1+\epsilon)$ for all $\lambda$.
The pair $(\lambda H,J)$ satisfies the maximum principle
with respect to $U^H_\lambda$
by \ref{lemma:maximumprinciple}.
Hence $(H,J)$ is $SH_*$ admissible.
We have that $(H,J)$ satisfies the bounded below property
and the Liouville vector field property because
$X_{\theta_M} = r \frac{\partial}{\partial r}$ in $U$, $h'>0$ when $h>0$
and $h''(r) > 0$ for $r \in (1,1+\epsilon_H)$.
We also have the that $(H,J)$ satisfies the
action bound property because
$-\theta_M(X_H) - H = rh'(r) - r + \frac{\epsilon}{4}= \frac{\epsilon}{4}$
outside a compact set (where $H = r - \frac{\epsilon}{4}$).
All of this means that $(H,J)$ is growth rate admissible.

Suppose that $(H_0,J_0)$, $(H_1,J_1)$ is growth rate admissible
with the property that $U^{H_0}_\lambda \subset U^{H_1}_\lambda$ for all $\lambda$
and such that $(H_1,J_1) = (\kappa H_0 + c_1,J_0)$ outside a closed subset of
of $U^{H_0}_\lambda$.
Suppose we have a non-decreasing family of pairs $(H_t,J_t)$
equal to $(\kappa_t H_0 + c_t,J)$ outside a closed subset of $U^{H_0}_\lambda$
for some smooth family of constants $\kappa_t \geq 1, c_t \in \R$.
Then we have a well defined continuation map from
$SH_*^{\#}(\lambda H_0,J_0)$ to 
$SH_*^{\#}(\lambda H_1,J_1)$ and this induces a morphism of filtered directed
systems.
We call such a morphism a {\it restricted continuation morphism}.
Suppose now we have two growth rate admissible pairs
$(H'_0,J'_0)$ and $(H'_1,J'_1)$ with the following properties:
\begin{enumerate}
\item $U^{H'_0}_\lambda = U^{H'_1}_\lambda$.
\item $(H'_0,J'_0) = (H'_1,J'_1)$ on a neighbourhood
of $\cup_{\lambda} U^{H'_0}_\lambda$.
\end{enumerate}
Then the filtered directed systems $(SH_*^{\#}(\lambda H'_0,J'_0))$
and $(SH_*^{\#}(\lambda H'_1,J'_1))$ are isomorphic
because all $1$-periodic orbits of non-negative action
and all Floer trajectories connecting them are identical.
We call such a morphism a {\it switch isomorphism}.

The main theorem of this section is:
\begin{theorem} \label{theorem:generalgrowthrateequivalence}
Let $(H,J),(H',J')$ be growth rate admissible.
Then the filtered directed system $(SH_*^{\#}(\lambda H,J))$ is isomorphic to
$(SH_*^{\#}(\lambda H',J'))$ as filtered directed systems.

The isomorphism between these filtered directed systems
is a composition of restricted continuation morphisms
and switch morphisms and inverses of these morphisms.
\end{theorem}

The reason why we mention restricted continuation morphisms
and switch morphisms is because we would
like this result to hold if $SH_*^{\#}$ had some additional algebraic structure
as well such as the pair of pants product so that it can be used
in future work.
The point is that if restricted continuation morphisms
and switch morphisms are also morphisms preserving this additional
algebraic structure then we immediately get an invariance result
for symplectic homology with this extra structure.

We need some preliminary lemmas before we prove this Theorem.
\begin{lemma} \label{lemma:positveactionhamiltonian}
Let $H$ be a Hamiltonian satisfying the Liouville vector field property, then
there is a growth rate admissible pair $(H^p,J^p)$
such that:
\begin{enumerate}
\item $H^p = H$ on a small neighbourhood of $H^{-1}(0)$
and $(H^p)^{-1}(0) = H^{-1}(0)$.
\item $-\theta_H(X_{H^p}) - H^p \geq 0$ everywhere and $-\theta_H(X_{H^p}) - H^p > 0$ when $H^p > 0$.
\item $-\theta_H(X_{H^p}) - H^p$ is greater than some constant $\delta'_{H^p}>0$
outside a large compact set.
\item The construction of $H^p$ only depends on $H$ near $H^{-1}(0)$.
\item There is a fixed compact set $K'$ such that all $1$-periodic
orbits of $\lambda H^p$ are contained in the interior of $K'$.
This set is an embedded codimension $0$ manifold with boundary
and we can ensure that it fits inside any open set containing $H^{-1}(0)$.
\item \label{item:strongermaximumprinciple}
Any solution $u: S \rightarrow M$ of the perturbed Cauchy-Riemann equations with
respect to $(\tilde{H},\tilde{J})$ where $\tilde{H}>0$,
$(\tilde{H},\tilde{J}) = (\lambda H^p,J)$ inside $K'$ and $u(\partial S) \subset K'$
is contained in $K'$.
\end{enumerate}
Here $\theta_H$ is the $1$-form that makes
$H$ satisfy the Liouville vector field property.
\end{lemma}
Note that this Lemma also tells us the following fact about $H$:
we have that $-\theta_H(H) - H > 0$ on some small neighbourhood
of $H^{-1}(0)$.
This will be useful later on.

\proof of Lemma \ref{lemma:positveactionhamiltonian}.
The main idea of this proof is to use bump functions to extend $(H,J)$
restricted to a small neighbourhood of $H^{-1}(0)$ to a pair which looks
like the pair described in the second example.
Here $\theta_H$ is a $1$-form such that $\theta_H - \theta_M$ is an exact $1$-form.
The vector field $X_{\theta_H}$ which is the $\omega_M$-dual of $\theta_H$
has the property that there is a function $f_H$ such that
$df_H(X_{\theta_H}) > 0$ inside a small neighbourhood of
the $\overline{M \setminus H^{-1}(0)}$.
We also have an $\epsilon_H$ such that $X_{\theta_H}(H) >0$
on $H^{-1}(0,\epsilon_H)$
by the Liouville vector field property.

Because we are only interested in what $H$ is near $H^{-1}(0)$
it can be anything we like outside a neighbourhood of this set. So from now on
(after changing $H$ outside a neighbourhood of $H^{-1}(0)$) we can assume that $H$ satisfies
the bounded below property.
Let $Q := f_H^{-1}(-\infty,C]$ be a manifold
with boundary such that $X_{\theta_H}$ is transverse to the
boundary of $Q$ and pointing outwards
and $df_H(X_{\theta_H}) > 0$ outside $Q$
and such that $H^{-1}(0)$ contains $Q$.
We have that $\partial Q$ is a contact manifold with contact form
$\alpha = \theta_H|_{\partial Q}$.
We flow $\partial Q$ along $X_{\theta_H}$ so that we have that
the set $\overline{(M \setminus Q)}$ is diffeomorphic
to some open subset $U$ of $[1,\infty) \times \partial Q$.
This is our partial cylindrical end.
We have that $\theta_H|_U = r\alpha$ where $r$ parameterizes $[1,\infty)$,
and $X_{\theta_H} = r \frac{\partial}{\partial r}$ in the region $U$.
We have:
\[-\theta_H(X_{H}) - H = r \frac{\partial H}{\partial r} - H\]
inside $U$.
Because $H$ satisfies the Liouville vector field property we have on a neigbourhood of $H^{-1}(0)$,
$d\left( dH(X_H) \right)(X_H) > dH(X_H) > 0$ outside $H^{-1}(0)$.
Because $X_H = r \frac{\partial}{\partial r}$ on our partial cylindrical end this condition becomes:
\[ r \frac{\partial H}{\partial r} + r^2 \frac{\partial^2 H}{\partial r^2} > r \frac{\partial H}{\partial r} > 0
\]
which in turn is equivalent to $\frac{\partial H}{\partial r}, \frac{\partial^2 H}{\partial r^2} > 0$.
We shrink $\epsilon_H$ so that $H^{-1}[0,\epsilon_H]$ is contained in $V$.
We also have that
\[r\frac{\partial^2 H}{\partial r} = \frac{\partial}{\partial r} (-\theta_H(X_H) - H).\]
Hence we have that
$-\theta_H(X_H) - H > 0$ in the region
$H^{-1}(0,\epsilon_H)$.
By the bounded below property, there exists a constant
$\delta_H>0$ such that $H>\delta_H$ outside $V$.
We shrink $\epsilon_H$ so that it is smaller than $\delta_H$.
Choose a function $q_1 : \R \rightarrow \R$
such that $q_1(x) = 1$ for $x \leq 2\epsilon_H / 3$
and $q_1(x) = 0$ for $x \geq \epsilon_H$.
Because $dH(X_{\theta_H}) > 0$ on the level set $W := H^{-1}(\epsilon_H/2)$,
we have that this level set is regular 
and is a contact manifold with contact form $\alpha_2 := \theta_H|_{W}$.
There exists a function $f : \partial Q \rightarrow [1,\infty)$
such that under the identification of $\overline{M \setminus Q}$
with a subset of $[1,\infty) \times \partial Q$,
\[W = \{(f(x),x) | x \in \partial Q\} \subset [1,\infty) \times \partial Q.\]
We will write this set as $(r/f)^{-1}(1)$.
We also have a new partial cylindrical end
which is the region $\{(r/f) \geq 1\}$.
This is diffeomorphic to some codimension $0$ submanifold with boundary of
$W \times [1,\infty)$ where the cylindrical
coordinate is $(r/f)$ and the contact form
is $\theta_H|_W$.
Let $q_2$ be a function on $\R$ satisfying:
\begin{enumerate}
\item 
$q_2(x) = 0$ for $x \leq 0$
\item
$q_2'(x),q_2''(x) \geq 0$.
\item
We let $q_2'(x)$ be constant and equal to $1$
for $x \geq \delta_r$.
Here $\delta_r$ is a constant such that
the level sets $(r/f)^{-1}(x)$ are all compact
for $1 \leq x \leq \delta_r + 1$.
So for $ x \geq \delta_r$, $q_2 = r - \iota_r$
where $\iota_r>0$ is a constant.
\item
For $x > 0$ we also choose $q_2$ so that
$q_2'(x)>0$.
\end{enumerate}
Let $H^p := q_1(H) H + \kappa q_2(r/f - 1)$ where $\kappa>0$  is a constant to be determined.
Here $q_1(H)$ has support inside the region $V$ because
$H \geq \epsilon_H$ outside $V$.
We have for $\kappa \gg 0$, that
$-\theta_H(H^p)(X_{\theta_H})- H^p \geq 0$.
This is because
$H^p = H + \kappa q_2(r/f - 1)$ in $H^{-1}(0,2\epsilon_H/3)$
and $-\theta_H(H) - H > 0$, $-\theta_H(q_2(r/f - 1)) - q_2(r/f - 1) \geq 0$
in this region because $q_2,q_2',q_2'' \geq 0$.
Also $-\theta_H(r/f - 1 - \iota_r) - r/f + \iota_r = \iota_r$ is greater than some
fixed constant outside $H^{-1}(-\infty,2\epsilon_H/3)$
hence by the action bound property we get that
$-\theta_H(X_{H^p})- H^p \geq 0$
for $\kappa$ large enough.
Hence in the region where $H^p > 0$, we have that
$-\theta_H(X_{H^p}) - H^p > 0$.

We also let  $J^p$ be an almost complex structure
such that it is cylindrical on the partial cylindrical end $\{(r/f) \geq 1\}$
viewed as a subset of $W \times [1,\infty)$.
We can perturb $f$ very slightly so that
the period spectrum of $W = (r/f)^{-1}(1)$ is discrete.
Define $K' := M \setminus (r/f)^{-1}(1,\infty)$.
This is a compact co-dimension $0$ manifold whose boundary is $W$.
Because the period spectrum of $W$ is discrete
and equal to $A \subset (0,\infty)$,
we have for all $\lambda \in (0,\infty) \setminus A$, $\lambda H^p$ has all
its $1$-periodic orbits contained in a compact subset of $M$.
Also \ref{lemma:maximumprinciple}
with boundary inside $K'$ must stay in $K'$.
This implies that $(H^p,J^p)$ satisfies the maximum principle.
We also have that 
that because $H^p \geq 0$ and $H^p= (r/f) - \iota_r$ near infinity,
$H^p$ satisfies the bounded below property.
Because $H$ satisfies the Liouville vector field
property, then so does $H^p$ as $H^p= H$
on a small neighbourhood of $H^{-1}(0)$
and ${H^p}^{-1}(0) = H^{-1}(0)$.
We also have that it satisfies the action bound property
because $-\theta_H(r/f - \iota_r) - (r/f - \iota_r) = \iota_r$
near infinity.

Finally property (\ref{item:strongermaximumprinciple})
is satisfied by Lemma \ref{lemma:maximumprinciple} as well.
\qed

Before we prove Theorem \ref{theorem:generalgrowthrateequivalence}
we will prove a couple of weaker versions of it. Here is one.
\begin{lemma} \label{lemma:veryrestrictedgeneralgrowthrateequivalence}
Let $(H,J),(H',J')$ be growth rate admissible pairs satisfying the following properties:
\begin{enumerate}
\item $H' = H$
in a neighbourhood of $H^{-1}(0)$.
\item $H^{-1}(0) = {H'}^{-1}(0)$.
\item $U^H_\lambda = U^{H'}_\lambda$.
\item $(H,J) = (H' + c,J')$
outside a closed subset of $U^H_1 = \cap_\lambda U^H_\lambda$.
\end{enumerate}
Then there is a restricted continuation morphism from
$(SH_*^{\#}(\lambda H,J))$ to
$(SH_*^{\#}(\lambda H',J'))$ which is an isomorphism.
Its inverse is also a restricted continuation morphism.
\end{lemma}
\proof of Lemma \ref{lemma:veryrestrictedgeneralgrowthrateequivalence}.
By the bounded below property
there exists a constant $\delta > 0$ with $H,H' > \delta$
outside a compact set $K$.
This implies that there exists constants $Q_1,Q_2 > 1$
such that $H \leq Q_1 H'$ and $H' \leq Q_2 H$.
Let $(H^s,J^s)$ $(s \in \R)$ be a family of pairs such that
\begin{enumerate}
\item $(H^s,J^s) = (H,J)$ for $s \ll 0$ and
$(H^s,J^s) = (Q_1 H',J')$ for $s \gg 0$.
\item We let $(H^s,J^s) = (\kappa_s H + c_s,J)$ outside
a closed subset of $\cap U^H_\lambda$
for some non-decreasing family of constants  $\kappa_s \geq 1$, $c_s \in \R$.
\item $\frac{\partial H^s}{\partial s} \geq 0$.
\end{enumerate}
The family of pairs $(\lambda H^s, J^s)$
induces a continuation map $\phi$
from the group $SH_*^{\#}(\lambda H,J)$ to $SH_*^{\#}(\lambda Q_1 H',J')$.
A similar family of pairs induces a continuation map
$\phi'$ from 
$SH_*^{\#}(\lambda H',J')$ to $SH_*^{\#}(\lambda Q_2 H',J')$.
We also have that $\phi' \circ \phi$ and $\phi \circ \phi'$
are directed system maps for
$(SH_*^{\#}(\lambda H',J'))$ and $(SH_*^{\#}(\lambda H',J'))$ respectively.
This is because the directed system maps are also continuation maps
induced by increasing families of Hamiltonians.
Hence $(SH_*^{\#}(\lambda H,J))$ and $(SH_*^{\#}(\lambda H',J'))$
are isomorphic as directed systems.
Also $\phi$ and its inverse $\phi'$ are restricted continuation morphisms.
\qed

\begin{lemma} \label{lemma:restrictivegeneralgrowthrateequivaalence}
Let $(H,J),(H',J')$ be growth rate admissible pairs such that
$H' = H$
in a neighbourhood of $H^{-1}(0)$
and $H^{-1}(0) = {H'}^{-1}(0)$.
Then $(SH_*^{\#}(\lambda H,J))$ is isomorphic to
$(SH_*^{\#}(\lambda H',J'))$ as filtered directed systems.
The isomorphism between these filtered directed systems
is a composition of restricted continuation morphisms
and switch morphisms and inverses of these morphisms.
\end{lemma}
This Lemma is similar to Lemma \ref{lemma:veryrestrictedgeneralgrowthrateequivalence}
with one constraint removed.
\proof of Lemma \ref{lemma:restrictivegeneralgrowthrateequivaalence}.
From Lemma \ref{lemma:positveactionhamiltonian} there is a pair $(H^p,J^p)$
constructed from $H$. This only depended on
what $H$ was on a small neighbourhood of $H^{-1}(0)$,
so the equivalent construction $({H'}^p,{J'}^p)$
is equal to $(H^p,J^p)$.
Hence all we need to do
in this section is to prove that the filtered directed system
$(SH_*^{\#}(\lambda H,J))$
is isomorphic to $(SH_*^{\#}(\lambda H^p,J^p))$.

There is a compact co-dimension $0$ submanifold $K'$
with boundary
containing $H^{-1}(0) = (H^p)^{-1}(0)$
and all the $1$-periodic orbits of $\lambda H^p$ for every
$\lambda \in (0,\infty) \setminus A$ (where $A$ is discrete).
We can assume that $K'$ is a subset of $U^H_1 = \cap_\lambda U^H_\lambda$.
It also has the property that any solution of the perturbed Cauchy-Riemann equations
with respect to $(H^p,J^p)$ with boundary in $K'$
is contained in $K'$.
Let $L : M \rightarrow \R$ be a Hamiltonian such that:
$L = H^p$ on a small neighbourhood of $K'$
and such that $L > 0$ outside $K'$.
We also require that $L = H$ outside a closed subset of $U^H_1$.
We let $J_L$ be an almost complex structure compatible with the symplectic form
such that it is equal to $J^p$ on a small neighbourhood
of $K'$ and equal to $J$ outside a closed subset of $U^H_1$.
Because $L = H$ on a small neighbourhood of
$H^{-1}(0) = L^{-1}(0)$ and $(L,J_L) = (H,J)$ outside a closed subset of $U^H_1$,
we have that $(L,J_L)$ is growth rate admissible.
Let $\theta$ be the $1$-form such that $\theta - \theta_M$
is exact
and $-\theta(X_H) - H$ is bounded.
Let $c > 0$ be a constant greater than the function
$-\theta(X_L) - L +1$.
Let $L' : M \rightarrow \R$ be a function which is
equal to $L$ on a small neighbourhood of $L^{-1}(0)$
and equal to $L + c$ on a small neighbourhood of $\overline{M \setminus K'}$
and is greater than $0$ everywhere else.
Because  $-\theta(X_{L'}) - L' < 0$ outside the interior of $K'$
and no $\lambda$-periodic orbits of $L'$ intersect the boundary of $K'$
for $\lambda$ outside some discrete set,
we have that all the $1$-periodic orbits of $\lambda L'$
of non-negative action are contained in the interior of $K'$.
Hence by part \ref{item:strongermaximumprinciple}
from Lemma \ref{lemma:positveactionhamiltonian},
we have that $(\lambda L',J_L)$ satisfies the maximum principle
with respect to $(K')^o$ (the interior of $K'$) for all $\lambda$ outside a discrete subset.
So we can define $U^{L'}_\lambda := (K')^o$.

Let $\tilde{L'}$ be equal to $L'$ inside $K'$
and equal to $H^p + c$ outside $K'$.
Again for all $\lambda$ outside a discrete subset,
all the $1$-periodic orbits of non-negative action
are contained in the interior $K'$.
Hence part (\ref{item:strongermaximumprinciple})
from Lemma \ref{lemma:positveactionhamiltonian},
ensures that any compact curve satisfying Floer's equations with
respect to $(\tilde{L}',J_L)$ whose boundary is contained in the interior of $K'$ must be contained
in $K'$.
Hence $(\tilde{L}',J_L)$ satisfies the maximum principle
with respect to $U^{\tilde{L'}}_\lambda := (K')^o$.
Because $U^{L'}_\lambda = U^{\tilde{L'}}_\lambda = (K')^o$
for all $\lambda$ and
$(L',J_L) = (\tilde{L}',J_L)$ on a small neighbourhood
of $K'$, we have a switch isomorphism from
$(SH_*^{\#}(\lambda L',J_L))$ to
$(SH_*^{\#}(\lambda \tilde{L'}, J^p))$.

Lemma \ref{lemma:veryrestrictedgeneralgrowthrateequivalence} tells us that
$(SH_*^{\#}(\lambda \tilde{L'}, J^p))$ and
$(SH_*^{\#}(\lambda H^p, J^p))$ are isomorphic by a restricted continuation morphism.
Lemma \ref{lemma:veryrestrictedgeneralgrowthrateequivalence} also tells us that
$(SH_*^{\#}(\lambda L', J))$ and
$(SH_*^{\#}(\lambda H, J))$ are isomorphic.
Hence $(SH_*^{\#}(\lambda H, J))$ and
$(SH_*^{\#}(\lambda H^p, J^p))$ are isomorphic.
\qed

If $(H_1,J_1),(H_2,J_2)$
are growth rate admissible such that 
$H_1^{-1}(0)$ contains $H_2^{-1}(0)$
then there is a natural morphism of filtered directed systems
from $(SH_*^{\#}(\lambda H_1,J_1))$ to $(SH_*^{\#}(\lambda H_2,J_2))$.
Here is how this morphism is constructed:
we have a non-decreasing family of Hamiltonians $H^s$ with
$H^s = H_1$ for $s \ll 0$ and such that $H^s = K$
for $s \gg 0$ where $K$ satisfies:
\begin{enumerate}
\item $K^{-1}(0) = H_2^{-1}(0)$ and $H_2 = K$ on a small neighbourhood of $K^{-1}(0)$.
\item $K = H_1$ outside a closed subset of $U^{H_1}_1$.
\item $(K,J_1)$ is growth rate admissible.
\end{enumerate}
We also require that $H^s = H_1$ outside a closed subset of $U^{H_1}_1$.
The family of pairs $(H^s,J_1)$ induces a continuation map from
$SH_*^{\#}(\lambda H_1,J_1)$ to $SH_*^{\#}(\lambda K,J_1)$.
This in turn induces a morphism of filtered directed systems from
$(SH_*^{\#}(\lambda H_1,J_1))$ to $(SH_*^{\#}(\lambda K,J_1))$ as the continuation
map commutes with the filtered directed system maps.
Our morphism is constructed by composing the above morphism
with the isomorphism from $(SH_*^{\#}(\lambda K,J_1))$
to $(SH_*^{\#}(\lambda H_2,J_2))$ from Lemma \ref{lemma:restrictivegeneralgrowthrateequivaalence}.
We call such a morphism a {\it growth rate admissible morphism}.
\begin{lemma} \label{lemma:growthratehamiltonianfunctoriality}
The composition of two growth rate admissible morphisms is a
growth rate admissible morphism.
\end{lemma}
\proof of Lemma \ref{lemma:growthratehamiltonianfunctoriality}.
The point is that composing switch morphisms gives us switch morphisms,
composing restricted continuation maps gives us restricted continuation maps and
if we have restricted continuation maps of the right form (so that some maximum principle
still applies) then restricted continuation maps can commute with switch morphisms.
Let $(H_1,J_1)$, $(H_2,J_2)$, $(H_3,J_3)$
be growth rate admissible pairs with $H_{i+1}^{-1}(0)$
contained in $H_i^{-1}(0)$.
We wish to show that the composition of the growth rate admissible morphisms
\[E_1 : (SH_*^{\#}(\lambda H_1,J_1)) \rightarrow (SH_*^{\#}(\lambda H_2,J_2))\]
\[E_2 : (SH_*^{\#}(\lambda H_2,J_2)) \rightarrow (SH_*^{\#}(\lambda H_3,J_3))\]
is a growth rate admissible morphism
\[E_3 : SH_*^{\#}(\lambda H_1,J_1) \rightarrow SH_*^{\#}(\lambda H_3,J_3).\]

Let $(H,J)$ be any growth rate admissible pair.
By using methods from the proof of Lemma
\ref{lemma:restrictivegeneralgrowthrateequivaalence}, we can construct Hamiltonians
$L_H^i$ and almost complex structures $J_H^i$ for $i=1,2,3$ with the following properties:
\begin{enumerate}
\item $(L_H^i)^{-1}(0) = H_i^{-1}(0)$ and $H_i = L_H^i$ near $H_i^{-1}(0)$.
\item $-\theta_M(X_{L_H^i}) - L_H^i \ll 0$ outside a small neighbourhood
$K_i$ of $H_i^{-1}(0)$.
\item Any map $s : S \rightarrow M$ with boundary inside $K_i$
satisfying Floer's equation with respect to $(L_H^i,J_H^i)$ must be contained in $K_i$.
This statement is true even if we change $(L_H^i,J_H^i)$ so that it equal to something
else outside $K_i$ as long as the Hamiltonian stays positive.
\item $K_3 \subset K_2 \subset K_1$.
\item $(L_H^2,J_H^2)$ is equal to $(L_H^3,J_H^3)$ outside $K_3$ and near $\partial K_3$
and $(L_H^1,J_H^1)$ is equal to $(L_H^2,J_H^2)$ outside $K_2$ and near $\partial K_2$.
\item $(L_H^i,J_H^i)$ is equal to $(H + \kappa,J)$ for some large constant $\kappa \gg 0$
near infinity.
\end{enumerate}
Note that the pair $(L_H^i,J_H^i)$ really depends on both $H$ and $J$
but we suppressed $J$ from the notation to make it less cluttered.
We now construct a smooth non-decreasing $1$-parameter family of pairs
$(L_H^t,J_H^t)$ for $t \in [1,3]$ such that
\begin{enumerate}
\item $(L_H^i,J_H^i)$ are equal to the pairs constructed above with the same name for $i = 1,2,3$.
\item $(L_H^t,J_H^t) = (L_H^3,J_H^3)$ outside $K_3$ and near $\partial K_3$.
\item For $t \in [1,2]$, $(L_H^t,J_H^t) =(L_H^2,J_H^2)$ outside $K_2$ and near $\partial K_2$.
\end{enumerate}
This family gives us continuation maps:
\[\Psi^H_{i,j} : SH_*^{\#}(\lambda L_H^i,J_H^i)) \rightarrow SH^+(\lambda L_H^j,J_H^j)\] that do not depend
on the choice of $(H,J)$ for any $i \leq j$,$i,j = 1 \text{ or } 2 \text{ or 3}$
because no Floer trajectories connecting non-negative action orbits escapes the
region $K_3$.
The point is that if we change $(H,J)$, we only do this outside $K_3$.
The map $\Psi^H_{i,j}$ induces a morphism of filtered directed systems
from $(SH_*^{\#}(\lambda L_H^i,J_H^i)))$ to $(SH^+(\lambda L_H^j,J_H^j))$.
The filtered directed systems are canonically isomorphic for any choice of $(H,J)$
via a switch morphism
and the maps $\Psi^H_{i,j}$ are exactly the same for any choice of $(H,J)$.
Hence we have filtered directed system maps $A_{i,j}^k$ that are isomorphisms
from $(SH_*^{\#}(\lambda L_{H_i}^k,J_{H_i}^k)))$ to $(SH^+(\lambda L_{H_j}^k,J_{H_j}^k))$.
And we also have the relations
\begin{equation} \label{equation:switchconjugate}
A_{k,k'}^j \circ \Psi^{H_k}_{i,j} \circ (A_{k,k'}^i)^{-1} = \Psi^{H_{k'}}_{i,j}.
\end{equation}

We also have a restricted continuation morphism $\Phi_i$ of from the
filtered directed system
$SH_*^{\#}(\lambda L_{H_i}^i,J_{H_i}^i))$ to $SH_*^{\#}(\lambda H_i,J_i))$
(see Lemma \ref{lemma:veryrestrictedgeneralgrowthrateequivalence}).
Its inverse $\Phi_i^{-1}$ is also a restricted continuation morphism.
We have that $\Psi^{H_2}_{2,3} \circ \Psi^{H_2}_{1,2}$
is induced entirely from continuation maps and so it is exactly the
same as the map $\Psi^{H_2}_{1,3}$.
Also composing $\Phi_i$ or $\Phi_i^{-1}$ with any of these maps (when possible)
also induces restricted continuation morphisms.
We get that the growth rate admissible morphism $E_1$ is equal to
\[\Phi_2 \circ A_{1,2}^2 \circ \Psi^{H_1}_{1,2} \circ \Phi_1^{-1} \]
and $E_2$ is equal to
\[\Phi_3 \circ A_{2,3}^3 \circ \Psi^{H_2}_{2,3} \circ \Phi_2^{-1}\]
by looking at the proof of Lemma \ref{lemma:restrictivegeneralgrowthrateequivaalence}
and the definition of growth rate admissible morphism.
Hence their composition $E_2 \circ E_1$ is equal to:
\[\Phi_3 \circ A_{2,3}^3 \circ \Psi^{H_2}_{2,3} \circ A_{1,2}^2 \circ \Psi^{H_1}_{1,2} \circ \Phi_1^{-1}.\]
By equation \ref{equation:switchconjugate} we have that this composition
is equal to
\[\Phi_3 \circ A_{2,3}^3 \circ A_{1,2}^3 \circ \Psi^{H_1}_{2,3} \circ \Psi^{H_1}_{1,2} \circ \Phi_1^{-1} =\]
\[\Phi_3 \circ A_{2,3}^3 \circ A_{1,2}^3 \circ \Psi^{H_1}_{1,3} \circ \Phi_1^{-1}.\]
The composition of the switch morphisms $A_{2,3}^3 \circ A_{1,2}^3$ is the switch morphism
$A_{1,3}^3$ hence we have that $E_2 \circ E_1$ is equal to:
\[\Phi_3 \circ A_{1,3}^3 \circ \Psi^{H_1}_{1,3} \circ \Phi_1^{-1}\]
which is equal to $E_3$.
Hence $E_2 \circ E_1 = E_3$ which gives us functoriality.
\qed

\begin{lemma} \label{lemma:growthrateadmissibleisomorphism}
Let $(H_s,J_s)$ be a smooth family of growth rate admissible
pairs parameterized by $s \in \R$ such that
for $s_1 < s_2$, $H_{s_1}^{-1}(0)$
contains $H_{s_2}^{-1}(0)$.
For $s_2$ greater
than $s_1$, the growth rate admissible morphism
from $(SH_*^{\#}(\lambda H_{s_1},J_{s_1}))$ to
$(SH_*^{\#}(\lambda H_{s_2},J_{s_2}))$ is an isomorphism.
\end{lemma}
\proof of Lemma \ref{lemma:growthrateadmissibleisomorphism}.
By a compactness argument and by functoriality
(Lemma \ref{lemma:growthratehamiltonianfunctoriality})
we only need to show that for each $s \in \R$, there is a $\delta_s > 0$
such that the result is true for $s_1 > s-\delta_s$
and $s_2 < s+\delta_s$.
We will prove this by changing our pairs $(H_s,J_s)$ to ones
similar to the one described in Lemma \ref{lemma:positveactionhamiltonian}.
Then we note that if all the orbits have non-negative action then
$SH_*$ is equal to $SH_*^{\#}$ hence we can construct an inverse to
our growth rate admissible morphism by using a decreasing family of Hamiltonians.

Fix the Hamiltonian $H_s$.
We can ensure that there is a continuous family of small neighbourhoods $U_{s'}$ of $H_{s'}^{-1}(0)$ in which
$-\theta(X_{H_{s'}}) - H_{s'} > 0$ on $U_{s'} \setminus H_{s'}^{-1}(0)$
with respect to $s'$
by a $1$-parameter version of the note after Lemma
\ref{lemma:positveactionhamiltonian}.
This means that there is a $\delta_s > 0$ and a neighbourhood $U$ of $H_{s'}^{-1}(0)$ with 
$-\theta(X_{H_{s'}}) - H_{s'} > 0$ on $U \setminus H_{s'}^{-1}(0)$
for all $|s-s'|<\delta_s$.

We first perturb $H_{s-\delta_s}$ by a $C^\infty$ small amount near its zero set to $K$
so that $(K,J_{s-\delta_s})$ is still growth rate admissible and so that
the interior of $K^{-1}(0)$ contains $H_{s-\delta_s}^{-1}(0)$ and $K^{-1}(0)$ is contained in $U$.
By Lemma \ref{lemma:positveactionhamiltonian} let $(K^p,J^p_{s-\delta_s})$ be a pair such that
\begin{enumerate}
\item $(K^p,J^p_{s-\delta_s})$ is growth rate admissible.
\item $K^p = K$ on a small neighbourhood of
$K^{-1}(0)$ and $(K^p)^{-1}(0)= K^{-1}(0)$.
\item $-\theta(K^p) - K^p \geq 0$ everywhere
and $-\theta(K^p) - K^p > 0$ when $K^p > 0$.
\end{enumerate}
We shrink $U$ so that it is also contained in $U^K_1$.
Choose an open set $U'$ whose closure is contained in $U$
and which still contains $K^{-1}(0)$.
Let $\rho : M \rightarrow \R$ be a bump function which
is equal to $1$ inside $U'$ and equal to $0$ outside $U$.
For $\kappa > 0$ large enough and for $|s - s'| \leq \delta_s$ we have that:
\[H'_{s'} := \rho H_{s'} + \kappa K^p\]
satisfies
$-\theta(H'_{s'}) - H'_{s'} \geq 0$ because
$-\theta(K^p) - K^p > 0$ in the
relatively compact region $U \setminus U'$
and so is bounded below by a positive constant.
We define a new family of almost complex structures $J'_{s'}$
to be equal to $J^p_{s-\delta_s}$ outside $U'$
and equal to $J_{s'}$ on a small neighbourhood of $H_{s'}^{-1}(0)$.
We have that $(H'_{s'},J'_{s'})$ is growth rate admissible.

We have that $SH_*(\lambda H'_{s'},J'_{s'})$
is a filtered directed system isomorphic to
$SH_*^{\#}(\lambda H'_{s'},J'_{s'})$
because $-\theta(H'_{s'}) - H'_{s'} \geq 0$.
Also we have that if $(K',Y')$ is any pair equal to
$(H'_{s'},J'_{s'})$ outside $U'$ then $SH_*(\lambda K',Y')$
is equal to
$SH_*(\lambda H'_{s'},J'_{s'})$ as the maximum principle
ensures that continuation maps between these Hamiltonians are well defined
(and that these continuation maps do not have to be from non-decreasing families
of Hamiltonians as we are not considering action).
Hence for all $|s'-s| \leq \delta_s$,
the filtered directed systems
$(SH_*(\lambda H'_{s'},J'_{s'}))$ are all isomorphic.

Let $s_1,s_2$ be such that $s- \delta_s < s_1 < s_2 < s + \delta_s$.
We have a morphism $\phi$ of filtered directed systems from
$(SH_*^{\#}(\lambda H'_{s_1},J'_{s_1}))$
to
$(SH_*^{\#}(\lambda H'_{s_2},J'_{s_2}))$.
By the previous discussion we also have a morphism $\phi'$ from
\[(SH_*^{\#}(\lambda H'_{s_2},J'_{s_2})) \cong (SH_*(\lambda H'_{s_2},J'_{s_2}))\]
to
\[(SH_*^{\#}(\lambda H'_{s_1},J'_{s_1})) \cong (SH_*(\lambda H'_{s_1},J'_{s_1})).\]
Because these morphisms are induced by continuation maps, we have
that $\phi \circ \phi'$ and $\phi' \circ \phi$ are filtered directed system maps
and hence by definition $\phi$ is an isomorphism of filtered directed systems.

By Lemma \ref{lemma:restrictivegeneralgrowthrateequivaalence}, we have an isomorphism $\Phi_i$
of filtered directed systems from $(SH_*^{\#}(\lambda H_{s_i},J_{s_i}))$ to
$(SH_*^{\#}(\lambda H'_{s_i},J'_{s_i}))$ for $i = 1,2$.
This isomorphism and its inverse are growth rate admissible morphisms.
We have an isomorphism $\Phi_2^{-1} \circ \phi \circ \Phi_1$
from $(SH_*^{\#}(\lambda H_{s_1},J_{s_1}))$ to
$(SH_*^{\#}(\lambda H_{s_2},J_{s_2}))$.
Because $\phi,\Phi_1,\Phi_2^{-1}$ are growth rate admissible morphisms, we have by
functoriality (Lemma \ref{lemma:growthratehamiltonianfunctoriality})
that $\Phi_2^{-1} \circ \phi \circ \Phi_1$
is also a growth rate admissible morphism.
Hence the natural morphism from
$(SH_*^{\#}(\lambda H_{s_1},J_{s_1}))$ to
$(SH_*^{\#}(\lambda H_{s_2},J_{s_2}))$ is an isomorphism.
\qed

\begin{lemma} \label{lemma:growthrateadmissiblemisomorphism2}
Let $(H_0,J_0),(H_1,J_1)$ be two growth rate admissible Hamiltonians
such that $H_1^{-1}(0) \subset H_0^{-1}(0)$.
Suppose that there is a smooth family of Hamiltonians
$G_t$ satisfying the Liouville vector field property
such that
\begin{enumerate}
\item $G_{s_2}^{-1}(0) \subset G_{s_1}^{-1}(0)$
if $s_1 \leq s_2$.
\item $G_0 = H_0$ on a neighbourhood of $H_0^{-1}(0)$
and $G_0^{-1}(0) = H_0^{-1}(0)$.
\item $G_1 = H_1$ on a neighbourhood of $H_1^{-1}(0)$
and $G_1^{-1}(0) = H_1^{-1}(0)$.
\end{enumerate}
Then the growth rate admissible morphism
\[(SH_*^{\#}(\lambda H_0,J_0)) \rightarrow (SH_*^{\#}(\lambda H_1,J_1))\]
is an isomorphism.
\end{lemma}
\proof of Lemma \ref{lemma:growthrateadmissiblemisomorphism2}.
By using construction similar to the one in the proof of Lemma \ref{lemma:growthrateadmissibleisomorphism}
combined with Lemma \ref{lemma:positveactionhamiltonian},
we have for each $s \in [0,1]$ there is a $\delta_s>0$
and a smooth family of growth rate admissible pairs $(G_s^t,J_s^t)$ 
($t \in (s-\delta_s,s-\delta_s)$) satisfying:
\begin{enumerate}
\item $G_s^t = G_t$ on a neighbourhood of $G_t^{-1}(0)$ and
 $G_t^{-1}(0) = (G^t_s)^{-1}(0)$.
\item $(G_0^0,J_0^0) = (H_0,J_0)$ on a neighbourhood
of $H_0^{-1}(0)$ and $(G_0^0)^{-1}(0) = H_0^{-1}(0)$.
\item $(G_1^1,J_1^1) = (H_1,J_1)$ on a neighbourhood
of $H_1^{-1}(0)$ and $(G_1^1)^{-1}(0) = H_1^{-1}(0)$.
\end{enumerate}
Hence by Lemma \ref{lemma:growthrateadmissibleisomorphism},
we have for $s-\delta_s < t_1 \leq t_2 < t+\delta_s$ the natural filtered directed system map from
$(SH_*^{\#}(G^{t_1}_s,J^{t_1}_s))$ to
$(SH_*^{\#}(G^{t_2}_s,J^{t_2}_s))$
is an isomorphism.

Hence by a compactness argument there is a sequence
of growth rate admissible pairs $(A_i,Y_i) := (G_{s_i}^{t_i},J_{s_i}^{t_i})$
for $i = 1,\cdots,k$
such that
\begin{enumerate}
\item  $A_i^{-1}(0)$ contains $A_{i+1}^{-1}(0)$.
\item $(A_k,Y_k) = (G_1^1,J_1^1)$ and $(A_0,Y_0) = (G_0^0,J_0^0)$.
\item The morphism from $(SH_*^{\#}(A_i,Y_i))$ to
$(SH_*^{\#}(A_{i+1},Y_{i+1}))$ is an isomorphism.
\end{enumerate}
We also have by Lemma \ref{lemma:restrictivegeneralgrowthrateequivaalence}
that the growth rate admissible morphisms
\[(SH_*^{\#}(H_0,J_0)) \rightarrow (SH_*^{\#}(A_0,Y_0))\]
and
\[(SH_*^{\#}(A_1,Y_1)) \rightarrow (SH_*^{\#}(H_1,J_1))\]
are isomorphisms.
Hence by functoriality of these morphisms we get that the growth rate admissible morphism from
$(SH_*^{\#}(\lambda H_0,J_0))$ to $(SH_*^{\#}(\lambda H_1,J_1))$
is an isomorphism.
\qed

\begin{lemma} \label{lemma:growthrateadmissiblelargezeroset}
Let $(H,J)$ be a growth rate admissible pair. For any compact set $K \subset M$,
there exists a growth rate admissible pair $(H_K,J_K)$ such that:
\begin{enumerate}
\item $H_K^{-1}(0)$ contains both $H^{-1}(0)$ and $K$.
\item The morphism from $(SH_*^{\#}(\lambda H_K,J_K))$ to $(SH_*^{\#}(\lambda H,J))$
is an isomorphism.
\end{enumerate}
\end{lemma}
\proof of Lemma \ref{lemma:growthrateadmissiblelargezeroset}.
For the purposes of this proof we may as well enlarge $K$
so that it contains $H^{-1}(0)$.
By Lemma \ref{lemma:growthrateadmissiblemisomorphism2}
all we need to do is create a
family of Hamiltonians $H_s$,($s \in [0,\infty]$) such that
\begin{enumerate}
\item $H_s$ satisfies the Liouville vector field property.
\item $H_{s_2}^{-1}(0)$ contains $H_{s_1}^{-1}(0)$ if $s_1 \leq s_2$.
\item $H_0 = H$ and $H_s^{-1}(0)$ contains $K$ for $s \gg 1$.
\end{enumerate}
Let $f_H,X_{\theta_H}$ be the function and Liouville vector field which
enables $(H,J)$ to satisfy the Liouville vector field property.
Let $\rho$ be a bump function such that it is equal to $1$
on a neighbourhood of $K$ and $0$ outside some larger compact set.
Because the vector field $V'_H := \rho X_{\theta_H}$ has compact support,
its flow $\phi^{V'_H}_t$ is well defined everywhere.
Let $H_s := ((\phi^{V'_H}_s)_* H)$.
For $\epsilon_H>0$ small enough we have that
$dH(X_{\theta_H}) > 0$ inside $H^{-1}(0,\epsilon_H)$.
Because $(\phi^{V'_H}_t)^* X_{\theta_H}$ is proportional to $X_{\theta_H}$
we have that $dH_s(X_{\theta_H}) =  (\phi^{V'_H}_t)_*\left(dH((\phi^{V'_H}_t)^* X_{\theta_H})\right) > 0$
in the region $H_s^{-1}(0,\epsilon_H)$.
Also because $H$ satisfies the Liouville vector field
property, there exists a $C_s$ such that
$f_H^{-1}(-\infty,C_s]$ is contained in $H_s^{-1}(0)$
and $df_H(X_{\theta_H}) >0$ in $f_H^{-1}(C_s,\infty)$.
Hence $H_s$ satisfies the Liouville vector field property for all $s \geq 0$.
There is a constant $C$ so that $f_H^{-1}(-\infty,C] \subset H^{-1}(0)$
and $df_H(V'_H)>0$ in the region $K \setminus f_H^{-1}(-\infty,C)$.
Hence the time $s$ flow of $f_H^{-1}(-\infty,C]$ contains $K$ for $s \gg 1$.
Because $f_H^{-1}(-\infty,C]$ is contained inside $H^{-1}(0)$
we have for $s \gg 1$ that $H_s^{-1}(0)$ contains $K$.

By Lemma \ref{lemma:positveactionhamiltonian}
there is a growth rate admissible pair $(H_K,J_K)$
such that $H_K^{-1}(0) = G_s^{-1}(0)$ for some $s \gg 0$
and such that $H_K=G_s$ on a neighbourhood of $H_K^{-1}(0)$.
We assume that $s$ is large enough so that $H_K^{-1}(0)$ contains $K$.
Hence by Lemma \ref{lemma:growthrateadmissiblemisomorphism2} we have that the growth rate admissible morphism
\[(SH_*^{\#}(\lambda H_K,J_K)) \rightarrow (SH_*^{\#}(\lambda H,J))\]
is an isomorphism.
\qed

\proof of Theorem \ref{theorem:generalgrowthrateequivalence}.
Let $K_1$ be a compact set
whose interior
contains both $H^{-1}(0)$ and ${H'}^{-1}(0)$.
By Lemma \ref{lemma:growthrateadmissiblelargezeroset}, there is a pair $(H_1,J_1)$
with $H_1^{-1}(0)$ containing $K_1$ such that
the growth rate admissible morphism from $(SH_*^{\#}(\lambda H_1,J_1))$ to
$(SH_*^{\#}(\lambda H,J))$ is an isomorphism.
Choose a compact set $K_2$ whose interior contains
$H_1^{-1}(0)$.
Let $(H'_1,J'_1)$ be a pair such that $(H'_1)^{-1}(0)$ contains
$K_2$ and such that the morphism from 
$(SH_*^{\#}(\lambda H'_1,J'_1))$ to
$(SH_*^{\#}(\lambda H',J'))$ is an isomorphism.
By repeating this process, we can find two more growth rate admissible pairs
$(H_2,J_2)$,$(H'_2,J'_2)$ such that:
\begin{enumerate}
\item $(H'_1)^{-1}(0)$ is contained in the interior of $H_2^{-1}(0)$.
\item The morphism from $(SH_*^{\#}(\lambda H_2,J_2))$ to
$(SH_*^{\#}(\lambda H_1,J_1))$ is an isomorphism.
\item $H_2^{-1}(0)$ is contained in the interior of $(H'_2)^{-1}(0)$.
\item The morphism from $(SH_*^{\#}(\lambda H'_2,J'_2))$ to
$(SH_*^{\#}(\lambda H'_1,J'_1))$ is an isomorphism.
\end{enumerate}
Hence we have the following sequence of morphisms:
\[(SH_*^{\#}(\lambda H'_2,J'_2)) \rightarrow
(SH_*^{\#}(\lambda H_2,J_2)) \rightarrow \]
\[(SH_*^{\#}(\lambda H'_1,J'_1)) \rightarrow
(SH_*^{\#}(\lambda H_1,J_1))\]
where composing any two of these morphisms gives an isomorphism.
By Theorem  \ref{lemma:equivalentfiltereddirectedsystems},
the middle morphism from the group 
$(SH_*^{\#}(\lambda H_2,J_2))$ to
$(SH_*^{\#}(\lambda H'_1,J'_1))$ is an isomorphism.
Hence $(SH_*^{\#}(\lambda H,J))$ is isomorphic
to $(SH_*^{\#}(\lambda H',J'))$.
This proves the theorem.
\qed

\begin{corollary} \label{cor:preliminarygrowthrateinvariance}
The filtered directed system $(SH_*^{\#}(\lambda H, J))$
is an invariant of $M$ up to exact symplectomorphism as long as this symplectomorphism
preserves our choice of trivialization $\tau$ of the canonical bundle and our class $b \in H^2(M,Z / 2\Z)$.
If $M$ is complete (i.e. the $\omega_M$-dual
of $\theta_M$ is an integrable vector field)
then it is an invariant up to general symplectomorphism
(again preserving $(\tau,b)$).
\end{corollary}
\proof
All the properties defining growth rate admissibility,
are invariants of $M$ up to exact symplectomorphism, hence
by Theorem \ref{theorem:generalgrowthrateequivalence} we have
that $(SH_*^{\#}(\lambda M,J))$ is an invariant up to exact symplectomorphism
preserving $(\tau,b)$.
Suppose that $M$ is complete.
Then it is the completion of some Liouville domain $\overline{M}$
hence by \cite[Lemma 1]{eliashberg:symplectichomology},
if $M,M'$ are symplectomorphic then they are exact symplectomorphic.
Hence $(SH_*^{\#}(\lambda M,J))$ is an invariant up to general
symplectomorphism preserving $(\tau,b)$ in this case.
\qed

Because its an invariant up to exact symplectomorphism preserving $(\tau,b)$, we will write:
\[(SH_*^{\#}(M,\theta,\lambda))\] for any filtered directed system
$(SH_*^{\#}(\lambda H,J))$ where we have chosen
some growth rate admissible pair $(H,J)$ and a pair $(\tau,b)$.

If $(N,\theta_N)$ is a Liouville domain, then the
interior of $N$, $N^0$ is a finite type convex symplectic manifold for the
following reason:
Let $X_{\theta_N}$ be the  $d\theta_N$-dual of $\theta_N$.
By flowing back $\partial N$ backwards along $X_{\theta_N}$,
we get that a collar neighbourhood of $\partial N$
is equal to $(1-\epsilon,1] \times \partial N$ with
$\theta_N = r_N \alpha_N$.
Here $r_N$ parameterizes the interval and $\alpha_N = \theta_N|_{\partial_N}$.
Let $g : (1-\epsilon,1) \rightarrow \R$ be a function which
is equal to $0$ near $1-\epsilon$ and tends to $+\infty$ near $1$
and also that its derivative is positive near $1$.
We let $f_N : N^0 \rightarrow \R$ be a function which is $0$
away from this collar neighbourhood and equal to $g(r_N)$
inside this collar neighbourhood.
This shows that $(N^0,\theta_N)$ has the structure of a finite type convex symplectic manifold.
\begin{lemma} \label{lemma:isomorphictoliouvillesubdomain}
Let $(M,\theta)$ be a finite type convex symplectic manifold
and let $X_{\theta}$ be the $d\theta$-dual of $\theta$.
Let $f_M : M \rightarrow \R$ be the exhausting function such that
$df_M(X_{\theta})>0$ in $f_M^{-1}[C,\infty)$ for some $C \gg 0$.
We define $N$ to be the Liouville domain $f_M^{-1}(-\infty,C]$.
We also assume that the period spectrum of the contact manifold $\partial N$
is discrete.
Then the filtered directed systems
$(SH_*(M,\theta,\lambda))$ and $(SH_*(N^0,\theta|_{N^0},\lambda))$
are isomorphic as long as the choice of trivialization $\tau$ and homology class $b$
for $N^0$ is equal to such a choice for $M$ restricted to $N^0$.
\end{lemma}
\proof of Lemma \ref{lemma:isomorphictoliouvillesubdomain}.
Let $U \subset [1,\infty) \times \partial N$ be the
partial cylindrical end of $M$ obtained by flowing
$\partial N$ along $X_{\theta}$.
Let $r$ be the coordinate parameterizing $[1,\infty)$.
We have that $\theta = r \alpha$ inside $U$ where $\alpha = \theta|_{\partial N}$.
By flowing $\partial N$ backwards along $X_{\theta}$
we can extend $U$ to $U'$ (containing $\overline{M \setminus N}$)
so that it is now a subset of $(0,\infty) \times \partial N$.
We also extend $r$ so it now parameterizes the larger interval
$(0,\infty)$.
Let $h : (0,\infty) \rightarrow \R$ be a function
such that $h(r) = 0$ near $r = 0$ and $h(r) = r$
for $r \geq 1-\delta$ where $\delta>0$ is small.
We also assume that $h \geq 0$ and $h' > 0$ for $h>0$.
We set $H$ to be equal to $h(r)$ where $r$ is well defined and $0$ otherwise.
We let $J$ be an almost complex structure on $M$
such that it is cylindrical on the region $\{r \geq 1-\delta\}$.
The pair $(H,J)$ is growth rate admissible
(this is because it is basically the same as the pair from
the second example mentioned earlier).
The pair $(H|_{N^0},J|_{N^0})$ is also growth rate admissible.
Let $A \subset (0,\infty)$ be the period spectrum
of $\partial N$.
If $\lambda \in (0,\infty) \setminus A$ then all the $1$-periodic
orbits of $\lambda H$ are contained in $N^0$.
Also because $J$ is cylindrical, by Corollary \ref{corollary:maximumprinciplerescaling}
all the Floer trajectories connecting orbits of $\lambda H$ or
continuation map Floer trajectories joining
$\lambda_1 H$ and $\lambda_2 H$ are contained in $N^0$.
This ensures that the filtered directed systems
$(SH_*^{\#}(\lambda H,J))$ and $(SH_*^{\#}(\lambda H|_{N^0},J|_{N^0}))$
are isomorphic.
This completes the Lemma.
\qed

\begin{lemma} \label{lemma:convexdeformationequivalence}
Suppose that $(M,\theta)$ and $(M',\theta')$ are convex deformation
equivalent, then
$(SH_*^{\#}(M,\theta,\lambda))$ and $(SH_*^{\#}(M',\theta',\lambda))$
are isomorphic as filtered directed systems
(again the choice of trivialization $\tau$
and homology class $b$ for $M'$ must be the same as that of $M$).
\end{lemma}
\proof of Lemma \ref{lemma:convexdeformationequivalence}.
Let $f_M,X_{\theta},N$ be as in the previous Lemma.
By Corollary \ref{lemma:finitetypecompleting}, we have that
$(M,\theta)$ is convex deformation equivalent to
the completion $(\widehat{N},\theta_N)$.
Let $r_N$ be the cylindrical coordinate of $\widehat{N}$.
We can extend the cylindrical end $[1,\infty) \times \partial N$
to $(0,\infty) \times \partial N$ inside $\widehat{N}$ by
flowing $\partial N$ backwards along the Liouville vector field $X_{\theta_N}$.
By applying Lemma \ref{lemma:isomorphictoliouvillesubdomain} twice
(once to $(M,\theta)$ and once to $(\widehat{N},\theta_N)$)
we get that $(SH_*^{\#}(M,\theta,\lambda))$
is isomorphic to $(SH_*^{\#}(N^0,\theta|_N=\theta_N|_N,\lambda))$
which is isomorphic to $(SH_*^{\#}(\widehat{N},\theta_N,\lambda))$.

Similarly we have that there is a Liouville domain $N'$
such that the filtered directed system $(SH_*^{\#}(M',\theta',\lambda))$
is isomorphic to
$(SH_*^{\#}(\widehat{N'},\theta_{N'},\lambda))$ and such that
$(M',\theta')$ is convex deformation equivalent to $(\widehat{N'},\theta_{N'})$.
Because convex deformation equivalence is an equivalence relation we get
that $(\widehat{N},\theta_N)$ is convex deformation equivalent to
$(\widehat{N'},\theta_{N'})$.
Also both $(\widehat{N},\theta_N)$ and $(\widehat{N'},\theta_{N'})$ are complete
so by Corollary \ref{cor:exactsymplectomorphicconvexdeformation}
they are exact symplectomorphic.
Hence by Lemma \ref{cor:preliminarygrowthrateinvariance} we have that
the groups
$(SH_*^{\#}(\widehat{N},\theta_N,\lambda))$ and
$(SH_*^{\#}(\widehat{N'},\theta_{N'},\lambda))$ are isomorphic.
This implies that $(SH_*(M,\theta,\lambda))$ and $(SH_*(M',\theta',\lambda))$
are isomorphic.
\qed

Motivated by Lemma \ref{cor:preliminarygrowthrateinvariance}
we have the following definition:
\begin{defn} \label{defn:flexiblegrowthrate}
For any finite type convex symplectic manifold $M$, we define
\[\Gamma(M,\theta) := \Gamma( (SH_*^{\#}(M,\theta,\lambda) ) ).\]
Sometimes we write $\Gamma(M)$ if it is clear what the Liouville form
$\theta$ is.
\end{defn}
We will show later in Corollary \ref{lemma:growthrateequalsactiongtrowthrate} that this is the
same as growth rate as in Definition \ref{defn:growthrate}.

In some cases, we wish to consider orbits of all actions
and not just ones of non-negative action.
Let $(H,J)$ be growth rate admissible
such that $-\theta_H(X_H) - H \geq 0$ for some
$\theta_H$ where $\theta_H - \theta_M$ is exact.
Let $(H_\lambda,J_\lambda)$ be a smooth family
of Hamiltonians parameterized by $\lambda \geq 1$
such that
$(H_\lambda,J_\lambda) = (\lambda H + c_\lambda,J)$
outside a closed subset of $U^H_\lambda$. Here $c_\lambda$ is a smooth family of constants.
Basically by the maximum principle we have that
$SH_*(H_\lambda,J_\lambda)$ is well defined and
for $\lambda_1 < \lambda_2$ there is a morphism
from $SH_*(H_{\lambda_1},J_{\lambda_1})$
to $SH_*(H_{\lambda_2},J_{\lambda_2})$.
Note that outside a closed subset of $U^H_\lambda$, $\lambda H$ has
no 1-periodic orbits of negative action and hence no 1-periodic orbits
outside this closed subset.
This morphism is induced by the smooth family of pairs
$(H_\lambda,J_\lambda)$ from $\lambda_1$ to $\lambda_2$.
Hence $(SH_*(H_\lambda,J_\lambda))$ forms a filtered directed system.
\begin{lemma} \label{lemma:positiveactionfiltereddirectedsystem}
Suppose that $-\theta(X_H) - H \geq 0$ then $(SH_*(H_\lambda,J_\lambda))$
is isomorphic to $(SH_*^{\#}(M,\theta,\lambda))$.
\end{lemma}
\proof of Lemma \ref{lemma:positiveactionfiltereddirectedsystem}.
Because the action of all the $1$-periodic orbits
of $\lambda H$ are non-negative, we have that
$SH_*(\lambda H,J) \cong SH_*^{\#}(\lambda H,J)$.
This isomorphism commutes with the filtered directed
system maps because they are continuation maps
induced from an increasing family of Hamiltonians.
Hence the filtered directed system $(SH_*(\lambda H,J))$
is isomorphic to $(SH_*^{\#}(\lambda H,J))$.

Let $q : \R \rightarrow [0,1]$ be a smooth function
such that $q(x) = 0$ for $x \leq 0$ and $q(x) = 1$ for $x \geq 1$.
By joining $(\lambda H,J)$ with $(H_\lambda,J_\lambda)$ via a smooth
family of pairs $(H_\lambda^s,J_\lambda^s)$ such that
$(H_\lambda^s,J_\lambda^s) = (\lambda H + q(s) c_\lambda,J)$ outside
a closed subset of $U^H_\lambda$, we have
by the maximum principle a well defined continuation isomorphism
from $SH_*(\lambda H,J)$ to $SH_*(H_\lambda,J_\lambda)$.
This monorphism commutes with the continuation maps.
Also the continuation map induced by the family
$(H_\lambda^{-s},J_\lambda^{-s})$ gives us an inverse to the above
morphism. This is because the composition of these two continuation
maps is a continuation map induced by some family of pairs
equal to $(\lambda H + q(s) c'_\lambda, J)$ near infinity
and these are homotopic through such families of pairs
to the constant pair. The constant pair gives us the identity map.

Hence we have that the filtered directed system
$(SH_*(\lambda H,J))$ is isomorphic to $(SH_*(H_\lambda,J_\lambda))$.
Hence we have that
$(SH_*^{\#}(\lambda H,J))$ is isomorphic to $(SH_*(H_\lambda,J_\lambda))$.
This proves the Lemma.
\qed

\begin{corollary} \label{cor:lineardirectedsystemhaspositiveaction}
Let $N$ be a Liouville domain and let $r_N$ be the cylindrical
coordinate of $\widehat{N}$.
Suppose also that the period spectrum of $\partial N$ is discrete.
Let $J$ be an almost complex structure that is cylindrical near infinity.
Let $(H_\lambda,J_\lambda)$ be a family of Hamiltonians such that
$H_\lambda = \lambda r_N + c_\lambda$ outside a large compact
set $K$ for some family of constants $c_\lambda$.
Suppose also that $J_\lambda$ is equal to $J$ outside $K$.
Then $(SH_*(H_\lambda,J_\lambda))$ is isomorphic as a filtered
directed system to $(SH_*^{\#}(\widehat{N},\theta_N,\lambda))$.
\end{corollary}
\proof of Corollary \ref{cor:lineardirectedsystemhaspositiveaction}.
Let $h : [1,\infty) \rightarrow \R$ be a function such that
$h(x) = 0$ for $x$ near $1$ and $h(x) = x - 2$ for $x \geq 3$.
Suppose also that $h'(x),h''(x) \geq 0$.
Let  $H$ be a Hamiltonian such that $H = 0$ inside $N \subset \widehat{N}$
and $H = h(r_N)$ outside $N$.
Inside $N$, $-\theta(X_H) - H = 0$.
Outside $N$,
\[-\theta(X_H) - H = r_N h'(r_N) - h(r_N) = \int_1^{r_N} x h''(x) dx \geq 0.\]
Hence $-\theta(X_H) - H \geq 0$.
This pair is growth rate admissible for reasons similar
to the reason why the second example mentioned earlier is growth rate admissible.
Hence by Lemma \ref{lemma:positiveactionfiltereddirectedsystem} we
get our result.
\qed

The problem is that the growth rate $\Gamma(M,b)$ has a different definition
to the one given in definition \ref{defn:growthrate}.
We recall the definition here:
In section \ref{section:symplectichomology}, we defined $SH_*^{\leq \lambda}(N)$ for
a Liouville domain $N$ which is the direct limit of $SH_*^{\leq \lambda}(H,J)$
where $(H,J)$ is a pair defined on $\widehat{N}$ and is cylindrical at infinity
and less than $0$ on $N$.
For $\lambda_1 \leq \lambda_2$, there is a natural map
$SH_*^{\leq \lambda_1}(N) \rightarrow SH_*^{\leq \lambda_2}(N)$.
This is a filtered directed system $(SH_*^{\leq \lambda}(N))$
whose direct limit is $SH_*(N)$.

\begin{lemma} \label{lemma:growthrateequalsactiongtrowthrate}
The filtered directed system $(SH_*^{\leq \lambda}(N))$ is isomorphic
to \\ $(SH_*(\widehat{N},\theta_N,\lambda))$.
Hence definitions \ref{defn:growthrate} and \ref{defn:flexiblegrowthrate}
are equivalent.
We assume that the period spectrum on $N$ is discrete.
\end{lemma}
\proof of Lemma \ref{lemma:growthrateequalsactiongtrowthrate}.
Let $(H,J)$ be the pair defined in the proof
of Corollary \ref{cor:lineardirectedsystemhaspositiveaction}.
We will first construct a family of pairs $(H_\lambda,J_\lambda)$.
Let $r$ be the cylindrical coordinate of $N$.
We will construct an isomorphism from
$SH_*(\lambda H,J)$ to $SH_*^{\leq \lambda}(N)$
that commutes with the filtered directed system maps.

Fix $\lambda \geq 1$.
We construct the Hamiltonian $H_a$ as follows:
$H_a$ is constant and equal to $-\frac{1}{a}$
inside $N$.
This construction only works if $a$ is sufficiently large.
In the region $\widehat{N} \setminus N$
we let $H_a = g_a(r)$ where 
$g_a(r) = -\frac{1}{a}$ near $r = 1$,
$g'_a(r),g''_a(r) \geq 0$ and
$g_a(r) = a( r - 1 - \frac{1}{a} )$
for $r \geq 1 + \frac{1}{a}$.
We also require that in a neighbourhood of the region
$1 + \frac{2}{3a \lambda} \leq r \leq 1 + \frac{5}{6a \lambda}$
that $g(a)$ is equal to
\[-\frac{1}{a} + (\lambda - \frac{3}{2a})(r - 1 - \frac{1}{2 a \lambda - 3}).\]
Here is a picture:
\begin{figure}[H]
\centerline{
   \scalebox{1.0}{
    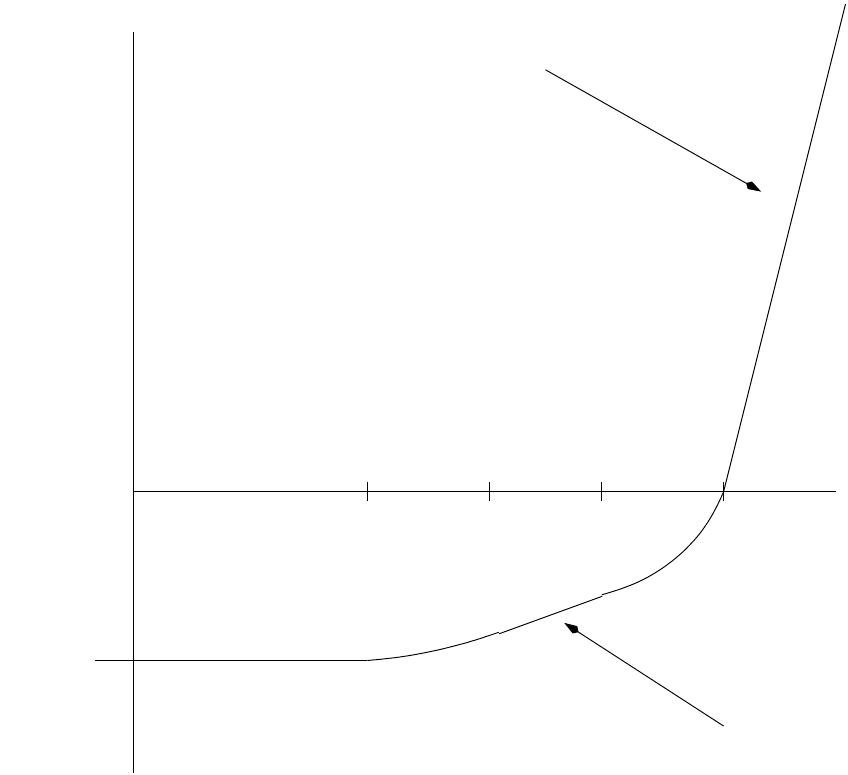
     }
      }
      \end{figure}
      
%
We let $J_a$ be anything we like inside $N$
and $J_a$ is cylindrical outside $N$.
If we have a periodic orbit inside the level set $r = c$,
then its action is $cg'(c) - g(c)$.
The derivative of this with respect to $c$
is $g'(c) + cg''(c) - g'(c) = cg''(c)$.
Because $g'' \geq 0$ we have that orbits in the level
set $r = c_1$ have action less than or equal to the ones
in the level set $r = c_2$ for $c_1 \leq c_2$.
Let $e = 1 + \frac{2}{3a \lambda}$.
All the $1$-periodic orbits in the region $r \leq e$
have action less than or equal to
$eg'(e) - g(e)$ and all the orbits in the region
$r \geq e$ have action greater than or equal to this same quantity.
Also because $g(r) = -\frac{1}{a} + (\lambda - \frac{3}{2a})(r - 1 - \frac{1}{2 a \lambda - 3})$
in the region $(e- \epsilon,e]$ for $\epsilon>0$ small,
\[eg'(e) - g(e)
= e(\lambda - \frac{3}{2a}) + \frac{1}{a} -(\lambda - \frac{3}{2a})(e - 1 - \frac{1}{2 a \lambda - 3})\]
\[= \frac{1}{a} + (\lambda - \frac{3}{2a})(1 + \frac{1}{2a \lambda -3})\]
\[= \frac{1}{a} + (\lambda - \frac{3}{2a})(1 + \frac{\frac{1}{2a}}{\lambda - \frac{3}{2a}}) \]
\[= \frac{1}{a} + \lambda - \frac{3}{2a} + \frac{1}{2a} = \lambda.\] 
Because $J$ is cylindrical and $H$ is linear in the region $(e- \epsilon,e]$,
we can use the maximum principle Corollary \ref{corollary:maximumprinciplerescaling}
to ensure that any Floer trajectory connecting orbits
of action $\leq \lambda$ must be contained in $\{r \leq e\}$.

We can construct another growth rate admissible pair $(H'_a,J'_a)$
where $(H'_a,J'_a) = (H_a,J_a)$ in the region $\{r \leq e\}$
and $H'_a =  -\frac{1}{a} + (\lambda - \frac{3}{2a})(r - 1 - \frac{1}{2 a \lambda - 3})$
outside this region.
Then $SH_*(H'_a,J'_a) = SH_*^{\leq \lambda}(H_a,J_a)$.
Also for $a_1 \leq a_2$, a similar maximum principle argument
ensures that the transfer map
\[SH_*(H'_{a_1},J'_{a_1}) \rightarrow SH_*(H'_{a_2},J'_{a_2})\]
coming from the family $(H'_{a_1 + t(a_2 - a_1)},J'_{a_1 + t(a_2 - a_1)})$
is identical to the transfer map
\[SH_*^{\leq \lambda}(H_{a_1},J_{a_1}) \rightarrow SH_*^{\leq
\lambda}(H_{a_2},J_{a_2})\]
coming from the family $(H_{a_1 + t(a_2 - a_1)},J_{a_1 + t(a_2 - a_1)})$.

Also because the slope of $H'_a$ is equal to $\lambda - \frac{3}{2a}$
near infinity we have that $SH_*(\lambda H,J)$ is isomorphic
to $\varinjlim_{a} SH_*(H'_a,J'_a)$.
Hence $SH_*(\lambda H,J) \cong SH_*^{\leq \lambda}(N)$
because $\varinjlim_{a} SH_*(H'_a,J'_a)$
is equal to $\varinjlim_{a} SH_*^{\leq \lambda}(H_a,J_a)$.
By ensuring that the pairs $(H_a,J_a)$,$(H'_a,J'_a)$ smoothly
increase with $\lambda$ and looking
at the resulting continuation maps
(and by Corollary \ref{corollary:maximumprinciplerescaling}
we get that the filtered directed system
maps for
$(SH_*(\lambda H,J))$
are identical to the ones for
$(SH_*^{\leq \lambda}(N))$.
Hence these directed systems are isomorphic.
Hence $(SH_*^{\#}(\widehat{N},\theta_N,\lambda))$ is isomorphic
to $(SH_*^{\leq \lambda}(N))$ by Corollary \ref{cor:lineardirectedsystemhaspositiveaction}.
\qed

We will now prove Theorem \ref{thm:growthrateinvariance}
which says that growth rate is an invariant of a finite
type Liouville manifold up to symplectomorphism.
\proof of Theorem \ref{thm:growthrateinvariance}.
By Corollary \ref{lemma:growthrateequalsactiongtrowthrate}
we have that growth rate is equal
to $\Gamma(SH_*^{\#}(M,\theta_M,\lambda))$.
By Lemmas \ref{cor:preliminarygrowthrateinvariance}
and \ref{lemma:growthrateuptoisomorphism},
we have that the growth rate $\Gamma(SH_*^{\#}(M,\theta_M,\lambda))$
is an invariant of $M$ up to symplectomorphism.
Hence $\Gamma(M)$ is an invariant of $M$ up to symplectomorphism preserving the class
$b \in H^2(M,\Z / 2\Z)$ and the choice of trivialization $\tau$.
\qed.

\subsection{Growth rate of cotangent bundles} \label{section:growthratecotangent}

We let $\K$ be a field.
Let $(Q,g)$ be a Riemannian manifold, and let
$L^\lambda := {\mathcal L}^{ \leq \lambda^2 }(Q,g)$
be the space of free loops of length $\leq \lambda^2$.
For $\lambda_1 \leq \lambda_2$ we have a natural inclusion
$L^{\lambda_1} \hookrightarrow L^{\lambda_2}$.
This gives us a filtered directed system
$(H_*(L^{\lambda},\K))$.
We define $\Gamma(Q,\K)$ to be equal to
$\Gamma( (H_*(L^{\lambda},\K)) )$.
Note that $Q$ has exponential growth if and only if
there is some field $\K$ such that $\Gamma(Q,\K) = \infty$.
We define $\omega_2 \in H^2(T^*Q,\Z / 2\Z)$ to
be the pullback of the second Stiefel-Whitney class of $Q$.
We also have a canonical choice of trivialization $\tau_Q$
of the canonical bundle of $T^*Q$ induced by the volume form on $Q$.
The cotangent bundle $T^*Q$ is the completion
of the unit cotangent bundle $D^* Q$ which is a Liouville domain.
Let $\theta_Q$ be the Liouville form on $T^*Q$.
This is locally equal to $\sum_i p_i dq_i$
where $p_i$ are momentum coordinates and $q_i$
are position coordinates in $Q$.
In this section we will prove:
\begin{theorem} \label{theorem:cotangentfiltereddirectedsystem}
The filtered directed systems $(H_*(L^{\lambda}))$
and \\ $(SH_*^{\#}(T^* Q,\theta_Q,\lambda,(\tau_Q,\omega_2)))$ are isomorphic.
\end{theorem}
This means we get the following Corollary:
\begin{corollary} \label{cor:loopspacegrowth}
$\Gamma(T^*Q,(\tau_Q,\omega_2)) = \Gamma(Q)$.
\end{corollary}

Before we prove this theorem we need a slightly different
definition of growth rate.
Let $(N,\theta_N)$ be a Liouville domain.
Let $r_N$ be the radial coordinate for the
cylindrical end $\partial N \times [1,\infty)$ of $\widehat{N}$.
We assume that the period spectrum of the contact boundary $\partial N$
is a discrete subset ${\mathcal P}$ of $\R$.
We say that a Hamiltonian $H : S^1 \times \widehat{N} \rightarrow \R$
is {\it quadratic admissible}
if:
\begin{enumerate}
\item
 there exists constants $b,b'\in \R$ such that
$H +b \leq \frac{1}{2}r_N^2 \leq H + b'$ in the region
$r_N \geq 1$.
We also require that:
$r_N + b \leq dH(\frac{\partial}{\partial r_N}) \leq r_N + b'$.
This ensures that if $H$ has some $1$-periodic orbit inside the region $\alpha \leq r_N \leq \beta$
then its action is greater than $\alpha(\alpha + b) + b' - \frac{1}{2}\beta^2$
and less than $\beta(\beta + b') + b - \frac{1}{2} \alpha^2$. 
\item
There is a sequence $l_1 < l_2 < \cdots$
tending to infinity such that $l_i \notin {\mathcal P}$
and $H = \frac{1}{2}r_N^2+d_i$ on a small neighbourhood of the hypersurface
$\{r_N = l_i\}$ where $d_i$ is a constant
(in particular $H$ has no $1$-periodic orbits near these hypersurfaces).
This also ensures that the $1$-periodic orbits starting at say $r_N = \alpha$
do not stray too far away from this hypersurface.
\item the sequence $l_i$ satisfies $l_i / l_{i-1} < \kappa$
where $\kappa$ is a constant.
\item We need that $H \geq l_i r_N -\frac{1}{2}l_i^2 + d_i$ in the region $r_N \geq l_i$.
\end{enumerate}
By Lemma \ref{lemma:perturbinghamiltonian} we can perturb any quadratic admissible
Hamiltonian so it becomes non-degenerate and remains quadratic admissible.
Let $H$ be such a Hamiltonian.
We also define an almost complex structure $J$
such that $J$ is equal to some cylindrical almost complex structure
on a neighbourhood of $\{r_N = l_i\}$ for each $i$.
We say that $J$ is compatible with $H$ if this is true.
Such a pair $(H,J)$ is said to be {\it quadratic admissible}.
Lemma \ref{lemma:maximumprinciple} tells us that we can define
$SH_*^{\leq \lambda}(H,J)$
(see \cite[Section 21.3]{Ritter:transfer}
for an alternative way of defining this).
These groups form a filtered directed system where all the filtered
directed system maps come from the natural inclusions.
Hence we have a filtered directed system
$(SH_*^{\leq \lambda^2}(H,J))$.
\begin{lemma} \label{lemma:fromonechaincomplex}
The filtered directed system $(SH_*^{\leq \lambda^2}(H,J))$
is isomorphic to $(SH_*^{\leq \lambda}(N))$.
\end{lemma}

We need a preliminary algebraic lemma.

\begin{lemma} \label{lemma:discreteisomorphism}
Let $(V_x),(V'_x)$ be two filtered directed systems and
let $1 \leq l_1 < l_2 < \cdots$ be a sequence tending
to infinity such that $l_i / l_{i-1} < K$ for some constant $K$.
Suppose for each $x \in [1,\infty)$ we have a map $p_{x,i}$ from
$V_x$ to $V'_{l_i}$ whenever $l_i \geq cx$ where $c \geq 1$ is some constant.
We assume that $p_x$ commutes with the filtered directed system maps
(i.e. if $a_{x,y} : V_x \rightarrow V_y$ and
$a'_{i,j} : V'_{l_i} \rightarrow V'_{l_j}$ are filtered directed system maps
then we require that $p_{y,j} \circ a_{x,y} = a'_{i,j} \circ a_{x,i}$).
Suppose we also have maps $p'_{i,x} : V'_{l_i} \rightarrow V_x$ for $l_i \leq c'x$
which also commute with the filtered directed system maps
(here $c' \geq 1$ is a constant).
Then $(V_x)$ is isomorphic to $(V'_x)$.
\end{lemma}
\proof of Lemma \ref{lemma:discreteisomorphism}.
Let $a_{x,y} : V_x \rightarrow V_y$ and
$a'_{x,y} : V'_x \rightarrow V'_y$ be the directed system maps.
We first construct a map $\phi : V_x \rightarrow V'_{Kx}$
as follows:
choose $l_i$ so that $1 \leq l_i / x < K$.
We define $\phi$ as $a'_{l_i,Kx} \circ p_{x,i}$.
We can define $\phi' : V'_x \rightarrow V'_{Kx}$ in a similar way
as $p'_{i,Kx} \circ a'_{x,l_i}$.
Because $p_{x,i}$,$p'_{i,x}$ commutes with the directed system maps,
we have that $\phi \circ \phi'$ and $\phi' \circ \phi$ are directed
system maps and hence we have an isomorphism.
\qed

\proof of Lemma \ref{lemma:fromonechaincomplex}.
We have a sequence $l_i$ satisfying $l_i / l_{i-1} < \kappa$
where $H = \frac{1}{2} r_N^2 + d_i$ in a neighbourhood of $r_N = l_i$.
Here $d_i$ is a constant.
We define $H_i$ as follows:
\begin{enumerate}
\item $H_i = H$ in the region
$r_N \leq l_i - \delta_i$ where $\delta_i$ is a very small constant.
\item $H_i = f_i(r_N)$ on a small neighbourhood of $l_i - \delta_i \leq r_N \leq l_i$
where $f_i',f_i'' \geq 0$.
We can then ensure that $H_i$ has no orbits in this neighbourhood.
\item $H_i = l_i r_N - \frac{1}{2}l_i^2 + d_i$ in the region $r_N \geq l_i$.
\item $H \geq H_i$ for all $i$ (this can be done by the last condtion
that a quadratic admissible Hamiltonian satisfies).
We also want that $H_{i-1} \leq H_i$ for all $i$.
\end{enumerate}
Let $K_s$, $s \geq 1$ be a smooth family of Hamiltonians such that:
\begin{enumerate}
\item $\frac{\partial K_s}{\partial s} \geq 0$.
\item $K_s$ is linear at infinity of slope $s$.
\item $K_{l_i} = H_i$.
\end{enumerate}
We have that $(SH_*(K_s,J))$ is a filtered directed system.

Because every orbit is contained in a region of the form
$l_{i-1} \leq r_N \leq l_i$ and
\[l_{i-1}(l_{i-1} + b) + b' - \frac{1}{2}l_i^2
\leq -\theta_N(X_H) - H \leq l_i(l_i + b') + b - \frac{1}{2} l_{i-1}^2\] 
we have constants $A,B$ so that any orbit starting on the level set $r_N = c$
must have action between $Ac^2$ and $Bc^2$.
Choose $l_i$ so that $Al_i^2 \geq \lambda^2$.
All orbits of $H$ and $K_{l_i}$ of action $\leq \lambda^2$
are contained inside the region $r_N \leq l_i$.
The maximum principle (Corollary \ref{corollary:parameterizedfloermax})
ensures that all Floer trajectories connecting these orbits must be contained
inside $r_N \leq l_i$.
Hence $SH_*^{\leq \lambda^2}(H,J) = SH_*^{\leq \lambda^2}(K_{l_i},J)$.
We have a natural map:
\[ SH_*^{\leq \lambda^2}(K_{l_i},J) \rightarrow SH_*(K_{l_i},J).\]
Composing these two maps gives us a natural map
\[ p_{\lambda,i} : SH_*^{\lambda^2}(H,J) \rightarrow SH_*(K_{l_i},J).\]
Because the filtered directed system maps for $H$ are induced by inclusions
and the filtered directed system maps for $K_s$ are induced by non-decreasing
families of Hamiltonians we have that $p_{\lambda,i}$ commutes with the filtered
directed system maps as described in the statement of Lemma \ref{lemma:discreteisomorphism}.

Now choose $\lambda$ so that $Bl_i^2 \leq \lambda^2$.
Then all orbits of $K_{l_i}$ have action $\leq \lambda^2$.
So the subcomplex of $CF_*^{\leq \lambda^2}(H,J)$ generated by orbits in the region
$r_N \leq l_i$ is isomorphic to the chain complex $CF_*(K_{l_i},J)$.
Hence there is a natural morphism
\[p'_{i,\lambda} : SH_*(K_{l_i},J) \rightarrow SH_*^{\leq \lambda^2}(H,J).\]
These morphisms also commute with the natural directed system maps.
So by Lemma \ref{lemma:discreteisomorphism},
$(SH_*^{\leq \lambda^2}(H,J))$ is isomorphic to
$(SH_*(K_s,J))$ which in turn by
Corollary \ref{cor:lineardirectedsystemhaspositiveaction}
and Lemma \ref{lemma:growthrateequalsactiongtrowthrate}
is isomorphic to $(SH_*(N,\lambda))$.
\qed

\proof of Theorem \ref{theorem:cotangentfiltereddirectedsystem}.
By Lemma \ref{lemma:fromonechaincomplex}, we get that
$(SH_*^{\leq \lambda^2}(H,J))$ is isomorphic to
$(SH_*^{\#}(T^*Q,\theta_Q,\lambda))$ for any quadratic admissible pair $(H,J)$.
Using the metric on $Q$, we have a functional $S$ defined on the loopspace
${\mathcal L} Q$ given by sending a loop $l : S^1 = \R / \Z \rightarrow Q$
to $\sqrt{\int_0^1 |l'(t)|^2 dt}$.
We have that $H_*(S^{-1}(-\infty,\lambda])$ is a filtered directed system
where the directed system maps come from the natural inclusion maps.
From \cite[Corollary 1.2]{salamonweber:loop} we have
a quadratic admissible pair $(H,J)$
with the property that there is an isomorphism
$SH_*^{\leq \lambda^2}(H,J,(\tau_Q,\omega_2)) \cong H_*(S^{-1}(-\infty,\lambda])$  for all $\lambda$.
This isomorphism commutes with the directed system maps as well.
Hence we get that $(SH_*^{\leq \lambda^2}(H,J,(\tau_Q,\omega_2)))$ is isomorphic
to $(H_*(S^{-1}(-\infty,\lambda])$ as directed systems.
Let $L_{\text{par}}^{\leq \lambda}Q$ be the set of loops $l$ of length $\leq \lambda$
such that $|l'(t)|$ is constant.
This space is homotopic to $L^\lambda Q$.
We have that $S(l)$ is equal to the length of $l$
when $l \in L_{\text{par}}^{\leq \lambda}$.
The result in \cite{Anosov:homotopiescurves} tell us that
the inclusion $L_{\text{par}}^{\leq \lambda}Q \hookrightarrow  S^{-1}(-\infty,\lambda]$
is a homotopy equivalence. Hence
$L^{\leq \lambda}$ is homotopy equivalent to $S^{-1}(-\infty,\lambda]$.
We now get that $(SH_*^{\leq \lambda^2}(H,J,(\tau_Q,\omega_2)))$ is isomorphic
to $(H_*(L^{\leq \lambda}Q))$ as directed systems.
Hence $(SH_*^{\#}(T^*Q,\theta_Q,\lambda,(\tau_Q,\omega_2)))$ is isomorphic
to $(H_*(L^\lambda))$ as directed systems.
\qed

\subsubsection{The growth rate and fundamental group}

Let $G$ be a finitely generated group, and let $A := \{g_1, \cdots, g_k\}$
be a set of generators for $G$.
Let $G^{\text{cong}}$ be the set of elements of $G$
modulo conjugacy classes.
Let $r_i$ be the number of elements of $G^{\text{cong}}$
which can be expressed as a product of at most $i$
generators.
We define $\Gamma^\text{cong}(G)$ to be
$\varlimsup_i \frac{\log{r_i}}{\log{i}}$.
This definition is independent of the choice of generators.
The reason is as follows:
Suppose $g'_1, \cdots, g'_{k'}$ is another set of generators,
then all we need to do is show that the growth rate
associated to $A$ is the same as the
growth rate associated to the union of the generators:
$B := \{g_1, \cdots, g_k,g'_1, \cdots, g'_{k'}\}$.
Let $r'_i$ be the number of elements of $G^{\text{cong}}$
which can be expressed as a product of at most $i$ elements of $B$.
Then $r_i \leq r'_i$ because $A \subset B$.
We have that there exists a $K \in \N$ such that
each element of $B$ can be expressed as a product of at most $K$
elements of $A$.
Therefore $r'_i \leq K r_i$ for all $i$
(where $K$ is independent of $i$).
Hence the growth rates are the same.

\begin{lemma} \label{lemma:growthratefundamentalgroup}
Let $Q$ be a compact oriented manifold,
then \[\Gamma(Q) \geq \Gamma^\text{cong}(\pi_1(Q)).\]
\end{lemma}
\proof
We need to show the following fact:
there exists a constant $P$ depending only on $Q$
such that for every pair of free loops $\gamma_1,\gamma_2$ in $Q$,
there exists a free loop $\gamma$ and
elements $k_1$ and $k_2$ of $\pi_1(Q)$ such that
$\gamma_i$ represents the conjugacy class $[k_i]$,
and $\gamma$ represents the conjugacy class $[k_1.k_2]$ with
$l(\gamma) \leq l(\gamma_1) + l(\gamma_2) + P$.
Here $l$ denotes the length of a loop with respect to the metric $g$.
This is done as follows:
Choose a constant $P$ such that for every pair
of points $q_1,q_2 \in Q$, there exists a path joining
them of length $\leq \frac{1}{2}P$.
Also choose a basepoint $a \in Q$.
Join $\gamma_1$ and $\gamma_2$ to $a$ using
paths of length $\leq \frac{1}{2}P$.
This gives us elements $k_1,k_2 \in \pi_1(Q,a)$.
The composition of such loops gives us a loop
of length $\leq l(\gamma_1) + l(\gamma_2) + P$.

Let $g_1, \cdots, g_k$ be generators for $\pi_1(Q,a)$.
Let $C$ be a constant greater than $l(g_i)$ for each $i$.
For $C \geq 1$ we have that $ \text{rank}(H_0^{\leq (C+P)c}({\mathcal L}Q))$
is an upper bound for the number of loops of word length $\leq c$.
Also the growth rate associated to $H_0^{\leq (C+P)c}({\mathcal L}Q)$
is the same as the growth rate of  $H_0^{\leq c}({\mathcal L}Q)$.
Hence $\Gamma(T^*Q,(\tau_Q,\omega_2))$ is an upper bound for
$\Gamma^\text{cong}(\pi_1(Q))$.
\qed

\begin{lemma} \label{lemma:triplefreeproduct}
Let $G$ be the free product of $3$ non-trivial groups
then \[\Gamma^{\text{cong}}(G) = \infty.\]
\end{lemma}
\proof of Lemma \ref{lemma:triplefreeproduct}.
Let $G = A \star B \star C$ where $A,B,C$ are non-trivial.
Let $a,b,c$ be non-trivial elements in $A,B,C$.
Choose a subset $I \subset \{1,\cdots,k\}$.
Let $q_I(i)$ be a function which is $1$ if $i \in I$
and $0$ otherwise.
Let $a_I := b c \prod_{i=1}^k (a^{q_I(i)}) (b c)^{1-q_I(i)})$.

%
Every element of $G$ can be written uniquely
in the form $f := \prod_{i=0}^l h_i$
where $h_i \in A$ or $B$ or $C$
and $h_{i-1}$ is not in the same group $A,B,C$ as $h_i$.
There is an $l' < \frac{l}{2}$
such that $h_i = h_{l - i}^{-1}$
for all $i \leq l'$.
The element $\prod_{i=l'+1}^{l-l'-1} h_i$
is called the {\it conjugation interior} of $f$.
We can conjugate $f$ by some element
$g$ so that the conjugation interior
is of the form $\prod_{i=l'+1}^{l-l'-1} h_i$
where $h_{l'+1}$ is not in the same
group $A$ or $B$ or $C$ as $h_{l-l'-1}$.
This is called a {\it standard conjugation interior}.
If we have some element
$\prod_{j=1}^{l'} h'_j$ as described above then a {\it rotation}
of $\prod_{j=1}^{l'} h'_j$ is the operation where we replace
it with $h'_l \prod_{j=1}^{l'-1} h'_j$ or we do the reverse.

If I conjugate $f$ by any element
$g$ then the standard conjugation interior can only change
by a sequence of rotations.

The conjugation interior of
$a_I$ is equal to $a_I$ and it is standard.
Hence if $a_{I'}$ is conjugate to $a_I$
then $I = I' +j$ where we view $\{1,\cdots,k\}$
as the cyclic group with $k$ elements.
We say that $I'$ is a rotation of $I$ if $I = I' + j$
for some $j \in \Z$.
Also the word length of all these $a_I$'s
with respect to the generators $a,b,c$
are all the same. This word length is between $2+ k$ and $2 + 2k$.
Hence the number of conjugacy classes
of elements of $G$ of word length between $2 + k$ and $2 + 2k$
in $a,b,c$ is at least the number of subsets
$I \subset \{1,\cdots k\}$ modulo rotation  which is
at least $\frac{2^k}{k}$.
Hence the growth rate of $G$ is bounded
below by the growth rate of $\frac{2^k}{k}$
and hence is infinite.
\qed

\section{Compactifications of algebraic varieties}
\label{secion:compactificationsofalgebraicvarieties}

\subsection{Making the divisors orthogonal}

Let $(M,\omega)$ be a symplectic manifold of dimension $2n$.
Let $S_1,\cdots,S_k$ be symplectic submanifolds of real co-dimension $2$.
For each $I \subset \{1,\cdots,k\}$ we define $S_I$ to be $\cap_{i \in I} S_I$.
We say that $S_1,\cdots,S_k$ are {\it symplectically intersecting}
if they intersect transversally and $S_I$ is a symplectic submanifold
for each $I \subset \{1,\cdots,k\}$.
For any symplectic submanifold $S \subset M$, we define
its $\omega${\it-orthogonal bundle} $NS$ to be the vector subbundle of $TM|_S$
given by vectors $u$ satisfying $\omega(u,v) = 0$ for each $v \in TS$.
\begin{defn}
We say that $S_1,\cdots,S_k$ {\it intersect positively}
if
\begin{enumerate}
\item They are symplectically intersecting.
\item Let $I \subset \{1,\cdots,k\}$ be a disjoint union $I_1 \sqcup I_2$,
$N_1$ the $\omega$-orthogonal bundle of $S_I$ inside $S_{I_1}$
and $N_2$ the $\omega$-orthogonal bundle of $S_I$ inside $S_{I_2}$.
The bundle $TS_I \oplus N_1 \oplus N_2$ is isomorphic to $TM|_S$.
Each bundle $TS_I,N_1$ and $N_2$ has an orientation induced by $\omega$
and hence their direct sum does.
We require that the natural orientation on this direct sum
matches the orientation induced by $\omega^n$ on $TM|_S$.
\end{enumerate}
\end{defn}

\begin{defn} \label{defn:orthogonaldivisors}
Let $Q_1,\cdots,Q_l$ be a collection of symplectic submanifolds
of any dimension. Let $U$ be any subset of $M$.
We say that $Q_1,\cdots,Q_k$ are {\it orthogonal} 
along $U$
if for each $i,j \in \{1,\cdots,k\}$ ($i \neq j$) and $x \in U \cap Q_i \cap Q_j$,
the $\omega$ orthogonal normal bundle to $Q_i$ at $x$ is contained in $TS_j$.
We just say they are orthogonal if they are orthogonal along $M$.
\end{defn}

The aim of this section is to deform positively intersecting submanifolds
so that they become orthogonal.
\begin{lemma} \label{lemma:deformingdivisors}
Let $S_1,\cdots,S_k$ be any finite set of positively intersecting symplectic submanifolds.
There is a smooth family of positively intersecting symplectic submanifolds
$S^t_i$ such that
\begin{enumerate}
\item $S^0_i = S_i$ for all $i \in \{1,\cdots,k\}$.
\item All the $S^1_i$ intersect orthogonally.
\item $S^t_k = S_k$.
\end{enumerate}
\end{lemma}

We need some preliminary lemmas and definitions before we prove the above lemma.

\begin{lemma} \label{lemma:deformingpositivelyiintersecting}
Suppose we have $k$ smooth families of symplectic submanifolds
$S_1^t,\cdots,S_k^t$ parameterized by $t \in [0,1]$
and that $S_1^t,\cdots,S_k^t$ are symplectically intersecting
for each $t$.
Suppose also that $S_1^0,\cdots,S_k^0$ are positively intersecting.
Then for each $t$, $S_1^t,\cdots,S_k^t$ are also positively intersecting.
\end{lemma}
\proof of Lemma \ref{lemma:deformingpositivelyiintersecting}.
For $I \subset \left\{ 1,\cdots,k \right\}$ let $S_I^t$ be equal to
$\cap_{i \in I} S_i^t$. If $I$ is the disjoint union of $I_1$ and $I_2$
then we define $N_i^t$ to be equal to be the $\omega$-orthogonal bundle of $S_I^t$
inside $S_{I_i}^t$.
We have a smooth family of bundles $TS_I^t \oplus N_1^t \oplus N_2^t$,
and the orientation of $TS_I^0 \oplus N_1^0 \oplus N_2^0$ agrees with that of
$TM|_{S_I^0}$. Hence because we have a smooth family of bundles,
the orientation of $TS_I^t \oplus N_1^t \oplus N_2^t$ agrees with the orientation of $TM|_{S_I^t}$.
\qed

From now on, $\R^{2n}$ is the standard symplectic vector space with coordinates
$x_1,y_1,\cdots,x_n,y_n$ and symplectic form
$\omega_{\text{std}} = \sum_{i=1}^n dx_i \wedge dy_i$.
\begin{lemma} \label{lemma:positivelyintersectingasinglehypersurface}
Let $S_1,\cdots,S_k$ be codimension $2$ symplectic vector subspaces of $\R^{2n}$
such that $S_k = \left\{ x_n,y_n = 0 \right\}$ and where
$S_1,\cdots,S_k$ are positively intersecting.
Then for all $\mu \geq 0$, $S_1,\cdots,S_k$ are positively intersecting with respect to the new
symplectic form $\omega_\mu := \omega_{\text{std}} + \mu dx_n \wedge dy_n$.
\end{lemma}
\proof of Lemma \ref{lemma:positivelyintersectingasinglehypersurface}.
Let $I \subset \left\{ 1,\cdots,k \right\}$.
If $k \in I$ then $S_I$ is a subset of $S_k$,
hence $\omega_\mu$ restricted to $S_I$
is equal to $\omega_{\text{std}}$ restricted to $S_I$ and hence $S_I$
is symplectic with respect to $\omega_\mu$ for all $\mu \geq 0$.
Now suppose that $k \notin I$.
Then $S_I$ is transverse to $S_k$.
Let $F$ be the $\omega_{\text{std}}$ orthogonal subspace of $S_I$
to $S_I \cap S_k \subset S_I$. This is a two dimensional symplectic subspace.
Let $N$ be the $\omega_{\text{std}}$ orthogonal bundle to $S_I \cap S_k$
inside $S_k$. We know that $N \oplus (S_I \cap S_k) \oplus F$ has the same orientation
as $\R^{2n}$ because $S_1,\cdots,S_k$ are positively intersecting.
Also $N \oplus (S_I \cap S_k)$ has the same orientation as $S_k$
because $N$ and $S_I$ are $\omega_{\text{std}}$ orthogonal.
Hence the orientation on $F$ is the same as the orientation induced
by restricting $dx_n \wedge dy_n$ to $F$.
Hence $\omega_\mu$ restricted to $F$ is a volume form on $F$ for all $\mu \geq 0$.
Also $F$ is still orthogonal to $S_I \cap S_k$ inside $S_I$ with respect to $\omega_\mu$
because $\omega_\mu(V,W) = \omega_{\text{std}}(V,W)$ for all vectors $V,W$ where $V$ is tangent to $S_k$.
Hence $\omega_\mu$ restricted to $S_I$ is still a symplectic form.

So $S_1,\cdots,S_k$ is symplectically intersecting with respect to $\omega_\mu$
for all $\mu$.
Because $S_1,\cdots,S_k$ are positively intersecting with respect to $\omega_0 = \omega_{\text{std}}$
we have by Lemma \ref{lemma:deformingdivisors} that
$S_1,\cdots,S_k$ are all positively intersecting with respect to $\omega_\mu$ for all $\mu \geq 0$.
\qed

\begin{lemma} \label{lemma:makingtingssymplectic}
Let $S_1,\cdots,S_k$ be transversally intersecting codimension $2$ vector subspaces of $\R^{2n}$
such that $S_k = \left\{ x_n,y_n = 0 \right\}$.
Suppose that $S_1 \cap S_k,\cdots,S_{k-1} \cap S_k$
are symplectically intersecting inside $S_k$, then for large enough $\mu \geq 0$
we have that $S_1,\cdots,S_k$ are symplectically intersecting with respect to the symplectic form
$\omega_\mu := \omega_{\text{std}} + \mu dx_n \wedge dy_n$.
\end{lemma}
\proof of Lemma \ref{lemma:makingtingssymplectic}.
Let $I \subset \left\{ 1,\cdots,k \right\}$.
If $k \in I$ then $S_I$ is a symplectic manifold with respect to $\omega_\mu$.
From now on we will assume that $k \notin I$.
We have that $S_I$ intersects $S_k$ transversally.
Let $F \subset S_I$ be the vector subspace of $S_I$ which consists
of vectors which are $\omega_{\text{std}}$ orthogonal to $S_I \cap S_k$.
First of all $F$ is of dimension at least $2$ because the set of vectors $\omega_{std}$
orthogonal to $S_I \cap S_k$ inside $\R^{2n}$ is $2(|I|+1)$ dimensional
and $S_I$ is $2(n-|I|)$ dimensional.
Here $|I|$ means the number of elements in $I$.
Also $F$ has dimension at most $2$ because it must be orthogonal to $S_I \cap S_k$
inside $S_I$. Hence $S_I$ is equal to $(S_I \cap S_k) \oplus F$ and
$\omega_{\text{std}}$ restricted to $S_I$ splits up under this direct sum
as $\omega_1 \oplus \omega_2$. Because $S_I \cap S_k$ is symplectic by assumption we have that
$\omega_1$ is a symplectic form. Also $dx_n \wedge dy_n$ restricted to $F$ is a non-degenerate
$2$-form because $F$ is $2$-dimensional and transverse to $S_k$.
So for large enough $\mu$, $\omega_2 + \mu dx_n \wedge dy_n|_F$ is a sympectic form on $F$.
Hence $\omega_{\text{std}} + \mu dx_n \wedge dy_n$ is a symplectic form on $S_I$
for large enough $\mu$.
Hence we have shown that $S_I$ is symplectic for all $I$ with respsect to
$\omega_{\text{std}} + \mu dx_n \wedge dy_n$ for large enough $\mu \geq 0$.
\qed

We have a parameterized version of this lemma that will be needed where
we have a continuous family $S_1^q,\cdots,S_{k-1}^q$ parameterized by $q \in Q$
where $Q$ is some compact topological space.
The result is that $S_1^q,\cdots,S_{k-1}^q,S_k$ are symplectically intersecting
with respect to $\omega_\mu$ for $\mu$ sufficiently large.

\begin{lemma} \label{lemma:preliminarylinearalgebra}
Let $S_k$ be a co-dimension $2$ symplectic
vector subspace of $\R^{2n}$.
Let $A_1,\cdots,A_{k-1}$ be symplectically intersecting symplectic vector subspaces
of $S_k$.
Let ${\mathcal B}$ be the space of $(k-1)$-tuples of $2n-2$ dimensional symplectic vector subspaces
$(S_1,\cdots,S_{k-1})$ of $\R^{2n}$
such that $S_1,\cdots,S_k$ are positively intersecting
and such that $S_i \cap S_k = A_i$.
If ${\mathcal B}$ is non-trivial then it
deformation retracts onto the point
$(A_1 + S_k^\perp,\cdots,A_{k-1} + S_k^\perp)$
where $S_k^\perp$ is the symplectic orthogonal subspace to $S_k$.

In particular we also get that 
$(A_1 + S_k^\perp,\cdots,A_{k-1} + S_k^\perp)$
are also positively intersecting and hence
$A_1,\cdots,A_{k-1}$ are positively intersecting inside $S_k$.
\end{lemma}
\proof of Lemma \ref{lemma:preliminarylinearalgebra}.
By a linear symplectic change of coordinates we can assume that $S_k = \left\{ x_n,y_n=0 \right\}$.
We will show that ${\mathcal B}$ is weakly contractible and hence contractible.
Let $S_1,\cdots,S_{k-1}$ be a point in ${\mathcal B}$.
The subspace $S_{\bracket{1,\cdots,k-1}}$ is at least two dimensional and transverse to $S_k$.
Let $W \subset S_{\bracket{1,\cdots,k-1}}$ be a two dimensional symplectic vector subspace transverse to
$S_k$. Such a subspace exists for the following reason:
We have that $S_1,\cdots,S_k$ intersect transversally and because these are codimension $2$
vector subspaces we have that $k \leq n$. Hence $S_{\bracket{1,\cdots,k-1}}$ has codimension at most $2(n-1)$
which means that it is at least $2$ dimensional. This space is also transverse to $S_k$.
Hence we can find such a two dimensional subspace $W$.
The vector subspace $W$ is also contained in $S_i$ for each $1 \leq i \leq k-1$.
Hence $S_i$ is a direct sum of vector spaces $(S_i \cap S_k) \oplus W$
for each $1 \leq i \leq k-1$.
Let $S_1^q,\cdots,S_{k-1}^q$ be a family of points in ${\mathcal B}$ continuously parameterized by
points $q$ in some sphere $S^m$.
Let $\phi^q : \R^2 \rightarrow \cap_{i=1}^{k-1} S_i^q$ be a continuous family of linear embeddings
parameterized by $q$ such that the image of $\phi^q$ is transverse to $S_k$.
Now choose a continuous family of maps $\phi_t^q : \R^2 \rightarrow \R^{2n}$
parameterized by $(q,t) \in S^m \times [0,1]$ such that they are linear embeddings
transverse to $S_k$, $\phi^q_0 = \phi^q$ and such that
$\phi^q_1$ is a linear isomorphism to $S_k^\perp$.
We now have a family of transversally intersecting (not necessarily symplectic)
vector subspaces $S_i^{q,t} := A_i \oplus \text{image} (\phi^q_t)$.
We have that $S_i^{q,0} = S_i^q$ and
$S_i^{q,1} = A_i \oplus S_k^\perp$.
By (a parameterized version of) Lemma \ref{lemma:makingtingssymplectic},
we have that $S_i^{q,t}$ is symplectically intersecting with respect to
$\omega_\mu := \omega_{\text{std}} + \mu dx_n \wedge dy_n$
where $\mu \geq 0$ is very large
(here we really need that $S_i^{q,t}$ is parameterized by a compact family
and this is why we need to prove weak contractibility first).
Let $\Phi_\mu$ be a linear automorphism sending
$(x_1,y_1,\cdots,x_{n-1},y_{n-1},x_n,y_n)$
to 
$(x_1,y_1,\cdots,x_{n-1},y_{n-1},x_n,\frac{1}{1+\mu}y_n)$.
Then $\Phi^*_\mu \omega_\mu = \omega_{\text{std}}$.
This automorphism also preserves $S_k$ and $S_k^\perp$,
so $\Phi^*_\mu S_i^{q,1} = S_i^{q,1}$.

By Lemma \ref{lemma:positivelyintersectingasinglehypersurface}
we have that $\Phi^*_{t\mu} S_i^{q,0}$ is positively intersecting for all $t \in [0,1]$.
We concatenate the isotopies $\Phi^*_{t\mu} S_i^{q,0}$
and $\Phi^*_\mu S_i^{q,t}$ giving us a new family of vector subspaces
$(\sigma_1^{q,t},\cdots,\sigma_{k-1}^{q,t})$ defined as follows:
\begin{enumerate}
\item For $t \in [0,\frac{1}{2}]$,
$\sigma_i^{q,t} = \Phi^*_{2t\mu} S_i^{q,0}.$
\item For $t \in [\frac{1}{2},1]$,
$\sigma_i^{q,t} = \Phi^*_\mu S_i^{q,2t-1}.$
\end{enumerate}
By Lemma \ref{lemma:deformingpositivelyiintersecting} we get that
$(\sigma_1^{q,t},\cdots,\sigma_{k-1}^{q,t},S_k)$ are positively intersecting
because $(\sigma_1^{q,0},\cdots,\sigma_{k-1}^{q,0},S_k)$ are.
Hence $(\sigma_1^{q,t},\cdots,\sigma_{k-1}^{q,t})$ are points in
${\mathcal B}$ parameterized by $q \in S^m$
starting at $S_1^q,\cdots,S_{k-1}^q$ and ending
at $A_1 \oplus S_k^\perp,\cdots,A_k \oplus S_k^\perp$.
Hence ${\mathcal B}$ is weakly  contractible and hence contractible.
\qed

\begin{lemma} \label{lemma:linearorthogonal}
Let ${\mathcal S}_k$ be the space of $k$-tuples of
positively intersecting $2n-2$ dimensional symplectic vectoor subspaces
$S_1,\cdots,S_k$ of $\R^{2n}$.Then ${\mathcal S}_k$
deformation retracts onto the space of positively
intersecting $\omega_{\text{std}}$ orthogonal
vector subspaces of $\R^{2n}$.
\end{lemma}
\proof of Lemma \ref{lemma:linearorthogonal}.
Throughout the proof of the Lemma, we use the following fact:
suppose that I have a fibration $p$ whose
fibers deformation retract to the fibers of a subfibration $p'$.
Then the total space of $p$ deformation
retracts to the total space of $p'$.

We proceed by induction on dimension.
Suppose this is true for all $\R^{2l}$
with $l < n$ and consider $\R^{2n}$.
We define ${\mathcal S}^{S_k}_{k-1}$
as the space of positively intersecting
$2n-4$ dimensional symplectic vector subspaces
$A_1,\cdots,A_{k-1}$ of $S_k$.
Then this deformation retracts to
the space of orthogonally intersecting
subspaces by our induction hypothesis.

We have a fibration:
\[P : {\mathcal S}_k \twoheadrightarrow {\mathcal S}_1\]
sending $(S_1,\cdots,S_k)$ to $S_k$.
We have an inclusion ${\mathcal S}^{S_k}_{k-1} \hookrightarrow P^{-1}(S_k)$
sending $A_1,\cdots,A_{k-1}$ to $(A_1 + S_k^\perp,\cdots,A_{k-1} + S_k^\perp)$
where $S_k^\perp$ is the symplectic orthogonal subspace to $S_k$.
So from now on we view 
${\mathcal S}^{S_k}_{k-1}$ as a subspace of $P^{-1}(S_k)$
and the union $\cup_{S_k \in {\mathcal S}_1} {\mathcal S}^{S_k}_{k-1}$
as a subfibration of $P$.
The fiber $P^{-1}(S_k)$ is also a fibration
$Q : P^{-1}(S_k) \twoheadrightarrow {\mathcal S}^{S_k}_{k-1}$
where $Q(S_1,\cdots,S_k) := (S_1 \cap S_k,\cdots S_{k-1} \cap S_k)$.
The fiber of $Q$ over a point $(A_1,\cdots,A_{k-1})$ consists of subspaces
$B_1,\cdots B_{k-1}$ positively intersecting $S_k$ such that
$B_i \cap S_k = A_i$.
This deformation retracts to the point $(A_1 + S_k^\perp,\cdots,A_{k-1} + S_k^\perp)$
by Lemma \ref{lemma:preliminarylinearalgebra}.
Hence $P^{-1}(S_k)$ deformation retracts to
the space ${\mathcal S}^{S_k}_{k-1}$ which by our induction
hypothesis deformation retracts to the space of orthogonally
intersecting subspaces of $\R^{2n-2}$.
This implies that ${\mathcal S}_k$ deformation retracts
to the space of orthogonally intersecting subspaces.
This proves our Lemma.
\qed

\begin{lemma} \label{lemma:neighbourhoodofdivisor}
Let $S \subset M$ be a symplectic submanifold
of $M$. Suppose that the $\omega$-orthogonal
bundle $\pi_S : NS \twoheadrightarrow S$ has structure group $G \subset U(n)$.
Then there is a neighbourhood $US$ and a projection
$p : US \twoheadrightarrow S$ whose fiber is symplectomorphic
to the ball $B_\epsilon$ of radius $\epsilon$
and such that its structure group is $G$ and whose
fibers are orthogonal to $S$.

There is also a vector field $L$
on $US$ whose flow $\phi_t$ is well defined for all negative $t$
and such that $e^t (\phi_{-t})^* \omega_{US} = \omega_{US} + (e^t -1)\pi^* \omega_S$
where $\omega_{US}$ and $\omega_S$ are the symplectic forms
on $US$ and $S$ respectively (here $S$ is identified with the zero section).
Inside  the fiber $B_\epsilon \subset \C^n$,
this vector field is equal to $\sum_i \frac{r_i}{2} \frac{\partial}{\partial r_i}$
where $(r_i,d\theta_i)$ are polar coordinates for the $i$th factor
of $\C$ in $\C^n$.
\end{lemma}
\proof of Lemma \ref{lemma:neighbourhoodofdivisor}.
This proof is very similar to the proof of \cite[Theorem
6.3]{McduffSalamon:sympbook}.
Let $(V_\alpha)_{\alpha \in A}$ be a finite covering
of $S$ with trivializations $V_\alpha \times \C^n$
of the symplectic fibration.
Choose a partition of unity $\rho_\alpha : V_\alpha \rightarrow \R$
subordinate to this cover.
Let $\pi_{\C^n}$ be the natural projection from
 $V_\alpha \times \C^n$ to $\C^n$.
We define $\sigma_\alpha := \sum_i r_i^2 d\theta_i$
where $(r_i,\theta_i)$ are polar coordinates for
the $i$th $\C$ factor of the fiber $\C^n$.
We write $\tau := \sum_\alpha d( (\rho_\alpha \circ \pi_S)  \sigma_\alpha)$
which is equal to
\[\tau = \sum_\alpha d(\rho_\alpha \circ \pi_S) \wedge \sigma_\alpha
+ \rho_\alpha \circ \pi_S d\sigma_\alpha.\]
Because $\sigma_\alpha$ restricted to the zero section is zero
and $TS$ is in the kernel of $d\sigma_\alpha$,
we have that
$\omega_{US} := \pi_S^* \omega|_S + \tau$ is a symplectic
form on a small neighbourhood of the zero section.
This is because $\|d(\rho_\alpha \circ \pi_S) \wedge \sigma_\alpha\|$
is small relative to $\|\rho_\alpha \circ \pi_S d\sigma_\alpha\|$
and $\|\pi_S^* \omega|_S \|$ near the zero section.
Let $NS_\epsilon$ be the open subset consisting
of vectors of modulus less than $\epsilon$
where $\epsilon$ is small enough so that $\omega_{US}$
is still a symplectic form on $NS_\epsilon$.
We have by a Moser theorem \cite[Theorem 3.3]{McduffSalamon:sympbook}
that for some $\epsilon$ small enough,
a neighbourhood of $S$ is symplectomorphic
to $NS_\epsilon$.
This has structure group $G$ and fibers $B_\epsilon$.

The fibers are orthogonal to $S$ because if we
have any vector $v$ tangent to $S$ and another one
$w$  tangent to the fiber, then
$\pi_S^* \omega|_S(v,w) = 0$
because $w$ is in the kernel of this symplectic form.
Also $d(\rho_\alpha \circ \pi_S) \wedge \sigma_\alpha(v,w) = 0$ because $\sigma_\alpha$ vanishes on $S$.
Finally $d\sigma_\alpha(v,w)$ vanishes because
$TS$ is in the kernel of $d\sigma_\alpha$.
Hence $\omega_{US}(v,w) = 0$ and hence the fibers are orthogonal to $S$.

Let $L$ be the vector field on $NS_\epsilon$ which is tangent to the fibers
of $\pi_S$ and equal to $\sum_i \frac{r_i}{2} \frac{\partial}{\partial r_i}$
on the fibers.
This is well defined because $\sum_i \frac{r_i}{2} \frac{\partial}{\partial r_i}$
is invariant under the $U(n)$ action on $\C^n$.
We have that the Lie derivative of
$\rho_\alpha \circ \pi_S . \sigma_\alpha$ with respect to $L$
is $\rho_\alpha \circ \pi_S . \sigma_\alpha$ (because the flow of $L$
does not change $\rho_\alpha \circ \pi_S$).
Hence $(\phi_t)^* \tau = e^t \tau$ where $\phi_t$ is the flow of $L$.
Hence
\[e^t (\phi_{-t})^* \omega_{US} = (e^t - 1) \omega|_S + \omega_{US}.\]
\qed

\begin{lemma} \label{lemma:orthogonallineardivisor}
Let $p : US \twoheadrightarrow S$ be a fibration
described in Lemma \ref{lemma:neighbourhoodofdivisor}
and let $N \subset N'$ be open subsets of the zero section
such that the closure of $N$ is contained in $N'$
($S$ here is compact and it may or may not have a boundary).
Let
$S_1,\cdots,S_k$ be positively intersecting co-dimension $2$ submanifolds
such that
\begin{enumerate}
\item
$\cap_i S_i$ is the zero section.
\item
$S_i$ intersects the fibers transversely.
\item
$S_i$ intersect orthogonally on $p^{-1}(N')$
\item
$S_i \cap p^{-1}(q)$ is a symplectic vector subspace of $B_\epsilon \subset \C^n$ for each $q \in S$.
\end{enumerate}
Then there is a family of manifolds
$S_i^t$ ($t \in [0,1]$)
\begin{enumerate}
\item $S_i^t$ are positively intersecting symplectic manifolds
with respect to the symplectic form $\omega_{US} + T p^* \omega_S$ for sufficiently
large $T > 0$.
\item $\cap_i S_i^t$ is the zero section.
\item $S_i^0 = S_i$.
\item $S_i^t \cap p^{-1}(N) = S_i \cap p^{-1}(N)$.
\item $S_i^1$ are all orthogonal along the zero section.
\item $S_i^t = S_i$ outside some neighbourhood of the zero section.
\end{enumerate}
\end{lemma}
\proof of Lemma \ref{lemma:orthogonallineardivisor}.

By Lemma \ref{lemma:linearorthogonal},
there is a family of submanifolds
$W_i^t$ of $US$ transverse to all the fibers such that
\begin{enumerate}
\item $W_i^t \cap p^{-1}(q)$ are positively intersecting symplectic vector subspaces of $p^{-1}(q)$
for all $q \in S$.
\item $W_i^t$ intersect each fiber of $p$ transversely.
\item $W_i^0 = S_i$.
\item \label{item:restrictiontoN}
$W_i^t = S_i$ on $p^{-1}(N)$
\item $W_i^1 \cap p^{-1}(q)$ are orthogonal inside $p^{-1}(q)$
for all $q \in S_I$.
\end{enumerate}
Note that $W_i^t$ become positively intersecting symplectic
submanifolds with respect to the symplectic form
$\omega_{US} + T p^* \omega_S$ for all sufficiently large $T > 0$.

For each $s \in [0,1]$,
we can construct (using a bump function)
a smooth family of co-dimension $2$ submanifolds
$U_i^{s,t}$ ($t \in [0,1]$) such that $U_i^{s,t} = W_i^s$
outside some closed subset of $US$ containing the zero section,
$U_i^{s,t} = W_i^t$ near the zero section
and $U_i^{s,s} = W_i^s$.
We can also assume that $U_i^{s,t} = W_i^s$ on $p^{-1}(N)$
because $W_i^s = W_i^t$ in this region.
For each $s$, there is a constant $\delta_s > 0$
such that for all $t \in [s-\delta_s,s+\delta_s]$,
$U_1^{s,t},\cdots,U_k^{s,t}$
are positively intersecting symplectic submanifolds
with respect to
$\omega_{US} + T p^* \omega_S$ for all sufficiently large $T > 0$.
We also assume that $U_i^{s,t}$ are transverse to all the fibers of $p$
for all $t \in [s-\delta_s,s+\delta_s]$.

Let $L$ be the vector field described in the statement of
Lemma \ref{lemma:neighbourhoodofdivisor}
and let $\phi_t$ be its flow.
For any $S \geq 0$ we have that
$\phi_{-S}(U_i^{s,t})$ is still a symplectic submanifold with respect to the symplectic
form $\omega_{US} + T p^* \omega_S$
because
\[e^S(\phi_{-S})^*(\omega_{US} +  T p^* \omega_S) = \omega_{US} + (T + 2e^S - 1)p^* \omega_S.\]
From now on we define
$\phi_{-T}(U_i^{s,t})$
to be equal to
$\phi_{-T}(U_i^{s,t})$ inside $\phi_{-T}(US)$
and equal to $U_i^{s,t}$ outside $\phi_{-T}(US)$.
This is a well defined symplectic manifold because
$U_i^{s,t} \cap p^{-1}(q)$
is linear outside some closed subset of $US$
containing the zero section for all $q \in S$ and the flow of $L$
is tangent to any submanifold with this property.
We also have that $\phi_{-T}(U_i^{s,t})$
are all positively intersecting symplectic submanifolds.
Also because $L$ is tangent to the fibers
of $p$, we have that $\phi_{-T}(U_i^{s,t})$
is transverse to all the fibers.

By compactness, there exists a finite sequence
$s_1 < s_2 < \cdots < s_l$ such that
the intervals $(s_i,s_i+\delta_{s_i})$
cover $(0,1]$.
This is because we can choose $s_1 = 0$
and then cover $[\delta_{s_1},1]$
with the open intervals $(s_i,s_i + \delta_{s_i})$
by compactness.
We assume that $s_l = 1$ and
$s_{i+1} \in (s_i,s_i + \delta_{s_i})$.

The submanifold $\phi_{-T}(U_i^{s,t})$
is still equal to $W_i^s$ outside some closed subset of $US$ containing the zero section
and equal to $W_i^t$ near the zero section.
Suppose we have found a family
$S_i^t$ ($t \in [0,s_{m-1}]$)
such that $S_i^t = S_i$ outside some relatively compact open set $O$
and such that $S_i^{s_{m-1}} = W_i^{s_{m-1}}$
near $0$.
We also suppose that $S_i^t = S_i$ on $p^{-1}(N)$.
We now consider $s_m$.
Let $O_m \subset O$ be an open set containing the zero section
such that $S_i^{s_{m-1}} = W_i^{s_{m-1}}$ on a neighbourhood
of the closure of $O_m$
for all $i$.
For $T$ large enough, we have
$G_t := \phi_{-T}(U_i^{s_{m-1},t})$
is equal to $W_i^{s_{m-1}}$
inside a neighbourhood of the closure of $US \setminus O_m$.
Note that $G_t$ is well defined
for $t \in [s_{m-1},s_m]$
because $(s_{m-1} - \delta_{s_{m-1}},s_{m-1} + \delta_{s_{m-1}})$
contains $[s_{m-1},s_m]$.
We define $S_i^t$ for $t \in [s_{m-1},s_m]$
to be equal to $S_i^{s_{m-1}}$ outside $O_m$
and equal to $G_t$ inside $O_m$.
This is a manifold because
these manifolds agree on the
boundary of $O_m$ and just outside $O_m$.
Also we have that $S_i^{s_m}$
is equal to $W_i^{s_m}$ near
$0$ because $U_i^{s_{m-1},s_m}$
does.
We have that $S_i^t = S_i$ on $p^{-1}(N)$
because $U_i^{s,t}$ has this property and
the flowing along $L$ preserves this property.

Hence by induction we have constructed a family $S_i^t$
equal to $S_i$ outside $O$ and
such that $S_i^{s_l} = S_i^1$ is equal
to $W_i^1$ near the zero section.
Because
$W_i^1$ are orthogonal at the zero section,
we get that $S_i^1$ are also orthogonal at the zero section.
Finally $S_i^t = S_i$ on $p^{-1}(N)$.
This proves the Lemma.
\qed

\begin{lemma} \label{lemma:neighbourhooduniqueness}
Let $S_1,\cdots,S_k$ be positively intersecting
symplectic manifolds in $M$
and $S'_1,\cdots,S'_k$ in $M'$.
Let $NS_i,NS'_i$ be their respective symplectic normal bundles.
Suppose that there is a series of symplectomorphisms
$\phi_i : S_i \rightarrow S'_i$
and symplectic bundle isomorphisms $N\phi_i : NS_i \rightarrow NS'_i$
covering $\phi_i$ such that
$\phi_i,\phi_j$ and $N\phi_i,N\phi_j$
agree when restricted to $S_i \cap S_j$ for all $i,j$.

Then there is a symplectomorphism $\Phi$ from
a neighbourhood of $\cup_i S_i$ to a neighbourhood
of $\cup_i S'_i$ such that $\Phi|_{S_i} = \phi_i$
and $D\Phi|_{NS_i} = N\phi_i$.

Suppose in addition that a neighbourhood $NS_I$ of $S_I$
is identical to a neighbourhood of $NS'_I$ and that for $i \in I$, $\phi_i$ and $N\phi_i$
restricted to $S_I$ is the identity map for all $i$. Then we may assume that
$\Phi$ is Hamiltonian near $S_I$ and $C^1$ small.
\end{lemma}
\proof of Lemma \ref{lemma:neighbourhooduniqueness}.
We will first construct a diffeomorphism
$\Psi$ from a neighbourhood of $\cup_i S_i$ to a neighbourhood
of $\cup_i S'_i$ such that $\Psi|_{S_i} = \phi_i$
and $D\Psi|_{NS_i} = N\phi_i$ (i.e.
every vector $v$ in $N_xS_i$ is sent to a vector in $N_{\phi_i(x)}S_j$
equal to $N\phi_i(v)$).
We do this by induction:
Suppose that we have found such a symplectomorphism
$\Psi'$ on a neighbourhood of $\cup_{i=1}^{m-1} S_i$
and consider $\cup_{i=1}^m S_i$.
By using an exponential map we have
a diffeomorphism $\tilde{\phi}_m$ from a
neighbourhood of $S_m$ to a neighbourhood
of $S'_m$ such that $D\tilde{\phi}_m$
induces the morphism $N\phi_m$.
On $\cup_{i=1}^{m-1} (S_i \cap S_m))$
we have that $\tilde{\phi}_m$
and $\Psi'$ agree.
Also $D\tilde{\phi}_m|_{\cup_{i=1}^{m-1} \cap s_m}$
is equal to
$D\Psi'|_{\cup_{i=1}^{m-1} \cap s_m}$.
Hence if we look at their graphs inside
$M \times M'$, they are tangent along
$(\cup_{i=1}^{m-1} (S_i \cap S_m) \times M' $
so we can perturb the graph of $\tilde{\phi}_m$
by a $C^1$ small amount so that it still stays the
graph of a diffeomorphism $\tilde{\phi}'$
and such that it agrees with $\Psi'$ on a neighbourhood
of $\cup_{i=1}^{m-1} (S_i \cap S_m)$.
Hence we can extend $\Psi'$ over $S_m$ with the appropriate properties.

Finally, we can use a standard Moser deformation argument
to deform $\Psi$ into a symplectomorphism
$\Phi$ such that the differential of $\Phi$
at $\cup_i S_i$ agrees with the differential of $\Psi$
at $\cup_i S_i$.
This ensures that $\Phi$ has the properties we want.

Now suppose in addition that a neighbourhood $NS_I$ of $S_I$
is identical to a neighbourhood of $NS'_I$ and that for $i \in I$, $\phi_i$ and $N\phi_i$
restricted to $S_I$ is the identity map for all $i$.
We can construct a smooth family $S_i^t$ $t \in [0,1]$
of positively intersecting submanifolds which are $C^1$ close to each other
so that on a small neighbourhood of $S_I$, $S_i^0$ is equal to $S_i$ for all $i$.
Also we want that globally, $S_i^1$ is equal to $S'_i$ for all $i$.
We can also assume that they are all symplectomorphic to each other,
so we have a smooth family of symplectomorphisms $\phi^t_i$ from $S_i$
to $S^t_i$  and also a smooth family of normal bundle maps $N\phi^t_i$
so that when $t = 1$ they all coincide with $\phi_i$ and $N\phi_i$.
By a parameterized version of the above discussion we have a smooth family of
symplectomorphisms $\Phi^t$ from a neighbourhood of $\cup_i S_i$ to a neighbourhood
of $\cup_i S^t_i$ such that $\Phi^t|_{S_i} = \phi^t_i$
and $D\Phi^t|_{NS_i} = N\phi^t_i$.
In addition we can assume that $\Phi^0$ is the identity map near $S_I$
(because we can choose our associated diffeomorphism $\Psi$ to have this property).
Because $\phi_i$ and $N\phi_i$ are the identity map on $S_I$ we have that $\Phi^1$
is $C^1$ small near $S_I$.
Let $V_t := \frac{d}{dt} \Phi^t$ be a smooth family of symplectic vector fields
defined near $S_I$. These are Hamiltonian vector fields near $S_I$ because $V_t = 0$
on $S_I$ and a small neighbourhood of $S_I$ deformation retracts onto $S_I$. 
Hence $\Phi^1$ is a Hamiltonian symplectomorphism near $S_I$ because it is Hamiltonian
isotopic to $\Phi^0$ near $S_I$ which is the identity map.
Hence $\Phi := \Phi^1$ has the required properties.
\qed

\begin{lemma} \label{lemma:extendingorthogonaltoneighbourhood}
Let $(M,\omega)$ be a symplectic manifold and let
$S$ be a compact symplectic submanifold
and $N,N'$ open subsets of $M$ such that the closure of $N$ is contained in $N'$.
Let $S_1,\cdots,S_k$ be positively intersecting symplectic submanifolds
with $\cap_i S_i = S$ and
such that they are orthogonal along $S$, and such that they
are also orthogonal in $N'$.
Then there exists a $C^1$ small perturbation of $S_1,\cdots,S_k$
to $S'_1,\cdots,S'_k$ such that
$S'_1,\cdots,S'_k$ are orthogonal on a small neighbourhood of $S$,
and $S'_1=S_1,\cdots,S'_k=S_k$ outside a small neighbourhood of $S$
and also inside $N$.
We may also arrange that $TS'_i = TS_i$ along $S$.
\end{lemma}
\proof of Lemma \ref{lemma:extendingorthogonaltoneighbourhood}.
We prove this by induction on the dimension of $M$.
So we suppose it is true in lower dimensions.
Choose a small neighbourhood $US$ of $S$
as in Lemma \ref{lemma:neighbourhoodofdivisor}
where $p : US \twoheadrightarrow S$ is a symplectic fibration
whose fibers are orthogonal to $S$.

Consider the positively intersecting submanifolds
$(S_i \cap S_k)_{i=1,\cdots,k-1}$ inside $S_k$.
By our induction hypothesis we can perturb them by
a $C^1$ small amount to $V_1,\cdots V_{k-1}$
such that:
\begin{enumerate}
\item $V_i = S_i$ inside $N'' \subset S_k$ and also outside a small
neighbourhood $M \subset S_k$ of $S \cap S_k$.
Here $N''$ is an open subset whose closure is contained in $N' \cap S$
and such that it contains closure of $N \cap S_k$.
\item $V_i$ intersect orthogonally along a small neighbourhood of $S \cap S_k$.
\end{enumerate}

Let $q$ be a projection from a small neighbourhood of $S_k$ to $S_k$
whose fibers are orthogonal to $S_k$.
Let $W_i := q^{-1}(V_i)$ for $i=1,\cdots,k-1$.
These symplectic submanifolds are tangent to $S_i$
along $S$ and orthogonal on a neighbourhood $p^{-1}(S)$ of $S$.
Because the normal bundle to $S_k$ is two dimensional
and $S_i,W_i$ are orthogonal to $S_k$ in the region $(S_k \cap N'')$ and of co-dimension $2$,
we have that their tangent spaces coincide in this region.
By Lemma \ref{lemma:neighbourhooduniqueness}
there is a $C^1$ small Hamiltonian symplectomorphism
(defined near $S_k$)
sending $W_i$ to $S_i$ for each $i$ in the region $p^{-1}(N'')$
and fixing $S_k$.
Choose a neighbourhood of $O$ of $S$ small enough so that
$N \cap O \subset p^{-1}(S_k \cap N'')$.
This means that $W_i = S_i$ inside $N \cap O$.
Hence we can perturb $S_i$ for $i = 1,\cdots,k-1$
by a $C^1$ small amount to $S'_i$ so that
it coincides with $W_i$ on a small neighbourhood of $S$ (inside $O$)
and is equal to $S_i$ outside $O$ and is unchanged in $N$.
This means that $S'_i = S_i$ inside $N$
and $S'_i$ is orthogonal on a neighbourhood of $S$.
Also $TS'_i = TS_i$ along $S$.
Hence $S'_i$ has the properties we want.
\qed

\begin{lemma} \label{lemma:orthogonaldivisorresultineuclideanspace}.
Let $S_1,\cdots,S_k$ be positively intersecting symplectic submanifolds inside some symplectic manifold $(M,\omega)$
and let $I \subset \{1,\cdots,k\}$.
Suppose that $S_1,\cdots,S_k$ intersect
orthogonally on some neighbourhood $N$ of $S_I \cap (\cup_{i \notin I} S_i)$.

Then there is a family of positively intersecting symplectic manifolds
$S_i^t$ ($t \in [0,1]$) with
\begin{enumerate}
\item $S_i^0 = S_i$.
\item $S_i^t = S_i$ for $i \notin I$.
\item $S_i^t = S_i$ on some open subset containing $S_I \cap (\cup_{i \notin I} S_i)$.
\item $S_i^1$ are all orthogonal on an arbitrarily small neighbourhood $NS_I$ of $S_I$.
\item $S_i^t = S_i$ outside an arbitrarily small neighbourhood of the closure of $NS_I$.
\item $\cap_{i \in I} S_i^t = S_I$.
\end{enumerate}
\end{lemma}
\proof of Lemma \ref{lemma:orthogonaldivisorresultineuclideanspace}.
By Lemma \ref{lemma:neighbourhoodofdivisor},
there is a neighbourhood $US$ of $S_I$
and a projection $p_I : US \twoheadrightarrow S_I$
such that the fibers of $p_I$ are symplectomorphic to $B_\epsilon$
and whose structure group is $U(n)$.
There is also a vector field $L$ tangent to the fibers
of $p_I$ such that $e^t (\phi_{-t})^* \omega = \omega + (e^t-1)p^*\omega|_S$
and such that it radial in the fibers.
We define $S := S_I$ and $\omega_S := \omega|_{S}$.

We prove the Lemma now by reducing it
to the linear case in Lemma \ref{lemma:orthogonallineardivisor}.
For each $i \in I$, consider the following manifold:
\[T_i := \cup_{q \in S} T_0 (S_i \cap p^{-1}(q)).\]
Here $T_0 (S_i \cap p^{-1}(q))$ means the
tangent space at $0$ of $S_i \cap p^{-1}(q)$
inside the tangent space of the linear fiber $p^{-1}(q)$
at $0$ which is canonically identified with $p^{-1}(q)$
(because $p^{-1}(q)$ is an open ball in $\C^n$).
We have that $T_i$ is symplectic near $S_I$.
By a parameterized version of Lemma \ref{lemma:neighbourhooduniqueness}
there is a $C^1$ small diffeomorphism preserving the fibers of $p_I$
sending $T_i$ to $S_i$ in the region $p^{-1}(N)$
for all $i \in I$
(possibly after shrinking $N$ slightly)
and such that it is a symplectomorphism when restricted to each fiber.
We push forward the $U(n-|I|)$ structure group of $p_I : US \twoheadrightarrow S_I$
via this fiberwise diffeomorphism as well so that $S_i$ restricted
to each fiber $p^{-1}(q)$ is linear for each $q \in N$.

For each $i \in I$, we can perturb $S_i$ by a $C^1$
small amount (without moving $S_I$) so that $S_i = T_i$
on a small neighbourhood $P$ of the zero section
and so that the $S_i$
are all still positively intersecting symplectic submanifolds.
We can also assume that this perturbation only happens
outside some neighbourhood of $S_I \cap (\cup_{i \notin I} S_i)$.

So from now on (after shrinking $US$)
we can assume that $S_i \cap p^{-1}(q)$
is linear inside $p^{-1}(q)$ for all $q$ and $i \in I$.
This means we have a codimension $0$ submanifold $\overline{S} \subset S_I$
with boundary disjoint from $S_I \cap (\cup_{i \notin I} S_i)$
such that $(S_i)_{i \in I}$ are orthogonal away from $\overline{S}$ and on
a neighbourhood ${\mathcal N}$ of the boundary of $\overline{S}$.
Let $U\overline{S}$ be equal to $p^{-1}(\overline{S})$.

By Lemma \ref{lemma:orthogonallineardivisor} (using $U\overline{S}$)
we can find a family of submanifolds $S_i^t$ so that
\begin{enumerate}
\item $S_i^t$ are positively intersecting symplectic manifolds
with respect to the symplectic form $\omega_{US} + T p^* \omega_S$ for sufficiently
large $T > 0$.
\item $\cap_i S_i^t$ is the zero section.
\item $S_i^0 = S_i$.
\item $S_i^t \cap p^{-1}({\mathcal N}) = S_i \cap p^{-1}({\mathcal N})$.
\item $S_i^1$ are all orthogonal along the zero section.
\item $S_i^t = S_i$ outside some small neighbourhood of $S$.
\end{enumerate}
We define $\phi_{-T}(S_i^t)$ to be equal
to $\phi_{-T}(S_i^t)$ inside $\phi^{-T}(US)$
and $S_i^t$ outside $\phi^{-T}(US)$.
These are also positively intersecting
submanifolds with the properties stated above
(because $\phi_{-T}^* \omega_{US} = \omega_{US} + (1 - e^{-T}) p^* \omega|_S$).


The problem is we want these symplectic manifolds
to be orthogonal on a {\it neighbourhood}
of $S_I$. But this can be done by perturbing $\phi_{-T}(S_i^t)$
by a $C^1$ small amount inside $U\overline{S}$
(by Lemma \ref{lemma:extendingorthogonaltoneighbourhood}),
hence $\phi_{-T}(S_i^t)$ has all the properties we want.
\qed

In the previous Lemma we have that
$\cap_{i \in I} S_i^t = S_I$.
We also have a smooth family of diffeomorphisms $\Psi_t$
from $S_i$ to $S_i^t$ which are the identity on $S_I$ and outside a small neighbourhood of $S_I$.
Because a neighbourhood of $S_I$ inside $S_i^t$ deformation
retracts to $S_I$, we have that
any $2$-cycle in $S_i^t$ near $S_I$
is homologous to a $2$-cycle in $S_I$.
This means that the integral of $\omega|_{S_i}$
over this $2$-cycle $C$ is the same as the integral
of $\omega|_{S_i^t}$ over $\Psi_*(C)$.
Also away from this neighbourhood we have that
$S_i^t=S_i$ so any $2$-cycle evaluated
on $\omega|_{S_i}$ is the same as the one
evaluated on $\Psi^*\omega_{S_i^t}$.
All of this implies that $\Psi_t^*([\omega]|_{S_i^t})$
is equal to $[\omega]|_{S_i}$
and by a Moser theorem we can then ensure that this is a symplectomorphism
if the $S_i$ are compact.
Hence (by using another Moser theorem) there is
a smooth family of symplectomorphisms
$P_i^t : \text{Nhd}(S_i) \rightarrow \text{Nhd}(S_i^t)$
where $\text{Nhd}$ means `a small neighbourhood of'.
These symplectomorphisms fix $S_I$ and hence
by a similar cohomological argument (as explained earlier),
we have that $P_i^t$ is in fact a Hamiltonian symplectomorphism,
so by using a cutoff function we can extend these symplectomorphisms
to Hamiltonian symplectomorphisms:
$P_i^t : M \rightarrow M$.

\proof of Lemma \ref{lemma:deformingdivisors}.
We basically proceed by induction on subsets $I$ where
$I \subset \{1,\cdots k\}$.
In order to do induction, we need a total order on this
finite set. Here is the following total order:
We say that $I \prec J$ when
\begin{enumerate}
\item $|I| > |J|$ or
\item $|I| = |J|$ with $I \neq J$ and the highest number in $J \setminus (J \cap I)$
is smaller than the highest number in $I \setminus (J \cap I)$.
\end{enumerate}

Fix some $I \subset \{1,\cdots,k\}$.
Suppose by our induction hypothesis, we have deformed
$S_1,\cdots,S_k$ through positively intersecting symplectic manifolds
so that they are orthogonal on a small neighbourhood $N$ of $\cup_{J \prec I} S_J$ and consider
$S_I$.
By Lemma \ref{lemma:orthogonaldivisorresultineuclideanspace},
we can deform $S_i$ through positively intersecting submanifolds
$S_i^t$ such that:
\begin{enumerate}
\item $S_i^0 = S_i$.
\item $S_i^t = S_i$ for $i \notin I$.
\item $S_i^1$ are all orthogonal on some small neighbourhood of $S_I$.
\item $S_i^t = S_i$ outside an arbitrarily small neighbourhood of $S_I$.
\end{enumerate}
Because $S_i^t = S_i$ outside an arbitrarily small neighbourhood $O$ of $S_I$,
we can assume that $(O \cap S_j) \setminus N$ is empty for all $j \notin I$.
This means that $S_i^1$ are still orthogonal along $N$
and also on some small neighbourhood of $S_I$.

Hence by induction we have proven that
we can deform $S_i$ through positively intersecting symplectic submanifolds
so that they are orthogonal on a neighbourhood of $\cup_i S_i$
and hence are orthogonal everywhere.

Let $(S'_i)^t$ be this deformation.
The problem is that $(S'_k)^t$ is not equal to $S_k$
for all $t$.
The paragraph before this proof tells us that
there is a sequence of symplectomorphisms
$P_k^t$ sending $S_k$ to $(S'_k)^t$
such that $P_k^0$ is the identity map.
So we can pull back $(S'_i)^t$ via $P_k^t$
for all $i$
and this ensures that we get a family
of positively intersecting submanifolds
$S_i^t$ such that $S_i^t$ are all positively intersecting,
$S_i^0 = S_i$ and $S_k^t = S_k$.
\qed

\subsection{Making the smooth affine variety nice at infinity}

\begin{lemma} \label{lemma:orthogonalneighbourhoods}
Suppose that $S_1,\cdots,S_k$ are positively intersecting codimension $2$
symplectic submanifolds of $(M,\omega)$ such that they are also orthogonal.
There exist small neighborhoods $US_i$ of $S_i$
and projections $\pi_i : US_i \twoheadrightarrow S_i$ such that
\begin{enumerate}
\item For $1 \leq i_1 < i_2 < \cdots < i_l \leq k$,
\[\pi_{i_l} \circ \cdots \circ \pi_{i_1} : \cap_{j = 1}^l US_{i_j}
\twoheadrightarrow S_{\{i_1,\cdots,i_l\}}\]
has fibers that are symplectomorphic to $\Pi_{j = 1}^l \D_{\epsilon}$
where $\D_\epsilon$ is the disk of radius $\epsilon$.
\item
If we look at a fiber $\Pi_{j = 1}^l \D_{\epsilon}$ of
$\pi_{i_l} \circ \cdots \circ \pi_{i_1}$, then for $1 \leq m \leq l$,
$\pi_{i_m}$ maps this fiber to itself.
It is equal to the natural projection
\[\Pi_{j = 1}^l \D_{\epsilon} \twoheadrightarrow  \Pi_{j = 1, j \neq m}^l \D_{\epsilon}\]
eliminating the $m$th disk $\D_\epsilon$.
\item
The symplectic structure on $US_i$ induces a natural connection for
$\pi_{i_l} \circ \cdots \circ \pi_{i_1}$ given by the $\omega$ orthogonal
vector bundles to the fibers.
We may require the associated parallel transport maps to be elements
of $U(1) \times \cdots \times U(1)$ where $U(1)$ acts on the disk $\D_\epsilon$
by rotation.
\end{enumerate}
\end{lemma}
\proof of Lemma \ref{lemma:orthogonalneighbourhoods}.
Suppose we have a bundle $p : V \rightarrow B$
with a $\prod_{i=1}^m U(1)$
structure where the fiber is $\D_\epsilon^m$ where
$\D_\epsilon$ is the $\epsilon$ ball in $\C$ and the $m$'th copy
of $U(1)$ rotates the $m$'th $\D_\epsilon$ factor in $\D_\epsilon^m$.
Suppose that the base $B$ has a symplectic structure $\omega_B$.
We can construct a symplectic structure on the total
space $V$ as follows:
Let $V_i \subset V$ be the subbundle whose
fiber is the subset is the $i$'th copy of $\D_\epsilon$ in $\D_\epsilon^m$.
This has a $U(1)$ structure group.
By Lemma \ref{lemma:neighbourhoodofdivisor},
there is a symplectic structure $\omega_{V_i}$ on $V_i$
such that the fibers of $V_i$ have the standard
symplectic structure on $\D_\epsilon \subset \C$.
We can ensure that the parallel transport maps are
in $U(1)$ as well as follows:
On $V_i$, there is an $S^1$ action
$A : S^1 \rightarrow \text{diffeo}(V_i)$
such that it fixes the map $p$ and rotates the fibers
of $V_i$
(i.e. it corresponds to the action given by rotating the fiber $\D_\epsilon
\subset \C$).
We define
\[\overline{\omega}_{V_i} := \int_{S^1} A(t)^* \omega dt.\]
This is a symplectic form if we shrink $\epsilon>0$ a bit.
Also the new symplectic form on $V_i$ ensures that the parallel
transport maps are in $U(1)$.
Let $P_i : V \twoheadrightarrow V_i$ be the natural projection
to $V_i$ where $(a_1,\cdots,a_m)$ in the fiber $\D_\epsilon^m$
is projected to the $i$'th $\D_\epsilon$ factor.
We define $\omega_V := \frac{1}{m} \sum_{i=1}^m P_i^*
\overline{\omega}_{V_i}$.
Any symplectic structure on a bundle with a
$\prod_{i=1}^m U(1)$ structure group as described above
is said to be {\it bundle compatible}.
The good thing about bundle compatible symplectic structures
is that the maps $P_i$ satisfy all the properties
as described in the start of this Lemma
(where $\pi_i$ is replaced by $P_i$).
Also if we have an open subset $U$
of $B$ and bundle compatible symplectic structures
on $p^{-1}(U)$
then by using similar methods and partitions
of unity, we can ensure that there is a bundle compatible
symplectic form on $V$ which coincides with the
symplectic form on $p^{-1}(U)$ possibly after shrinking
$p^{-1}(U)$ a tiny bit.
We will call the maps $P_i$ the {\it divisor projections}.

In order to prove our Lemma,
we basically proceed by induction on subsets $I$ where
$I \subset \{1,\cdots k\}$.
In order to do induction, we need a total order on this
finite set. Here is the following total order:
We say that $I \prec J$ when
\begin{enumerate}
\item $|I| > |J|$ or
\item $|I| = |J|$ with $I \neq J$ and the highest number in $J \setminus (J \cap I)$
is smaller than the highest number in $I \setminus (J \cap I)$ or
\end{enumerate}
Let $I_L$ be the $L$th subset in this total order.
If $I \subset \{1,\cdots,k\}$ then we write
$p_I$ to mean the map 
\[\pi_{i_l} \circ \cdots \circ \pi_{i_1} : \cap_{j = 1}^l US_{i_j}
\twoheadrightarrow S_{\{i_1,\cdots,i_l\}}\]
where $I = \{i_1,\cdots,i_l\}$
and $i_1,\cdots,i_l$ are distinct.
Suppose by our induction hypothesis
there exists a neighbourhood $U_{L-1}$ of $\cup_{i \leq L-1} S_{I_i}$
such that for each $i$ with $1 \leq i \leq k$, we have a neighbourhood $U_{L-1} \cap US_i$
of $S_i$ and a projection map $\pi_i : U_{L-1} \cap US_i \twoheadrightarrow S_i \cap U_{L-1}$
such that for each $I \subset \{1,\cdots,k\}$,
the map $p_I$ has the structure of
a $\prod_{i \in I} U(1)$ bundle which is bundle compatible
with the symplectic structure.
We will also assume that the maps $\pi_i$ are the associated
divisor projections for $p_I$ locally around $S_I$ inside
$U_{L-1}$.
Even though $US_i$ has not been constructed yet, we will use
the notation $U_{L-1} \cap US_i$ for the part of $US_i$ that has been constructed inside $U_{L-1}$.
We write $I := I_L$.
We want to extend $U_{L-1}$ to $U_L$ containing $S_I$
and $US_i \cap U_{L-1}$ to $US_i \cap U_L$ along with the maps $\pi_i$
so that they satisfy the properties as stated above.
The normal bundle of $S_I$ has a natural
$\prod_{i \in I} U(1)$ structure group because its tangent bundle
is the intersection $\cap_i TS_i$.
Hence by the previous discussion
and using an exponential map, we can extend the maps
$\pi_i$ over a neighbourhood and we can extend
the symplectic structure $\omega|_{U_{L-1}}$
to $\omega'$
over some neighbourhood of $S_I$ so that:
\begin{enumerate}
\item
$\pi_I$ has has the structure of a $\prod_{i \in I} U(1)$
bundle with fiber $\prod_{i \in I} \D_\epsilon$.
\item The symplectic structure $\omega'$
is bundle compatible.
\item the maps $\pi_i$ are the associated bundle projections.
\end{enumerate}
The symplectic structure $\omega'$
coincides with the symplectic structure
$\omega_M$ on $TS_I$, and the symplectic normal bundles are the
same although the symplectic form on one normal bundle
is a positive scalar multiple of the other.
We can make the symplectic forms coincide exactly along $S_I$
by pushing forward $\omega'$ by a diffeomorphism induced by a vector field which is
$0$ on $S_I$ and which is tangent to the normal bundle so that it rescales the symplectic
form $\omega'$ on this normal bundle so it coincides with $\omega_M$.
Hence by a Moser theorem we can ensure that
$\omega'$ coincides with $\omega_M$
(we have to deform the maps $\pi_i$ as well).
The problem is that because we have deformed the maps $\pi_i$,
the image of $\pi_i$ might not coincide with $S_i$
away from $\cup_{j \notin I} S_j$ anymore.
But by Lemma \ref{lemma:neighbourhooduniqueness},
there is a symplectomorphism which is the identity
near $\cup_{j \notin I} S_j$ moving all the $\pi_i$'s
so that their image is in $S_i$.
Hence $\pi_i$ has all the properties we want
and is defined on a neighbourhood $U_L$ of $S_I$.

Hence by induction we have maps $\pi_i$
defined on some neighbourhood of $\cup_i S_i$
with the properties stated as in the Lemma.
\qed

Let $S_1^t,\cdots,S_k^t$ $(t \in [0,1])$ be a smooth family of codimension $2$
symplectic submanifolds such that for each fixed $t$, $S_1^t,\cdots,S_k^t$
are positively intersecting.
We write $S^t := \cup_{i,j} (S_i^t \cap S_j^t)$.

\begin{lemma} \label{lemma:convexdeformationequivalent}
There is a smooth family of symplectomorphisms
$\Phi_t$ from $M \setminus S^0$
to $M \setminus S^t$ sending $S_i^0 \setminus S^0$
to $S_i^t \setminus S^t$ for each $i$.

In particular this means that there is a smooth family
of symplectomorphisms induced by $\Phi_t$
from $M \setminus (\cup_i S_i^0)$
to $M \setminus (\cup_i S_i^t)$.
\end{lemma}

We first need a preliminary Lemma:
\begin{lemma} \label{lemma:deformingasubmanifold}
Suppose $M$ is any symplectic manifold (open or closed
without boundary)
and let $\iota_t : S \hookrightarrow M$
($t \in [0,1]$) be any smooth family of proper symplectic embeddings of the symplectic
manifold $S$ (without boundary).
Then there is a family of Hamiltonian symplectomorphisms
$\phi_t : M \rightarrow M$ such that
$\phi_0$ is the identity and
$\phi_t$ sends the image of $\iota_0$ to the image
of $\iota_1$.
\end{lemma}
\proof of Lemma \ref{lemma:deformingasubmanifold}.
By a Moser theorem, there is a neighbourhood $U$
of $\iota_0(S)$ (which gets very thin near infinity)
and a smooth family of symplectic embeddings
$\iota'_t : U \hookrightarrow M$ such that
$\iota'_t|_S = \iota_t$.
We can also ensure that there is a smooth submersion
$\pi : U \twoheadrightarrow S$
whose fibers are all $\omega_M$-orthogonal to $S$.
We will also assume that $U$ deformation retracts to $S$.

For each $t$, there is a symplectic vector field
$V_t$ on $\iota'_t(U)$
defined as $\frac{d}{dt}(\iota'_t(p))$
at $p \in \iota'_t(U)$.
This vector field has the following property:
for any compact subset $\kappa$ of $S$,
there is an $\epsilon_\kappa>0$
such that for all $T \in [0,1]$
the flow $\phi^{V_t}_t(\iota'_T(\kappa))$ of $V_t$
is well defined and satisfies
$\phi^{V_t}_{t-T}(\iota'_T(\kappa)) \subset \iota'_t(\kappa)$
for all $|t - T| < \epsilon_\kappa$.
We will say that any vector field
satisfying this property {\it satisfies property Q}.

Let $\theta_{V_t}$ be the $\omega_M$-dual of $V_t$.
This is some closed $1$-form.
Because $U$ is homotopic to $S$, there exists
another closed $1$-form $\nu$ on $S$
such that $\nu + (\iota'_t)^* \theta_{V_t}$ is exact.
Let $\tilde{\nu} = \pi^* \nu$
and let $X_{\tilde{\nu}}^t$ be a smooth family of vector fields
defined only on $\iota'_t(U)$
whose $\omega_M$-dual is $\tilde{\nu}$.
We have that $X_{\tilde{\nu}}^t$ is tangent to $\iota'_t(S)$
because the fibers of $\pi$ are symplectically orthogonal
to $S$.
Because $X_{\tilde{\nu}}^t$ is tangent to $\iota'_t(S)$,
we have that $X_{\tilde{\nu}}^t + V_t$ satisfies property Q.
Because the $\omega_M$-dual of this vector field
is exact, we have a smooth family of Hamiltonians
$H_t : \iota'_t(U) \rightarrow \R$
whose associated Hamiltonian vector field is
$X_{\tilde{u}}^t + V_t$.
By using a cutoff function, we can assume that
$H_t$ smooth family of Hamiltonians on the whole of $M$.

The problem with this family of Hamiltonians is that
the associated Hamiltonian vector field may not be integrable.
In order to make it integrable (while still ensuring that
it satisfies property Q), we do the following:
Let $f$ be some exhausting function on $M$.
Then $f|_{\iota_t(S)}$ is also exhausting because
$\iota_t$ is a proper embedding.
We can perturb $f$ by a $C^0$ small amount near $\iota_t$
to create a new function $f_t$ which is equal to
$\pi^* f$ near $\iota_t(S) \subset \iota'_t(U)$
(here by abuse of notation we view $\pi$
as the pushforward of the map $\pi : U \twoheadrightarrow S$).
We can also assume that $f_t$
is a smooth family of exhausting functions.
The Hamiltonian flow of $f_t$ preserves $\iota_t(S)$
again because the fibers of $\pi$ are orthogonal to $S$.
Also because $f_t$ is a small perturbation of $f$,
we can assume that it is an exhausting function.

Choose some rapidly increasing positive function $g : \R \rightarrow \R$
such that $g \circ f_t + H_t$ is exhausting for all $t \in [0,1]$.
Because the Hamiltonian flow of $g \circ f_t$
preserves $\iota_t(S)$, we have that the flow of
$g \circ f_t + H_t$ satisfies property Q.
Also this Hamiltonian vector field is integrable for the
following reason:
Let $K_t$, $t \in [0,1]$ be a smooth family of exhausting Hamiltonians.
We will show that its Hamiltonian flow
is integrable.
Let $g : \R \rightarrow \R$ be a positive smooth function
such that
\[g(x) > \sup_{t \in [0,1], y \in K_t^{-1}(-\infty,x]} \frac{d
K_t}{dt}(y).\]
We choose $G$ to be any function with $G' = \frac{1}{g}$.
Let $p(t)$ be a path in $M$ satisfying $\frac{dp(t)}{dt} = X_{K_t}$.
We have
\[\frac{d}{dt}(K_t(p(t)) = \frac{dK_t}{dt}(p(t)) + dK_t(X_{K_t})
= \frac{dK_t}{dt}(p(t)) < g(K_t(p(t))).\]
Hence
\[\frac{d}{dt}(G(K_t(p(t))) < 1.\]
This implies that the function $t \rightarrow G(K_t(p(t)))$
is less than $t +C$ for some constant $C$.
This means that if $p(0)$ starts in $K_0^{-1}(a)$
then $p(t)$ must be contained in $K_t^{-1}(-\infty,G^{-1}( t + C)]$.
This implies that $X_{K_t}$ is integrable.

This means that $g \circ f_t + H_t$ has an integrable
Hamiltonian flow that satisfies property Q.
In particular the Hamiltonian flow
sends $\iota_0(S)$ to $\iota_t(S)$ at time $t$.
\qed

\proof of Lemma \ref{lemma:convexdeformationequivalent}.
For induction purposes, we will assume that this
Lemma is true in all dimensions less than $2n$.
For $j = 1,\cdots,n$, we will write $A_j^t$ as the union
$\cup_{|I| = n-j} (\cap_{i \in I} S_i^t)$.
This is the union of all dimension $2j$ strata in $\cup_i S_i^t$.
Suppose for $j' < j$, there is a smooth family
of symplectomorphisms $\Phi_{t,j'}$ from
$M \setminus A_{j'}^0$ to $M \setminus A_{j'}^t$
sending $A_{j'+1}^0 \setminus A_{j'}^0$
to $A_{j'+1}^t \setminus A_{j'}^t$.
We now wish to show the same thing for $j$.
We have (by our previous assumption) in particular
a smooth family of symplectomorphisms induced by
$\Phi_{t,{j-1}}$ sending the complement of $A_{j-1}^0$
to the complement of $A_{j-1}^t$ in $M$.
Note that $A_j^t \setminus A_{j-1}^t$
is a disjoint union of manifolds inside $M \setminus A_{j-1}^t$.
All these manifolds are compactified by some smooth normal
crossing divisor smoothly depending on $t$.
Because we have assumed that the Lemma is true
in all dimensions less than $2n$,
there is also a smooth family of symplectomorphisms $\Psi_t$
from $A_j^0 \setminus A_{j-1}^0$
to $A_j^t \setminus A_{j-1}^t$
starting at the identity symplectomorphism.
We have a smooth family of embeddings
$\iota_t := \Phi_{j-1,t}^{-1} \circ \Psi_t$ of
$A_j^0 \setminus A_{j-1}^0$ into $M \setminus A_{j-1}^0$.
By Lemma \ref{lemma:deformingasubmanifold},
there is a smooth family of Hamiltonian symplectomorphisms
$F_t : M \setminus A_{j-1}^0 \rightarrow M \setminus A_{j-1}^0$
where $F_0$ is the identity and $F_t$
sends the image of $\iota_0$ to the image of $\iota_t$.
Hence the smooth family of symplectomorphisms
$\Phi_{j-1,t} \circ F_t$ sends the image
of $\iota_0$ inside $M \setminus A_{j-1}^0$
to $A_j^t \setminus A_{j-1}^t$ inside $M \setminus A_{j-1}^t$.

Hence by induction, we have found a smooth family of symplectomorphisms
from $M \setminus A_{n-1}^0$ to $M \setminus A_{n-1}^t$ starting at the identity
sending $A_n^0 \setminus A_{n-1}^0$ to $A_n^t \setminus A_{n-1}^t$.
Because $A_{n-1}^t = S^t$ and
$A_n^t \setminus A_{n-1}^t$ is the disjoint union
of $S_i^t \setminus S^t$ for all $i$,
the previous statement is the statement of the Lemma.
Hence we have proven the Lemma.
\qed

Let $S_1,\cdots,S_k$ be positively intersecting inside
$M$. We can assume that $S_i \setminus \cup_{i \neq j} S_j$
is connected because we can replace a disconnected such
manifold with a union of connected manifolds.
Let $\theta$ be a $1$ form on the complement
$M \setminus (\cup_i S_i)$ such that $d\theta = \omega_M$.
Let $W$ be a small symplectic submanifold of $M$
of dimension $2$ symplectomorphic to some small disk $\D_\delta \subset \C$
of radius $\delta>0$.
Suppose that $W$ is disjoint from $S_j$, $j \neq i$
and $W$ intersects $S_i$ orthogonally at $0 \in \D_\delta$
(i.e. it is tangent to the normal bundle of $S_i$).
Let $(r,\vartheta)$ be polar coordinates on $\D_\delta$.
Then $\theta$ pulls back to a $1$-form on $W \setminus \{0\}$.
We have that $d\theta = rdr \wedge d\vartheta$
on the punctured disk $\D_\delta \setminus \{0\}$.
Hence $\theta$ is cohomologous to $\frac{r^2}{2} d\vartheta + \kappa_i d\vartheta$
for some constant $\kappa_i$.
Suppose that I had some other disk $W'$ intersecting $S_i$ orthogonally
and disjoint from the other $S_j$'s.
Because $S_i \setminus \cup_{j \neq i} S_j$ is connected,
there is a smooth family $W_t$ joining $W'$ and $W$.
Hence the constant $\kappa_i$ associated to $W'$
is the same as the constant $\kappa_i$ associated to $W$.
Hence $\kappa_i$ is an invariant of $\theta$ and $S_i$.
We call $\kappa_i$ the {\it wrapping of $\theta$ around $S_i$}.

\begin{lemma} \label{lemma:wrappingorthogonal}
Suppose that $S_1,\cdots,S_k$ are orthogonal positively intersecting
submanifolds of a compact symplectic manifold $(M,\omega_M)$
and suppose that there is a $1$-form $\theta$
such that
\begin{enumerate}
\item $d\theta = \omega_M$.
\item The wrapping of $\theta$ around $S_i$ for each $i$ is negative.
\end{enumerate}
Then there exists a function $f$ defined on $M \setminus \cup_i S_i$
such that $(M \setminus \cup_i S_i,\theta + df)$
has the structure of a finite type  convex symplectic manifold.

The form $\theta + df$ restricted to a fiber $\D_\epsilon$
of the maps $\pi_i$ minus the origin as described in
Lemma \ref{lemma:orthogonalneighbourhoods}
is $(\frac{r_i^2}{2} + \kappa_i) d\vartheta_i$ where $(r_i,\vartheta_i)$ are polar
coordinates for $\D_\epsilon$ and $\kappa_i<0$ is the wrapping number
of $\theta$ around $S_i$.
\end{lemma}
\proof of Lemma \ref{lemma:wrappingorthogonal}.
We proceed by induction on the ordering $\prec$
mentioned in Lemma \ref{lemma:orthogonalneighbourhoods}.
Let $\pi_I : US_I \rightarrow S_I$
be the maps as described
in the proof of Lemma \ref{lemma:orthogonalneighbourhoods}.
Suppose inductively we have for all $I' \prec I$ we have a function
$f'$ such that $\theta + df'$ is equal to
$(\sum_{i \in I'} (\frac{r_i^2}{2} + \kappa_i) d\vartheta_i)$ on the fiber
$\prod_{i \in I'} \D_\epsilon$.
We now wish to modify $f'$ on $US_I$ away from a very small neighbourhood
of $\cup_{j \notin I} S_j$ so that it has the properties we want.

We do this as follows:
If we look near $\cup_{j \notin I} S_j$, $\theta + df'$ has the form
we want because $I \cup \{j\} \prec I$ and
the map $\pi_I$ inside $US_{I \cup \bracket{j}}$ preserves the fibers
of $\pi_{I \cup \bracket{j}}$ and corresponds to the projection
\[\prod_{i \in I \cup \bracket{j}} \D_\epsilon \twoheadrightarrow \prod_{i \in I} \D_\epsilon.\]
eliminating the $\D_\epsilon$ corresponding to $j$.
Because $\theta + df'$ restricted to a fiber of $\pi_I$
is cohomologous to $\sum_{i \in I} (\frac{r_i^2}{2} + \kappa_i) d\vartheta_i$,
there is a smooth function $f''$ from $US_I$ to $\R$ such that
$\theta + df' + df''$ restricted to each fiber is of the form
$\sum_{i \in I} (\frac{r_i^2}{2} + \kappa_i) d\vartheta_i$.
Such a function exists because $\theta + df' - \sum_{i \in I}
(\frac{r_i^2}{2} + \kappa_i) d\vartheta_i$
restricted to each fiber of $\pi_I$ is exact, so for
each fiber $\pi_I^{-1}(q)$ we have a function $f_q$
which is unique up to adding a constant
such that this $1$-form is equal to $df_q$.
We can assume that $f_q$ smoothly depends on $q$
by adding some function of the form $\kappa \circ \pi_I$
(this basically follows from a parameterized version of the
Poincar\'{e} Lemma).
Hence we can view $f_q$ as a function $f''$ which has the properties we want.

Near $\cup_{j \notin I} S_j$ we have that $\theta + df'$
already has the form we want, hence
$f''$ must be of the form $\pi_I^* h$ in this region.
By extending $h$ to the whole of $S_I$, we can subtract
$\pi_I^* h$ from $f''$ so that $f'' = 0$ near $\cup_{j \notin I} S_j$.
This means that $\theta + df' + df''$ restricted
to each fiber of $\pi_I$ is equal to $\sum_{i \in I} (\frac{r_i^2}{2} + \kappa_i) d\vartheta_i$,
and also because $f'' = 0$ near $\cup_{j \notin I} S_j$, it also
has this form for all other maps $\pi_{I'}$ where $I' \prec I$.
Hence by induction we have shown that there is some function $f$
such that $\theta + df = \sum_{i \in I} (\frac{r_i^2}{2} + \kappa_i) d\vartheta_i$
inside the fibers of $\pi_I$ for all $I$.
This is the second part of the Lemma.

We now wish to show that we have the structure of a finite type convex symplectic manifold.
Let $\nu : (0,\epsilon) \rightarrow \R$ be a function which is equal
to $0$ near $\epsilon$ and tends to infinity as we reach $0$ and which has non-positive derivative.
We assume it has negative derivative near $0$.
Then we view $\nu(r_i)$ as a function on the complement of $S_i$
as we can extend it by zero away from $US_i$.
We have that $\sum_i \nu(r_i)$ is an exhausting function.
Let $X_{\nu(r_i)}$ be the Hamiltonian vector field of $\nu(r_i)$.
Because the $\omega$-orthogonal plane bundle to the fibers $\D_\epsilon$
of $US_i$ are contained in the level sets of $r_i$,
we have that $X_{\nu(r_i)}$ is tangent to the fibers
of $\pi_i$ inside $US_i$ and is zero elsewhere.
Also $X_{\nu(r_i)}$ restricted to some fiber $\D_\epsilon$
is equal to $X_{\nu(r_i)|_{\D_\epsilon}}$
which is equal to $-\frac{\nu'(r_i)}{r_i} \frac{\partial}{\partial \vartheta_i}$.
Hence $X_{\nu(r_i)} (\theta + df) < 0$ for $r$ small enough because $\kappa_i$ is negative.
If $X_{\theta + df}$ is the $\omega_M$-dual of $\theta+df$
then $X_{\theta+df}( d(\sum_i \nu(r_i)) ) = -(\sum_i X_{\nu(r_i)})(\theta + df) \geq 0$.
This is greater than zero near $\cup_i S_i$ because for each point $p$ near $\cup_i S_i$,
we have that $\nu'(r_i) > 0$ for some $i$.
Hence $(M \setminus \cup_i S_i, \theta +df)$ has the structure of a finite
type convex symplectic manifold.
This proves the Lemma.
\qed

Let $P$ be a smooth projective variety, and $D$ some effective divisor so that
\begin{enumerate}
\item $D$ is a smooth normal crossing divisor.
i.e. it is a union of smooth complex hypersurfaces which are transversally intersecting.
\item $P \setminus D$ is isomorphic to $A$ where $A$ is a smooth affine variety.
\end{enumerate}
Let $L$ be a line bundle on $P$ and $\|\cdot\|$ a metric on $L$ so that if $F$
is its curvature form then $-iF$ restricted to $A$ is a K\"{a}hler form
but this may not be true along $D$.
Let $s$ be a section of $L$ so that $s^{-1}(0) = D$.
We will call $d^c \log{\|s\|}$ a {\it compactification $1$-form}
associated to $A$.
Recall that we have another $1$-form on $A$ coming from an embedding $\iota : A \hookrightarrow \C^N$.
This is described before Lemma \ref{lemma:affineLiouvilledomain}.
This is given by $\iota^* \sum_i\frac{r_i^2}{2} d\vartheta_i$ where
$(r_i,\theta_i)$ are polar coordinates on the $i$th factor of $\C^N$.
The following lemma is almost exactly the same as
\cite[Lemma 4.4]{Seidel:biasedview}.
\begin{lemma} \label{lemma:formsonasmoothaffinevariety}
We have that $(A,d^c\log{\|s\|})$ makes $A$ into a finite type convex symplectic manifold
convex deformation equivalent to $(A,\theta_A)$.
\end{lemma}
\proof of Lemma \ref{lemma:formsonasmoothaffinevariety}.
This lemma will be done in three steps. In step 1 we show that
$(A,d^c\log{\|s\|})$ is a finite type convex symplectic manifold.
In step 2 we will show that if $d^c\log{\|s'\|'}$ is another compactification $1$-form
then $(A,d^c\log{\|s\|})$ is convex deformation equivalent to
$(A,d^c\log{\|s'\|'})$. Finally in step 3 we will show that
$(A,d^c\log{\|s\|})$  is convex deformation equivalent to $(A,\theta_A)$.

{\it Step 1:}
Let $p \in D$ and choose local holomorphic coordinates
$z_1,\cdots,z_n$ and a trivialization of $L$ around $p$ so that
$s = z_1^{w_1} \cdots z_n^{w_n}$ ($w_i \geq 0$).
The metric $\|.\|$ on $L$ is equal to $e^{\psi}|.|$
for some function $\psi$ with respect to this trivialization
where $|.|$ is the standard metric on $\C$.
So \[d\log{\|s\|} = -\psi - (\sum_i w_i d\log{|z_i|}).\]
If we take the vector field
$Y := -z_1 \partial_{z_1} \cdots - z_n \partial_{z_n}$,
then $d\log{(|z_j|)}(Y) = -1$ and $Y.\psi$ tends to zero.
Hence $d\log{\|s\|}$ is non-zero near infinity.
If $X_{d^c \log{\|s\|}}$ is the $dd^c \log{\|s\|}$-dual of $d^c \log{\|s\|}$
then $-d\log{(\|s\|)}(X_{d^c \log{\|s\|}}) = \|d\log{\|s\|}\|^2 > 0$ near infinity.
Hence $(A,d^c \log{\|s\|})$ is a finite type convex symplectic manifold.

{\it Step 2:}
Suppose now that $P',D',L',s',\|\cdot\|'$ are different choices of compactification,
effective divisor, line bundle, section and metric satisfying similar properties
to $P,D,L,s$ and $\|\cdot\|$. Here $P' \setminus D'$ is isomorphic to $A$.
By applying the Hironaka resolution of singularities theorem again
\cite{hironaka:resolution} we have a third compactification $P''$
and morphisms $\pi : P'' \rightarrow P$, $\pi' : P'' \rightarrow P$.
We pull back $L$ and $L'$ to $P''$  as well as the sections $s$, $s'$ giving us new line bundles
$\widetilde{L},\widetilde{L}'$ and sections $\widetilde{s}$ and $\widetilde{s}'$.
By abuse of notation we write $\|\cdot\|$ and $\|\cdot\|'$ for the metrics on $\widetilde{L}$
and $\widetilde{L}'$ which are the pullbacks of $\|\cdot\|$ and $\|\cdot\|'$ respectively.
Again we look at some local coordinate chart $z_1,\cdots,z_n$ and vector field $Y$ around some
point $p \in D''$ where $D'' = P'' \setminus A$.
Using the same arguments as before we have
\[ \left(t \log{\|\widetilde{s}\|} + (1-t) \log{\|\widetilde{s}'\|'}\right)(Y) > 0\]
near infinity for all $t \in [0,1]$.
Hence $(A,d^c \log{\|s\|})$ and $(A,d^c \log{\|s'\|'})$ are convex deformation equivalent.

{\it Step 3:}
In this step we will need to refer to the proof of Lemma
\ref{lemma:affineLiouvilledomain} so the reader must be familiar with this lemma first.
Let $R := \sum_i \frac{r_i^2}{4}$.
We have that $\theta_A = -d^c R$. Basically from the proof of Lemma
\ref{lemma:affineLiouvilledomain} we have a compactification $1$-form
equal to $-d^c f(R)$ where $f$ is a non-decreasing function.
We have that
$(A,-d^c(tR + (1-t) f(R)))$ is a convex deformation equivalence from
from $(A,\theta_A)$ to $(A,-d^c f(R))$.
The reason why this is a convex deformation equivalence is that the level
sets of $R$ and $f(R)$ coincide (up to reparameterization).
By Step 2 we have $(A,\theta)$ is convex deformation equivalent to $(A,-d^c f(R))$.
Hence $(A,\theta)$ is convex deformation equivalent to $(A,\theta_A)$.
\qed

\begin{lemma} \label{lemma:smoothaffinevarietywrappingnumbers}
Let $M$ be a smooth projective variety and let $\cup_i S_i$
be a smooth normal crossing divisor such that
$\sum_i a_i S_i$ is ample for some $a_1,\cdots,a_k \in \N \setminus \{0\}$.
The complement $A := M \setminus \cup_i S_i$ is an affine variety
and it has a natural $1$-form $\theta$ with $d\theta = \omega_M$.
We have that $(A,\theta)$ is convex deformation equivalent to
$A$ with the natural $1$-form $\theta_A$ coming from some embedding of $A$ into $\C^N$
(see Section \ref{section:liouvilldomaindefinition}).
Then $S_1,\cdots,S_k$ are positively intersecting
and the wrapping number of $\theta$ around $S_i$ is negative.
\end{lemma}
\proof of Lemma \ref{lemma:smoothaffinevarietywrappingnumbers}.
The natural $1$-form $\theta$ is a compactification $1$-form as described earlier.
We have an ample line bundle $L$ coming from $\sum_i a_i S_i$.
There is a metric $\|.\|$ on $L$ whose associated curvature
form $F$ has the property that $iF$ is a positive $(1,1)$-form.
This is our symplectic form $\omega_M$ on $M$.
Choose a holomorphic section $s$ of $L$ whose zero set is exactly
$\sum_i a_i S_i$ (i.e. the zero set has multiplicity $a_i$ at $S_i$).
Then our $1$-form $\theta$ is $d^c\log{\|s\|}$ where $d^c = d \circ J$
where $J$ is our complex structure on $M$.

First of all, the divisors $S_i$ are positively intersecting submanifolds
because $J$ is compatible with the symplectic form $\omega_M$
and all the $S_i$'s are holomorphic submanifolds.
So the hard part is proving that $\theta$ wraps negatively around $S_i$.
Let $p$ be a point in $S_i$ disjoint from $\cup_{j \neq i} S_j$ and
choose a small holomorphic disk $\D_\delta$ disjoint from $\cup_{j \neq i} S_j$
and intersecting $S_i$ at $p$ only.
We can assume that this disk is intersecting $S_i$ orthogonally
because the $\omega_M$-orthogonal bundle to $S_i$ is holomorphic.
The line bundle $L$ restricted to this disk $\D_\delta \subset \C$
is equal to $a_i \{0\}$ where $\{0\}$ means the divisor at $0$
and the section $s_i$ has a zero of multiplicity $i$ at $0$.
This means we can choose a trivialization of $L$
so that $s_i$ is equal to $z^{a_i}$ where $z$ is the complex coordinate on $\D_\delta$.

We have that $dd^c \log{\|z^{a_i}\|} = dd^c \log{\|1\|}$
so the wrapping number of $\theta$ around $S_i$
corresponds to the cohomology class of
\[d^c \log{\|z^{a_i}\|} - d^c \log{\|1\|} = a_i d^c \log{|z|} + d^c \log{\|1\|} - d^c \log{\|1\|} \]
\[= a_i d^c \log{|z|} = - a_i d\vartheta.\]
Here $|.|$ is the standard euclidean metric on $\D_\delta \subset \C$
and $\vartheta$ is the angular coordinate on $\D_\delta$.
This is a negative multiple of $d\vartheta_i$ and
hence we have that $\theta$ wraps around $S_i$ negatively.
By Lemma \ref{lemma:formsonasmoothaffinevariety} we get that $(A,\theta_A)$
is convex deformation equivalent to $(A,\theta)$.
\qed

Putting all these Lemmas together gives us the following theorem:
\begin{theorem} \label{theorem:affinedeforming}
Let $A$ be a smooth affine variety.
Then $A$ is convex deformation equivalent to a
finite type convex symplectic manifold $(W,\theta_W)$
with the following properties:
\begin{enumerate}
\item $W$ is symplectomorphic to $M \setminus \cup_i S_i$
where $S_i$ are co-dimension $2$ symplectic submanifolds transversely intersecting.
\item There are neighborhoods $US_i$ of $S_i$ and fibrations
$\pi_i : US_i \twoheadrightarrow S_i$ satisfying the properties
stated in Lemma \ref{lemma:orthogonalneighbourhoods}.
\item $\theta_W$ restricted to the fiber $\D_\epsilon$
of $\pi_i$ is equal to $(r_i^2 + \kappa_i) d\vartheta_i$
for some $\kappa_i < 0$.
\end{enumerate}
\end{theorem}
\proof of Theorem \ref{theorem:affinedeforming}.
We first compactify $A$ to a smooth projective variety
and then by using \cite{hironaka:resolution}, we blow up this projective variety so that
$A = M \setminus \cup_i S'_i$
for some transversely intersecting complex hypersurfaces $S'_i$ inside a new smooth projective variety $M$.
These are positively intersecting symplectic submanifolds of $M$
and the natural $1$-form $\theta$ ($d\theta = \omega_M$) on $A$ wraps negatively
around $S_i$ for each $i$ by Lemma \ref{lemma:smoothaffinevarietywrappingnumbers}.

By Lemma \ref{lemma:deformingdivisors},
there is a smooth family of positively intersecting symplectic submanifolds
$S_i^t$ such that $S_i^0 = S'_i$
and $S_i^1$ are orthogonal.
Let $S^t := \cup_{i,j} (S_i^t \cap S_j^t)$.
By Lemma \ref{lemma:convexdeformationequivalent},
there is a smooth family of symplectomorphisms $\Phi_t$
from $M \setminus S^0$ to $M \setminus S^t$
starting at the identity and sending
$S_i^0 \setminus S^0$ to $S_i^t \setminus S^t$.

This means that the symplectomorphisms $\Phi_t$ also induce symplectomorphisms
from $M \setminus \cup_i S_i^0$
to $M \setminus \cup_i S_i^t$.
Also $\theta_t := (\Phi_t|_{M \setminus \cup_i S_i^0})_* \theta$
wraps around $S_i^t$ negatively
because the symplectomorphism $\Phi_t|_{M \setminus \cup_i S_i^0}$
extends smoothly over $S_i$ to $\Phi_t$.
This means by Lemma \ref{lemma:wrappingorthogonal},
there is an exact $1$-form $df$ such that
$\theta_1 + df$ makes $M \setminus \cup_i S_i^1$
into a finite type convex symplectic manifold with the property
that $\theta_1 + df$ restricted to each fiber $\D_\epsilon$ of $\pi_i$
is equal to $(r_i^2 + \kappa_i) d\vartheta_i$ where $\kappa_i <0$.
So if we set $W := M \setminus \cup_i S_i^1$
and $\theta_W := \theta_1 +df$ then we have proven the Lemma.
\qed

\section{The growth rate of affine varieties}
\label{subsection:smoothaffinevarieities}

\begin{defn} \label{defn:exactsymplectomorphicatinfinity}
Let $(M_1,\theta_1),(M_2,\theta_2)$ be two exact symplectic
manifolds.
We say that $M_1,M_2$ are {\it exact symplectomorphic at infinity}
if there are open subsets $U_i \subset M_i$
and a symplectomorphism
\[\Psi : M_1 \setminus U_1 \rightarrow M_2 \setminus U_2\]
such that
\begin{enumerate}
\item
The closure of $U_i$ is compact inside $M_i$.
\item
$\Psi$ sends compact sets in $M_1 \setminus U_1$
to compact sets in $M_2 \setminus U_2$.
\item
$\Psi$ is an exact symplectomorphism
(i.e. $\Psi^*\theta_2 = \theta_1 + df$).
\end{enumerate}
\end{defn}

\begin{defn} \label{defn:convexdeformationequivalentatinfinity}
Let $(W_1,\theta_{W_1}),(W_2,\theta_{W_2})$ be two convex symplectic manifolds.
We say that they are {\it convex deformation equivalent at infinity}
if there is a sequence of convex symplectic manifolds $(Q_i,\theta_{Q_i})$
for $i = 1,\cdots k$ such that:
\begin{enumerate}
\item  $Q_1$ is convex deformation equivalent to $W_1$.
\item $Q_k$ is convex deformation equivalent to $W_2$.
\item $Q_i$ is convex deformation equivalent to $Q_{i+1}$
or they are exact symplectomorphic at infinity.
\end{enumerate}
\end{defn}
Note that being {\it convex deformation equivalent at infinity}
is an equivalence relation.

\begin{theorem} \label{thm:polynomialgrowthrate}
Suppose that a finite type convex symplectic manifold $B$
and a smooth affine variety $A$ are convex deformation equivalent at infinity.
Then filtered directed system $(SH_*^{\#}(B,\theta_B,\lambda))$
is isomorphic to a filtered directed system
$(V_x)$ where the dimension $|V_x|$ satisfies
$|V_x| \leq P(x)$ for some polynomial $P$.
The degree of this polynomial is less than or equal to $m_A$
where $m_A$ is defined in the introduction.
\end{theorem}

This has the following direct corollary:
\begin{theorem} \label{theorem:polynomialgrowthrate}
We have $\Gamma(B) \leq m_A$.
\end{theorem}

Theorem \ref{thm:maingrowthratebound} follows directly
from Theorem \ref{theorem:polynomialgrowthrate}
where we set $B = A$.
Theorem \ref{thm:mainfillingobstruction}
also follows from Corollary \ref{thm:polynomialgrowthrate}.
Here is a statement of this theorem:
{\it Suppose that the boundary of a Liouville domain $M$
is algebraically Stein fillable, then $\Gamma(M) \leq m_A$.}
\proof of Theorem \ref{thm:mainfillingobstruction}.
The boundary of $M$ is contactomorphic to the boundary of $\overline{A}$.
This means that we can deform $M$ through a family of Liouville
domains to a new Liouville domain $M'$ such that
the contact form on $\partial M'$ coincides exactly with
the contact form on $\partial \overline{A}$.
This means that the completions $\widehat{M'}$
and $A := \widehat{\overline{A}}$ are exact symplectomorphic at infinity
because their cylindrical ends are identical.
Also $\widehat{M'}$ is convex deformation equivalent to $\widehat{M}$
which means that $\widehat{M}$ and $\widehat{A}$
are convex deformation equivalent at infinity.
By Theorem \ref{theorem:polynomialgrowthrate},
we have that $\Gamma(\widehat{M}) < \infty$
which proves Theorem \ref{thm:mainfillingobstruction}.
\qed

Before we prove Theorem \ref{thm:polynomialgrowthrate},
We need some lemmas:
The first Lemma gives us a slightly more general definition of growth
rate which will be useful for our purposes.
Let $(H,J)$ be a growth rate admissible Hamiltonian
on some finite type convex symplectic manifold $(W,\theta_W)$ and let
$\kappa>0$ be any constant.
Let $\lambda_1 \leq \lambda_2$. Then there is a natural morphism
from $SH_*^{[0,\kappa \lambda_1]}(\lambda_1 H,J)$
to $SH_*^{[0,\kappa \lambda_2]}(\lambda_2 H,J)$
given by first composing the morphism
\[SH_*^{[0,\kappa \lambda_1]}(\lambda_1 H,J)
\rightarrow SH_*^{[0,\kappa \lambda_1]}(\lambda_2 H,J)\]
induced by the continuation map followed by the morphism
\[SH_*^{[0,\kappa \lambda_1]}(\lambda_2 H,J)
\rightarrow SH_*^{[0,\kappa \lambda_2]}(\lambda_2 H,J)\]
given by the natural inclusion map of chain complexes.
We could have done this the other way around
by starting with an inclusion of chain complexes:
\[SH_*^{[0,\kappa \lambda_1]}(\lambda_1 H,J)
\rightarrow SH_*^{[0,\kappa \lambda_2]}(\lambda_1 H,J)\]
and then composing it with a continuation map.
But this gives us exactly the same map.
Because we can swap inclusion maps and continuation maps,
we have that these maps satisfy functoriality properties
which means we get a filtered directed system:
\[(SH_*^{[0,\kappa \lambda]}(\lambda H,J))\]

\begin{lemma} \label{lemma:growthrateactionbound}
The filtered directed system described above is isomorphic
to $(SH_*^{\#}(W,\theta_W))$.
\end{lemma}
\proof of Lemma \ref{lemma:growthrateactionbound}.
We have a natural map $\phi$ from
$SH_*^{[0,\kappa \lambda]}(\lambda H,J)$
to
$SH_*^{\#}(\lambda H,J)$ given by inclusion of chain complexes.
This commutes with the filtered directed system maps
so is a morphism of filtered directed systems.
There is a constant $K>1$ and a $1$-form $\theta$
such that $\theta-\theta_W$ is exact, $\kappa K > 1$ and $-\theta(X_H) - H < \kappa K$
by the action bound property.
We have a map $\phi'$ from
$SH_*^{\#}(\lambda H,J)$
to
$SH_*^{[0,\kappa K \lambda]}(K\lambda H,J)$
constructed as follows:
First of all we have an isomorphism from
$SH_*^{\#}(\lambda H,J)$
to
$SH_*^{[0,\kappa K \lambda)}(\lambda H,J)$
because all the $1$-periodic orbits have action less than $\kappa K \lambda$.
We then compose this isomorphism with
the continuation map:
\[SH_*^{[0,\kappa K \lambda]}(\lambda H,J) \rightarrow
SH_*^{[0,\kappa K \lambda]}(K\lambda H,J).\]
Again this commutes with the filtered directed system maps
because continuation maps and action inclusion maps commute.

The map $\phi \circ \phi'$ is the composition:
\[SH_*^{\#}(\lambda H,J) \rightarrow
SH_*^{[0,\kappa K \lambda]}(\lambda H,J) \rightarrow\]
\[SH_*^{[0,\kappa K \lambda]}(K\lambda H,J) \rightarrow
SH_*^{\#}(K\lambda H,J).\]
This is the natural continuation map because action inclusion
and continuation maps commute so the above composition
is equal to the following composition:
\[SH_*^{\#}(\lambda H,J) \rightarrow
SH_*^{[0,\kappa K \lambda]}(\lambda H,J) \rightarrow\]
\[SH_*^{\#}(\lambda H,J) \rightarrow
SH_*^{\#}(K\lambda H,J).\]
Also because all the $1$-periodic orbits of $\lambda H$ have action
less than $\kappa K \lambda$, we have that the composition
of the first two maps is the identity map.
Hence $\phi \circ \phi'$ is a continuation map.

We also have that $\phi' \circ \phi$ is a continuation map
because it is equal to the natural composition:
\[SH_*^{[0,\kappa \lambda]}(\lambda H,J) \rightarrow
SH_*^{\#}(\lambda H,J) \rightarrow\]
\[SH_*^{[0,K \kappa \lambda]}(\lambda H,J) \rightarrow
SH_*^{[0,K \kappa \lambda]}(K\lambda H,J).\]
Because all $1$-periodic orbits of $\lambda H$ have action less
than $\kappa K \lambda$ the composition of the first two maps
is identical to the natural inclusion map.
Hence $\phi' \circ \phi$ is a directed system map
and so $\phi,\phi'$ give us our isomorphism of filtered directed systems.
\qed

\begin{lemma} \label{lemma:hamiltoniangrowthratebound}
Suppose we have a Hamiltonian $H : W \rightarrow \R$, a function $P : \R \rightarrow \R$
and a small open neighbourhood ${\mathcal N}$ of $H^{-1}(0)$ such that:
\begin{enumerate}
\item $H$ satisfies the Liouville vector field property.
\item For every $\lambda$ outside some discrete subset,
there is a $C^2$ small perturbation $H_\lambda$
of $\lambda H$ such that all the $1$-periodic orbits
of $H_\lambda$ inside ${\mathcal N}$ are non-degenerate and the number of
such orbits is bounded above by $P(\lambda)$.
\end{enumerate}
Then there is a filtered directed system $(V_\lambda)$ isomorphic
to $(SH_*^{\#}(W,\theta_W))$
such that the rank of $V_\lambda$ is bounded above by $P(\lambda)$.
\end{lemma}
\proof of Lemma \ref{lemma:hamiltoniangrowthratebound}.
If $H$ was growth rate admissible and ${\mathcal N} = W$ then this lemma would
be fairly straightforward as the rank of $SH_*$ is bounded
above by the number of non-degenerate orbits.
The problem is that it may not be growth rate admissible.
Instead we will construct a growth rate admissible Hamiltonian $H^p$
such that all the $1$-periodic orbits of sufficiently low action
are the same as the ones of $H$ and then invoke Lemma \ref{lemma:growthrateactionbound}.

By Lemma \ref{lemma:positveactionhamiltonian},
there is a growth rate admissible pair $(H^p,J^p)$ such that
\begin{enumerate}
\item $H^p = H$ on a small neighbourhood of $H^{-1}(0)$
and $(H^p)^{-1}(0) = H^{-1}(0)$.
\item $-\theta_H(X_{H^p}) - H^p \geq 0$ everywhere.
\item $-\theta_H(X_{H^p}) - H^p > 0$ when $H^p > 0$.
\item $-\theta_H(X_{H^p}) - H^p$ is greater than some constant $\delta_H>0$
outside some compact set.
\end{enumerate}
Here $\theta_H$ is the $1$-form that makes
$H$ satisfy the Liouville vector field property.

Let $A$ be the function $-\theta_H(X_{H^p}) - H^p$.
The level sets of $H^p$ near $(H^p)^{-1}(0)$
are compact because $H^p$ satisfies the bounded below property.
Let $\nu_H>0$ be a constant so that
for all $x \in (0,\nu_H]$ we have that $(H^p)^{-1}(x)$
is compact and regular and contained inside ${\mathcal N} \cap \{H^p = H\}$.
Also because $A(x) > \delta_H$ outside some compact set, $A \geq 0$ and
$A(y) > 0$ if and only if $H^p(y) > 0$,
there is a constant $\epsilon_A>0$ with the property that $A^{-1}([0,\epsilon_A])$
is contained in the region $(H^p)^{-1}([0,\nu_H])$.
For $\lambda \geq 1$, any $1$-periodic orbit of $\lambda H^p$ of action
$\leq \lambda \epsilon_A$ must be contained in the region $(H^p)^{-1}([0,\nu_H])$.
This orbit is also contained inside ${\mathcal N}$ and the region where $H^p$
is equal to $H$.

Choose a sequence of perturbations $K_\lambda^i$
that $C^2$ converge to $\lambda H$ where the number of $1$-periodic orbits of $K_\lambda^i$
contained inside ${\mathcal N}$
is bounded above by $P(\lambda)$ and so that all of these
orbits are non-degenerate.
We make the perturbation small enough so that
orbits of action less than or equal to $\lambda \epsilon_A$
are contained in some fixed open subset $U$ of $\{H^p = H\} \cap {\mathcal N}$
such that $U$ contains $(H^p)^{-1}([0,\epsilon])$.
Let $\rho : W \rightarrow \R$ be a bump function which is $0$
on a small neighbourhood of the closure of $U$ and $1$ in the region where $H^p \neq H$.
Let $H_\lambda^i := (1-\rho) K_\lambda^i + \rho \lambda H^p$.
These Hamiltonians $C^1$ converge to $\lambda H^p$.
Hence for large enough $i$, we have that the action of all the $1$-periodic
orbits of $H_\lambda^i$ that are not entirely contained in $U$
is greater than $\frac{\lambda \epsilon_A}{2}$.
For large enough $i$ we also have a sequence of constants $\delta_i$
tending to $0$ such that $(H_\lambda^i)^{-1}([0,\lambda \epsilon-\delta_i])$
contains all the orbits of action $\leq \frac{\lambda \epsilon_A}{2}$
and so that there are no orbits on the boundary of
$U_i := (H_\lambda^i)^{-1}([0,\lambda \epsilon_A-\delta_i])$.
We can also assume that $U_i \subset U$.

By Lemma \ref{lemma:perturbinghamiltonian},
we perturb $H_\lambda^i$ again to $H_\lambda^{'i}$ by a $C^2$ small amount outside the closed set
$U_i$ so that $H_\lambda^{'i} = K_\lambda^i$ inside $U_i$,
all of its $1$-periodic orbits outside $U_i$ are non-degenerate
of action greater than $\frac{\lambda \epsilon_A}{2}$
and such that it is equal to $H_\lambda^i$ near infinity.
For large enough $i$, $(H_\lambda^{'i},J^p)$
is growth rate admissible where $H_\lambda^{'i}$ has only non-degenerate orbits and
such that all of its orbits of action less than or equal to
$\frac{\lambda \epsilon_A}{2}$
are in the region where this Hamiltonian is equal to $K_\lambda^i$.
Hence the number of $1$-periodic orbits of $H_\lambda^{'i}$
of action less than or equal to $\frac{\lambda \epsilon_A}{2}$
is bounded above by $P(\lambda)$.
This implies that the rank of
$V_\lambda := SH_*^{[0,\frac{\lambda \epsilon_A}{2}]}(H_\lambda^{'i},J^p)$ is bounded
above by $P(\lambda)$.
By Lemma \ref{lemma:growthrateactionbound},
the filtered directed system
$(V_\lambda)$
is isomorphic to $SH_*^{\#}(W,\theta_W)$.
This proves the Lemma.
\qed

Suppose we have a Hamiltonian $H$ and a function $P$
such that there is a small neighbourhood ${\mathcal N}$
of $H^{-1}(0)$ where $H,P,{\mathcal N}$
satisfy the properties
of the previous Lemma, then we say that $H$ is {\it $P$ bounded}.

\begin{lemma} \label{lemma:boundednessinvariance}
Suppose that $(W,\theta_W)$ and $(W',\theta_{W'})$
are finite type convex symplectic manifolds that
are convex deformation equivalent at infinity.
Suppose also that we have a Hamiltonian $H$ on $W$
and a function $P : \R \rightarrow \R$ such that:
\begin{enumerate}
\item $H$ satisfies the Liouville vector field property.
\item $H$ is $P$ bounded.
\end{enumerate}
Then $W'$ also admits a Hamiltonian $H'$
satisfying the Liouville vector field property and that is $P$ bounded.
\end{lemma}
\proof of Lemma \ref{lemma:boundednessinvariance}.
Let $X_{\theta_W}$ be the $d\theta_W$-dual of $\theta_W$.
Let $f_W$ be an exhausting function such that $df_W(X_{\theta_W}) > 0$
when $f_W \geq C_W$ for some $C_W$.
By Lemma \ref{lemma:makingthingscomplete},
there is a family of $1$-forms $\theta_W^s$ such that:
\begin{enumerate}
\item $d\theta_W^s$ is symplectic.
\item $\theta_W = \theta_W^0$.
\item $\theta_W^s = \theta_W$ in the region $f^{-1}(-\infty,C_W]$.
\item $\theta_W^1$ is complete.
\item The $d\theta_W^s$-dual $X_{\theta_W^s}$ of $\theta_W^s$ 
satisfies $X_{\theta_W^s} = g_sX_{\theta_W^0}$ for some positive function $g_s$
which smoothly depends on $s$.
\end{enumerate}
We have a similar family of $1$-forms $\theta_{W'}^s$ on $W'$
and a similar function $f_{W'}$ such that these $1$-forms
are all equal on $f_{W'}^{-1}(-\infty,C_{W'}]$.
The pairs $(W,\theta_W^s)$ (resp. $(W',\theta_{W'}^s)$)
are all finite type convex symplectic manifolds because
their Liouville vector fields are all of the form $g_s X_{\theta_W}$
(resp. $h_s X_{\theta_{W'}}$) for some smooth family of functions $g_s$ (resp.
$h_s$)
and the same reason ensures that this is a convex deformation.
We can assume that the closure  of $H^{-1}(0)$
is contained in $f_W^{-1}(-\infty,C_W)$.
Hence $\theta_W = \theta_W^s$ on some compact set whose interior
contains $H^{-1}(0)$.
This implies that $H$ satisfies the Liouville vector field property
and is $P$ bounded on $(W,\theta_W^s)$ for all $s$
and in particular for $(W,\theta_W^1)$.
Let $X_{\theta_W^1}$ be the $d\theta_W^1$-dual of $\theta_W^1$.
If $\Phi_t$ is the flow of this vector field then $\Phi_t^*\theta_W^1 = e^t \theta_W^1$
which implies that $H_t := \Phi_t^*(H)$ is still $P$ bounded.
In particular we can ensure that the compact set $H_t^{-1}(0)$
is arbitrarily large.

Because $(W,\theta_W^1)$ and $(W',\theta_{W'}^1)$
are complete finite type convex symplectic manifolds that
are convex deformation equivalent at infinity,
we have by Lemma \ref{lemma:convexdeformationequivalentatinfinitycomplete}
that they are exact symplectomorphic at infinity.
Let $\phi : W \setminus K_W \rightarrow W' \setminus K_{W'}$
be this exact symplectomorphism.
By the previous discussion, we can assume that the interior of $H_t^{-1}(0)$
contains $K_W$.
We define $K'$ on $W'$ by pushing $H_t$ forward via $\phi$
and then defining $K'$ to be zero inside $K_{W'}$.
This Hamiltonian still satisfies the Liouville vector field property
and is $P$ bounded on $(W',\theta_{W'}^1)$.
Let $\Phi'_t$ be the flow of the vector field
$X_{\theta_{W'}^1}$.
For large enough $T$, $H' := (\Phi'_T)^* K'$
has the property that $(H')^{-1}(0)$ is a subset of
$f_{W'}^{-1}(-\infty,C_{W'})$.
Because $\theta_W^s = \theta_W$ inside $f_{W'}^{-1}(-\infty,C_{W'})$,
we have that $H'$ is $P$ bounded on $(W',\theta_{W'}^s)$
for all $s$ and in particular on $(W',\theta_{W'})$.
\qed

\bigskip
Given some smooth affine variety $A$, we will construct a Hamiltonian $H$ on
$A$ satisfying the Liouville vector field property.
By Theorem \ref{theorem:affinedeforming},
the smooth affine variety $A$ is exact symplectomorphic
to a finite type convex symplectic manifold $W$ described as follows:
We start with a compact symplectic manifold
$M$ and a of set of transversely intersecting
symplectic submanifolds $S_i$ ($i=1,\cdots,k$).
We write $S_{\bracket{i_1,\cdots,i_l}}$ for the intersection
$\cap_{j} S_{i_j}$.
There are small neighborhoods $US_i$ of $S_i$
and projections $\pi_i : US_i \twoheadrightarrow S_i$ and
a positive integer $n_W$ (the number of such $S_i$'s) such that
\begin{enumerate}
\item For $1 \leq i_1 < i_2 < \cdots < i_l \leq n_W$,
\[\pi_{i_l} \circ \cdots \circ \pi_{i_1} : \cap_{j = 1}^l US_{i_j} \twoheadrightarrow S_{\{i_1,\cdots,i_l\}}\]
has fibers that are symplectomorphic to $\Pi_{j = 1}^l \D_{\epsilon}$
where $\D_\epsilon$ is the disk of radius $\epsilon$.
\item
If we look at a fiber $\Pi_{j = 1}^l \D_{\epsilon}$ of
$\pi_{i_l} \circ \cdots \circ \pi_{i_1}$, then for $1 \leq m \leq l$,
$\pi_{i_m}$ maps this fiber to itself.
It is equal to the natural projection
\[\Pi_{j = 1}^l \D_{\epsilon} \twoheadrightarrow  \Pi_{j = 1, j \neq m}^l \D_{\epsilon}\]
eliminating the $m$th disk $\D_\epsilon$.
\item
The symplectic structure on $US_i$ induces a natural connection for
$\pi_{i_l} \circ \cdots \circ \pi_{i_1}$ given by the $\omega$ orthogonal
vector bundle to the fibers.
We require the associated parallel transport maps to be elements
of $U(1) \times \cdots \times U(1)$ where $U(1)$ acts on the disk $\D_\epsilon$
by rotation.
\end{enumerate}
The finite type convex symplectic manifold $W$ we want is symplectomorphic
to $M \setminus \cup_i S_i$.
Let $r_i : US_i \rightarrow \R$ be the function such that when
we restrict $r_i$ to a fiber $D_\epsilon$ of $\pi_i$,
we get that $r_i$ is the distance from the origin.
This uniquely determines $r_i$ because $\pi_i$ has a $U(1)$ structure
group.

We will now construct our Hamiltonian
on the Liouville manifold $W$.
Let $\nu : [0,\epsilon) \rightarrow [0,\infty)$
be a smooth function
satisfying:
\begin{enumerate}
\item Near $t = \epsilon$,
$\nu(t) = 0$ and near $0$, $\nu(t) = \frac{\epsilon^2}{4} - t^2$.
\item $\nu'<0$ when $\nu > 0$.
\item There is one point $x$ where $\nu''(x)=0$ and $\nu(x)>0$
(i.e. the graph of $\nu$ has exactly one point of inflection).
In particular, $\nu''(x)$ is negative when $\nu(x)>0$ is small.
\item If $\nu(x)>0$ is small then we require that
$x < -\kappa_i$ for each $i$ and $x < 1$.
\end{enumerate}

\begin{figure}[H]
\centerline{
 \scalebox{1.0}{
  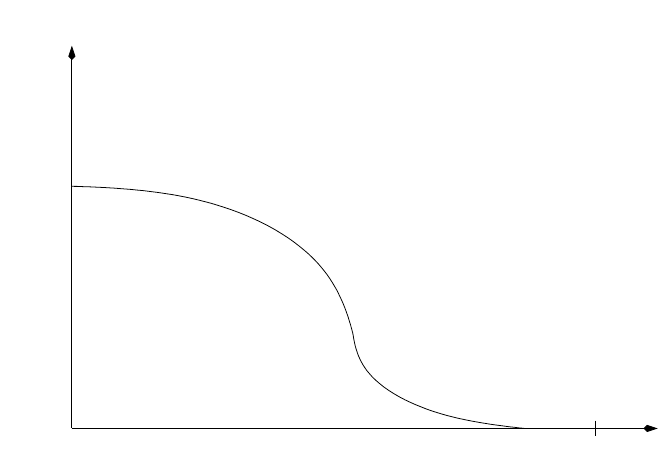
  }
  }
 \end{figure}

The function $\nu(r_i) : US_i \rightarrow \R$
can be extended smoothly by $0$ to the whole of $M$.
From now on we will write $\nu(r_i)$
as this function from $M$ to $\R$.
Let $H := \sum_i \nu(r_i)$.
Theorem \ref{theorem:affinedeforming} also says that
the Liouville form $\theta_W$ on $W$ satisfies:
$\theta_W$ restricted to a fiber
$\D_\epsilon$ of $\pi_i$ is equal to $(r_i^2 + \kappa_i) d\vartheta_i$
where $r_i,\vartheta_i$ are polar coordinates on $\D_\epsilon$ and $\kappa_i < 0$
is some constant.
Because the level sets of $\nu(r_i)$ contain
the $\omega_A$-orthogonal vector bundle to the fibers of $\pi_i$,
we have that its Hamiltonian flow is contained in the fibers
and is equal to $-\frac{1}{2r_i} \nu'(r_i) \frac{\partial}{\partial \vartheta_i}$.
We have that $\theta_W(X_H) = -\frac{r_i^2 + \kappa_i}{2r_i} \nu'(r_i)$.
Hence if $X_{\theta_W}$ is the $d\theta_W$-dual of  $\theta_W$
then $dH(X_{\theta_W}) = -\theta_W(X_H) \geq 0$
and $dH(X_{\theta_W}) = -\theta_W(X_H) > 0$ when $H > 0$.
The derivative with respect to $r_i$ of $\theta_W(X_H)$
is $\left(- \frac{1}{2} + \frac{\kappa_i}{2r_i^2}\right) \nu'(r_i)
-\frac{r_i^2 + \kappa_i}{2r_i} \nu''(r_i)$
which is negative when $\nu(r_i)>0$
is small because $\nu''(r_i) \gg \nu'(r_i)$ and $\kappa_i$ is negative.
Hence $d\left(dH(X_{\theta_W})\right)(X_{\theta_W})
= \theta_W(X_{\theta_W(X_H)})>0$.
Hence $H$ satisfies the Liouville vector field property.
We will call $H$ a {\it Hamiltonian compatible with a compactification of $A$}.

\begin{lemma} \label{lemma:showingHisPbounded}
The Hamiltonian $H$ is $P$ bounded for some polynomial
of degree at least $d$ where
\[d = n - \text{min}(\frac{1}{2}\text{dim}_\R \{S_I | S_I \neq \emptyset\} ).\]
\end{lemma}
\proof of Lemma \ref{lemma:showingHisPbounded}.
The basic idea here is that the orbits of $\lambda H$ form
manifolds with corners coming from the divisors.
The number of such manifolds is bounded above by some constant times $\lambda^d$
and the number of such manifolds up to diffeomorphism is finite.
We use a Morse function on each of these manifolds with corners
to perturb them into non-degenerate orbits and this gives us our bound.
From now on we will assume that $W$ is connected.

We write $US_I$ for $\cap_{i \in I} US_i$.
The Hamiltonian flow of $\lambda H$
in the region $US_I \setminus \cup_{J \frac{\subset}{\neq} I} US_J$
fixes the fibers $\prod_{i \in I} \D_\epsilon$ and for each $i \in I$ it rotates the
$i$'th fiber in this product
by an angle of $-\lambda \frac{1}{2r_i} \nu'(r_i)$.
Hence the fixed points of the time-$\lambda$ flow $H$ form manifolds
with corners.
These are diffeomorphic to $\cap_{i \in I} US_i$ minus the union of
the interiors of $US_i$ for all $i \notin I$.
%
%
There is also one such connected manifold of co-dimension $0$ which is diffeomorphic
to $W$ minus the union of
the interiors of $US_i$ for all $i$.
Let $\tau$ be a constant greater than the supremum
of $-\frac{4\nu'(t)}{t}$. This is bounded because
$\nu'(t)$ is a multiple of $t$ near $0$.
There are at most
$\lfloor \tau \lambda / 2\pi \rfloor^m \binom{n_w}{m}$ such manifolds of codimension $m$
for $1 \leq m \leq d$ and $0$ such manifolds for $m > d$.
Also there are at most $\sum_{m=0}^d  \binom{n_W}{m}$
diffeomorphism types of such manifolds.

Let $h : M \rightarrow \R$
be a function such that $h$ restricted to
$\cap_j S_{i_j}$ is Morse for any $i_1<i_2<\cdots<i_l$.
We also assume it has the following properties:
\begin{enumerate}
\item $h$ is a $C^{\infty}$ small
perturbation of $H$.
\item 
On a small neighbourhood of the closure
of $\nu(r_i) > 0$, we have that $h = \nu(r_i(1 - \delta_h))$
for some small $\delta_h>0$.
\end{enumerate}
We now perturb $\lambda H$ as follows:
Consider a manifold $Y$ of fixed points of the Hamiltonian flow
of $\lambda H$ in the region where
there is an $i_1 < i_2 < \cdots < i_l$ such that
$\nu'(r_k)\neq 0$ if and only if $k = i_j$
for some $j$.
This manifold is contained in the region $\cap_j US_{i_j}$.
Let $N_Y$ be a small neighbourhood of $Y$
such that the only fixed points of $\lambda H$
inside $N_Y$ are ones in $Y$.
The manifold is of co-dimension $l$ and is a co-dimension $0$ submanifold with corners
inside the closed manifold $\overline{Y} := \cap_{j=1}^l \{r_{i_j} = s_j\}$ for some $s_1,\cdots,s_l \in \R$.
In fact $Y = \overline{Y} \cap \cap_{k \neq i_j \forall j} \{\nu(r_k) = 0\}$.
Let $\rho : N_Y \rightarrow \R$ be a bump function which is $0$
outside some compact set in $W$ and $1$ near $Y$.
Let \[\tilde{h} := h|_{S_{i_1,\cdots,i_l}} \circ \pi_{i_1} \circ \cdots \circ \pi_{i_j}.\]
We perturb $\lambda H$ to
\[\tilde{H} := \lambda H +
\delta_Y \rho \tilde{h}.\]
For $\delta_Y$ sufficiently small, there are no fixed points in the region
where the derivative on $\rho$ is non-zero by a compactness argument.
The point is that if there were such orbits for $\delta_Y$ small,
then there would a sequence of $\delta_Y$'s converging to zero
with a fixed point in the region where $d\rho \neq 0$. But this would imply
by compactness that $\lambda H$ would also have such an orbit which is a contradiction.
There are no $1$-periodic orbits of $\tilde{H}$ near the
boundary of $Y$ for the following reason:
%
We have that $Y$ is a manifold with corners.
Near a codimension $k \geq 1$ corner of $Y$,
\[H = \nu(r_{i_1}) + \cdots + \nu(r_{i_l}) + \nu(r_{j_1}) +
\cdots + \nu(r_{j_k})\]
for some $j_1 \cdots j_k$ not equal to any of $i_1,\cdots i_k$.
Hence in this region
$\tilde{H}$ is equal to
\[\nu(r_{i_1}) + \cdots + \nu(r_{i_l}) +
\nu(r_{j_1}) + \cdots + \nu(r_{j_k}) + \]
\[\delta_Y (\nu(r_{i_1}-\delta_h) + \cdots + \nu(r_{i_l}-\delta_h) +
\nu(r_{j_1}(1-\delta_h)) + \cdots + \nu(r_{j_k}(1-\delta_h))).\]
Let $J := \{i_1,\cdots,i_l,j_1,\cdots j_k\}$.
Near this codimension $k$ corner, the Hamiltonian $\tilde{H}$
preserves the fibers $\prod_{j \in J} \D_\epsilon$ of:
\[\pi_{j_1} \circ \cdots \circ \pi_{j_k} \circ
\pi_{i_1} \circ \cdots \circ \pi_{i_l}.\]
These Hamiltonians also split up as a sum:
\[\sum_{j \in J} \nu(r_j) + \delta_Y(\nu(r_j(1-\delta_h)))\]
on this fiber $\prod_{j \in J} \D_\epsilon$.
Near the boundary of the closure of $\nu(r_j)>0$,
$\nu(r_j) + \delta_Y(\nu(r_j(1-\delta_h)))$ has no $1$-periodic orbits
for $\delta_Y,\delta_h$ sufficiently small.
Hence $\tilde{H}$ has no $1$-periodic orbits near
this codimension $k$ corner of $Y$.
This implies that $\tilde{H}$ has no $1$-periodic
orbits on a small neighbourhood of $\partial Y$.

%
Hence all the fixed points of $\tilde{H}$
are contained in the region where $\nu(r_k) = 0$
for all such $k$ satisfying $k \neq i_j \forall j$.
In this region, we have that the Hamiltonians
$\tilde{h}$ and $\lambda H$
Poisson commute because the Hamiltonian vector field
of $\tilde{h}$ is contained in the horizontal plane
distribution of the fibration $\pi_{i_1} \circ \cdots \circ \pi_{i_l}$
which is contained in the level set of $\lambda H$.
Let $p \in Y$ be a fixed point of $\tilde{H}$.
The symplectic tangent bundle splits into $\R^{2l} \oplus \R^{2n-2l}$.
Here $\R^{2l}$ is the tangent space of the fiber of
$\pi_{i_1} \circ \cdots \circ \pi_{i_j}$ and
$\R^{2n-2l}$ is the symplectic complement (this is the
Horizontal plane bundle coming from the natural connection
on this fibration).
The linearized return map restricted to $\R^{2l}$
is the same as the linearized return map
of the Hamiltonian $\sum_{j=1}^l \nu(r_{i_j})$
restricted one of the fibers of $\pi_{i_1} \circ \cdots \circ \pi_{i_j}$.
Because this fiber is symplectomorphic to a product of disks,
we get that the flow is equal to the flow of
$\sum_{j=1}^l \nu(R_j)$
on $(\D_\epsilon)^l$ where $R_j$ is the $j$'th
radial coordinate on $\D_{\epsilon}$.
Because this Hamiltonian system splits as a product of autonomous
non-degenerate Hamiltonians, we have that this Linearized return map
has an $l$ dimensional eigenspace with eigenvalue $1$. This eigenspace
is spanned by the Hamiltonian vectors $X_{\nu(R_j)}$ $j = 1,\cdots,l$ at this point.
The linearized return map of $\tilde{h}$
when restricted to the symplectic complement $\R^{2n-2l}$ is
conjugate to the linearized return map of the Hamiltonian
$\delta_Y h|_{\cap_j S_{i_j}}$ inside $\cap_j S_{i_j}$.
Because this is a $C^2$ small Morse function (for $\delta_Y>0$ sufficiently small),
this has no eigenvalue equal to $1$.
This implies that the linearized return map at $p$ of the autonomous Hamiltonian
$\tilde{H}$ has a $2l$ dimensional
eigenspace with eigenvalue $1$ which is spanned by the Hamiltonian vectors $X_{\nu(R_j)}$ at $p$.
Also the number of critical points created from the manifold $Y$ of fixed
points is at most the number of critical points of
$ h|_{\cap_j S_{i_j}}$.
If we perturb all such manifolds of fixed points of
$\lambda H$  we get an non-degenerate
autonomous
Hamiltonian with at most $ C \lambda^d$ $(S^1)^l$ families of fixed points
where $C$ is some constant greater than
the sum over all strata $\cap_j S_{i_j}$
of the number of critical points of $h|_{\cap_j S_{i_j}}$
multiplied by $2\pi . 2^{n_W}$.
We will also perturb the fixed points where $\lambda H = 0$
so that it becomes non-degenerate in this region.
Let $M_H$ be the number of such fixed points (such a number can be independent of $\lambda$).

Finally using work from \cite{CieliebakFloerHoferWysocki:SymhomIIApplications},
we perturb each $(S^1)^l$ family of orbits away from the divisors $\cup_i S_i$ into $2^l$ orbits creating
a Hamiltonian
$\varrho : S^1 \times W \rightarrow \R$
where all of its orbits are non-degenerate.
Really \cite{CieliebakFloerHoferWysocki:SymhomIIApplications} deals with $S^1$ families
of periodic orbits and not $(S^1)^l$ families, but because our Hamlitonian $\tilde{H}$
restricted to each fiber $\R^{2l}$ splits up as a sum $\sum_j \nu(R_j)$ we can perturb
$\nu(R_j)$ first (if we view such a function as a function on one of the $\R^2$ factors of $\R^{2l}$).
Hence the number of orbits of $\varrho$ is bounded above by $2^{\text{dim}_\R A}C \lambda^d + M_H$.
This implies that $H$ is $P$ bounded where $P$ is the polynomial
$2C\lambda^d$.
\qed

\proof of Theorem \ref{thm:polynomialgrowthrate}.
From the discussion before Lemma \ref{lemma:showingHisPbounded}
we can find a Hamiltonian $H$ which is compatible with a compactification of $A$.
We can choose an appropriate compactification $M$ of $A$ so that the constant $d$
in Lemma \ref{lemma:showingHisPbounded} is equal to $m_A$.
Hence by this lemma we have that $H$ satisfies the Liouville vector field property and is
$P$ bounded where $P$ is a polynomial of degree at most $m_A$.

Because $A,B$ are convex deformation equivalent at infinity,
by Lemma \ref{lemma:boundednessinvariance} we have that $B$
also admits a Hamiltonian $H'$
that satisfies the Liouville vector field property and is $P$ bounded.
By Lemma \ref{lemma:hamiltoniangrowthratebound} we then have
that $(SH_*^{\#}(B,\theta_B))$ is isomorphic to a filtered directed system
$(V_\lambda)$ where $|V_\lambda| \leq P(\lambda)$.
This proves the Theorem.
\qed

\section{Finite covers of smooth affine varieties}
\label{section:coversofsmoothaffinevarieties}

In this section we will prove Theorem \ref{thm:coverfillingobstruction}.
We will prove the following which has Theorem \ref{thm:coverfillingobstruction}
as a direct corollary:
\begin{theorem} \label{thm:coverpolynomialgrowthrate}
Suppose that a finite type convex symplectic manifold $B$
and a finite cover $\widetilde{A}$ of a smooth affine variety $A$ are convex deformation equivalent at infinity.
Then filtered directed system $(SH_*^{\#}(B,\theta_B,\lambda))$
is isomorphic to a filtered directed system
$(V_x)$ where the dimension $|V_x|$ satisfies
$|V_x| \leq Q(x)$ for some polynomial $Q$.
The degree of this polynomial is less than or equal to $m_A$
where $m_A$ is defined in the introduction.
\end{theorem}
\proof of Theorem \ref{thm:coverpolynomialgrowthrate}.
By looking at the proof of Theorem \ref{thm:polynomialgrowthrate}
we see that $A$ admits a Hamiltonian $H$ which satisfies the Liouville
vector field property and which is $P$ bounded
for some polynomial $P$ of degree $m_A$.

Let $p : \widetilde{A} \twoheadrightarrow A$ be the covering
of degree $k$.
We can pull back this Hamiltonian $H$ to $\widetilde{H}$.
This also satisfies the Liouville vector field property
and is also $kP$ bounded.
It satisfies the Liouville vector field property because
we can lift any Liouville vector field on $A$ to a Liouville vector field
on $\widetilde{A}$.
It is $kP$ bounded for the following reason:
If $H'$ is any Hamiltonian and $\widetilde{H'}$ its lift to $\widetilde{A}$
then any orbit of $\widetilde{H'}$ projects down to an orbit of $H'$.
This orbit is non-degenerate if and only if the projected orbit
is non-degenerate.
There are at most $k$ orbits of $\widetilde{H'}$ that project to a given
orbit of $H$.
Hence the number of non-degenerate orbits of $\widetilde{H'}$
is bounded above by $k$ times the number of non-degenerate orbits of $H'$.
Also the action of this orbit is equal to the action of the projected orbit.

Hence by Lemma \ref{lemma:boundednessinvariance},
$B$ admits a Hamiltonian which is $kP$ bounded.
Hence by Lemma \ref{lemma:hamiltoniangrowthratebound}
we get that $(SH_*^{\#}(B,\theta_B,\lambda))$
is isomorphic to
$(V_x)$ where the dimension $|V_x|$ satisfies
$|V_x| \leq kP(x)$.
Hence we have proven this theorem where $Q = kP$.
\qed

\section{Appendix A: convex symplectic manifolds}

Recall that a convex symplectic manifold is a manifold $M$
together with a $1$-form $\theta$
\begin{enumerate}
\item $\omega := d\theta$ is a symplectic form.
\item The $\omega$-dual $X_{\theta}$ of $\theta$ is a vector field satisfying
$df_M(X_{\theta}) > 0$ in the region $f_M^{-1}(A_M)$ where $A_M \subset \R$ is an unbounded subset
and $f_M : M \rightarrow \R$ an exhausting function.
\end{enumerate}
Exhausting means that $f_M$ is bounded from below and the preimage
of every compact set is compact.

We now need to describe various ways of deforming $(M,\theta)$.
The first kind of deformation is the most restrictive:
Two convex symplectic manifolds $(M,\theta),(M',\theta')$
are {\it strongly deformation equivalent} if there is
a diffeomorphism $\phi : M \rightarrow M'$, an exhausting function
$g : M \rightarrow \R$ and a smooth family
of $1$-forms $\theta_t$ ($t \in [0,1]$) on $M$
such that:
\begin{enumerate}
\item $\omega_t := d\theta_t$  is symplectic.
\item If $X_{\theta_t}$ is the $\omega_t$-dual of $\theta_t$ then
$dg(X_{\theta_t}) > 0$ on $g^{-1}(A_g)$ where $A_g \subset [0,\infty)$ is unbounded.
\item $\theta_0 = \theta$ and $\theta_1 = \phi^* \theta'$.
\end{enumerate}
Two convex symplectic manifolds $(M,\theta),(M',\theta')$
are {\it convex deformation equivalent} if there is a finite sequence
$(M_1,\theta_1),\cdots,(M_k,\theta_k)$ such that
\[(M_1,\theta_1) = (M,\theta), (M_k,\theta_k) = (M',\theta')\]
and such that for each $i<k$, $(M_i,\theta_i)$ is strongly deformation
equivalent to $(M_{i+1},\theta_{i+1})$.

An important class of convex symplectic manifolds are the complete ones.
A {\it complete convex symplectic manifold} is a convex symplectic
manifold such that the associated Liouville vector field
$V$ (the $\omega$-dual of $\theta$) is integrable for all time.
A (strong) convex deformation $(M,\theta_t)$ is called {\it complete}
if the associated Liouville vector fields $X_{\theta_t}$ are integrable.

\begin{lemma} \label{lemma:makingthingscomplete}
Let $(M,\theta)$ be a convex symplectic manifold
and let $K \subset M$ be any compact set.
There exists a smooth family of $1$-forms $\theta_t$,($t \in [0,1]$)
on $M$ such that:
\begin{enumerate}
\item $\theta_0 = \theta$
\item $\omega_t := d\theta_t$ is a symplectic form.
\item If $X_{\theta_t}$ is the $\omega_t$-dual of $\theta_t$, then there is
a smooth family of functions $f_t : M \rightarrow \R$ such that $X_{\theta_t} = f_t X_{\theta_0}$
and $0 < f_t \leq 2$.
\item The vector field $X_{\theta_1}$ is integrable.
\item $\theta_t|_K = \theta_0|_K$.
\end{enumerate}
\end{lemma}
\proof of Lemma \ref{lemma:makingthingscomplete}.
Let $X_{\theta}$ be the $d\theta$-dual of $\theta$.
Let $A \subset [0,\infty)$ be an unbounded subset and $f$ an exhausting function
such that $df(X_{\theta})>0$ on $f^{-1}(A)$.
Note that each $a  \in A$ is a regular value of $f$ because
$df(X_{\theta})>0$ on $f^{-1}(a)$ which implies in particular that $df \neq 0$
on $f^{-1}(a)$.
Because $f$ is smooth, there exists a sequence $c_1 < c_2,\cdots$
tending to infinity such that $c_i$ is in $A$
and $X_{\theta}$ is transverse to $f^{-1}(c_i)$ and pointing outwards.
We choose $c_i$ large enough so that $f^{-1}(-\infty,c_i)$ contains $K$.
By flowing along $X_{\theta}$
there is a neighbourhood of the manifold  $C_i := f^{-1}(c_i)$
diffeomorphic to $[1-\epsilon_i,1+\epsilon_i] \times C_i$
with $\theta = r_i \alpha_i$ where $r_i$ parameterizes the interval
and $\alpha_i = \theta |_{C_i}$.
Again we make these neighborhoods small enough so that they are disjoint from $K$.
Let ${\mathcal N}_i$ be this neighbourhood and let ${\mathcal N}$
be the union of all these neighborhoods.
We can make these neighborhoods small enough so that they are disjoint
and also so that $\epsilon_i \leq \frac{1}{2}$.

We will now construct $\theta_t$:
Let $g_i : [1-\epsilon_i,1+\epsilon_i] \rightarrow \R$ be smooth functions
such that:
\begin{enumerate}
\item $g_i' > 0$.
\item In the region $[1-\epsilon_i,1-2\epsilon_i/3]$, we have $g_i(r_i) = r_i$.
\item In the region $[1-\epsilon_i/2,1+\epsilon_i/2]$ we have $g_i/g'_i \leq \epsilon_i$.
\item $g_i/g_i' \leq 1$
or equivalently $\log(g_i)' \geq 1$.
\item In the region $[1+2\epsilon_i/3,1+\epsilon_i]$ we have $g_i(r_i) = \kappa_i r_i$
where $\kappa_i \geq 1$ is some large constant.
\end{enumerate}
We define $g_i^t(r_i) := (1-t)r_i + tg_i(r_i)$.
We define $c_{0} := -\infty$.
Let $\Xi_i := \prod_{j=1}^{i-1} \kappa_i$ for $i > 1$
and $\Xi_1 := 1$.
In the region $f^{-1}(c_{i-1},c_i) \setminus {\mathcal N}$,
we define $\theta_t := ((1-t) + t\Xi_i) \theta$.
In the region ${\mathcal N}_i$ we define $\theta_t$
as $((1-t) + t\Xi_i) g_i^t(r_i) \alpha_i$.
This definition ensures that $\theta_t$ ($t \in [0,1]$)
is a smooth family of $1$-forms.

Outside ${\mathcal N}_i$, we have that $\theta_t$
is equal to some locally constant function multiplied by
$\theta$, hence $d\theta_t$ is still symplectic.
In the region ${\mathcal N}_i$, we have that
$d\theta_t = ((1-t) + t\Xi_i)\left( (g')_i^t dr_i \wedge \alpha_i + r_i d\alpha_i\right)$.
This is symplectic because $\alpha_i$ is a contact form
and ${g_i^t}' > 0$.
Hence $d\theta_t$ is symplectic for all $t$.
Outside ${\mathcal N}_i$ we have that the $\omega_t$-dual $X_{\theta_t}$
of $\theta_t$ is equal to $V$ because rescaling $\theta$ by a locally constant function
does not change the associated dual vector field.
Inside ${\mathcal N}_i$ we have that
$X_{\theta_t} = g_i^t / {g_i^t}' \frac{\partial}{\partial r_i}$
which is equal to $V = r_i \frac{\partial}{\partial r_i}$
multiplied by some positive function $f_t$.
Because $g_i / g_i' \leq \epsilon_i$ in the region
$[1-\epsilon_i/2,1+\epsilon_i/2]$, if we flow any point $p$ in the
region where $r_i = 1-\epsilon_i/2$ along the vector field
$X_{\theta_1}$ for time $1$ to a point $q$, then $q$ is still contained
in $[1-\epsilon_i/2,1+\epsilon_i/2] \times C_i$.
This ensures that the vector field $X_{\theta_1}$ is complete
and hence $(M,\theta_1)$ is complete.
We also have that $f_t \leq 2$ because $f_t = 1$ outside ${\mathcal N}$
and is equal to
$\frac{g_i^t}{{g_i^t}' r_i} \leq \frac{1}{1-\epsilon_i} \leq 2$
inside ${\mathcal N}_i$.
Finally we have that $\theta_0 = \theta$ by definition
and $\theta_t|_K = \theta_0|_K$ because ${\mathcal N}$ is disjoint from $K$.
\qed

We also  have a $1$-parameter version of this Lemma as follows
(with almost exactly the same proof):
Suppose $(M,\theta_s)$ ($s \in [0,1]$) is a smooth family of convex symplectic manifolds
such that we have a function $f : M \rightarrow \R$ and an
unbounded $A \subset [0,\infty)$ with $df(X_{\theta_s} > 0$ on $f^{-1}(A)$
where $X_{\theta_s}$ is the $d\theta_s$-dual of $\theta_s$.
Then there is a two parameter family of $1$-forms $\theta_{s,t}$
$(s,t) \in [0,1]^2$ satisfying:
\begin{enumerate}
\item $\theta_{s,0} = \theta_s$
\item $\omega_{s,t} := d\theta_{s,t}$ is a symplectic form.
\item If $X_{\theta_{s,t}}$ is the $\omega_{s,t}$-dual of $\theta_{s,t}$, then there is
a smooth family of functions $f_{s,t}$ such that $X_{\theta_{s,t}} = f_{s,t} X_{\theta_{s,0}}$
and $0 < f_{s,t} \leq 2$.
\item The vector field $X_{\theta_{s,1}}$ is integrable for each $s \in [0,1]$.
\end{enumerate}

\begin{corollary} \label{cor:completing}
Every convex symplectic manifold
is strongly deformation equivalent to a complete
convex symplectic manifold.
\end{corollary}
\proof
Let $f_M,X_{\theta}$ be the associated function and Liouville
vector field of $M$. There is an unbounded
$A \subset [0,\infty)$ such that $df_M(X_{\theta}) >0$ in $f^{-1}(A)$.
Let $(M,\theta_t)$ be the family of convex symplectic manifolds
described in the Lemma \ref{lemma:makingthingscomplete} above.
If $X_{\theta_t}$ is the associated Liouville vector field
of $(M,\theta_t)$ then because $X_{\theta_t} = f_t X_{\theta}$ for some smooth
family of functions $f_t > 0$, we have that
$df_M(X_{\theta_t}) > 0$ in $f^{-1}(A)$.
Hence we have that $(M,\theta_t)$ is our strong convex deformation.
\qed

\begin{corollary} \label{lemma:finitetypecompleting}
Let $f_M,X_{\theta}$ be the associated function and Liouville
vector field of $M$.
Suppose $M$ is of finite type which means that $df_M(X_{\theta}) > 0$
in the region $f^{-1}[C,\infty)$, then
$M$ is convex deformation equivalent to the completion of
the Liouville domain $f^{-1}(-\infty,C]$.
\end{corollary}
\proof
Let $(M,\theta_t)$ be the family of $1$-forms as in Lemma
\ref{lemma:makingthingscomplete} and let $X_{\theta_t} = f_t X_{\theta}$
be the associated family of Liouville vector fields.
Then $df_M(X_{\theta_1}) > 0$ in the region
$f_M^{-1}[C,\infty)$.
Let $D := f_M^{-1}(-\infty,C]$ and let $\alpha = \theta_1|_{\partial D}$..
By flowing the contact manifold $(\partial D,\alpha)$
along $X_{\theta_1}$ (which is integrable) we obtain a diffeomorphism $\phi$ from
$\partial D \times [1,\infty)$ to $f_M^{-1}[C,\infty)$
such that $\phi^* \theta_1 = r \alpha$ where $r$ parameterizes $[1,\infty)$.
Also by the previous corollary \ref{cor:completing} we have
that $(M,\theta)$ is convex deformation equivalent to $(M,\theta_1)$.
Hence $(M,\theta_1)$ is convex deformation equivalent to the
completion $(\widehat{D},\theta_1)$.
\qed

Suppose we have a convex deformation equivalence $(M,\theta_t)$.
We say that this is a {\it complete convex deformation equivalence}
if the $d\theta_t$-dual $X_{\theta_t}$ of $\theta_t$ is an integrable
vector field for each $t \in [0,1]$.

\begin{lemma} \label{lemma:completingconvexdeformation}
Let $(M,\theta),(M',\theta')$ be complete convex symplectic manifolds
that are convex deformation equivalent.
Then there is a complete convex deformation equivalence between
$(M,\theta)$ and $(M',\theta')$.
\end{lemma}

We need a preliminary Lemma first:
\begin{lemma} \label{lemma:strongohomotopycompleting}
Suppose $(M,\theta),(M',\theta')$ are compete convex
symplectic manifolds that are strongly deformation equivalent.
Then there is a complete strong deformation equivalence between them.
\end{lemma}
\proof of Lemma \ref{lemma:strongohomotopycompleting}.
Let $(M,\theta_s)$ be the strong deformation equivalence.
Let $X_{\theta_s}$ be the associated Liouville vector fields
and let $f : M \rightarrow \R$ and $A \subset [0,\infty)$ (unbounded)
be such that $df(X_{\theta_s}) > 0$ on $f^{-1}(A)$.
By Lemma \ref{lemma:makingthingscomplete}
there is a two parameter family of $1$-forms $\theta_{s,t}$($(s,t) \in [0,1]^2$) such that:
\begin{enumerate}
\item $\theta_{s,0} = \theta_s$
\item $\omega_{s,t} := d\theta_{s,t}$ is a symplectic form.
\item If $X_{\theta_{s,t}}$ is the $\omega_{s,t}$-dual of $\theta_{s,t}$, then there is
a smooth family of functions $f_{s,t}$ such that $X_{\theta_{s,t}} = f_{s,t} X_{\theta_{s,0}}$
and $0 < f_{s,t} \leq 2$.
\item The vector field $X_{\theta_{s,1}}$ is integrable for each $s \in [0,1]$.
\end{enumerate}
Let $p : [0,1] \rightarrow [0,1]^2$
be a smooth path starting at $(0,0)$ and ending at $(1,0)$
whose image is equal to $\{0\} \times [0,1] \cup [0,1] \times \{1\} \cup \{1\} \times [0,1]$.
Then $(M,\theta_{p(t)})$ is our complete strong convex deformation.
This is because $X_{\theta_{s,0}}$ (resp. $X_{\theta_{s,1}}$) is integrable
because it is equal to $f_{s,0} X_{\theta_{0,0}}$ (resp. $f_{s,1} X_{\theta_{0,1}}$)
with $0 < f_{s,t} \leq 2$ for all $s,t$.
This is also because $f_{1,t}$ is integrable for all $t$.
This completes the Lemma.
\qed

\proof of Lemma \ref{lemma:completingconvexdeformation}.
Let $(M,\theta_s)$ be the convex deformation between $(M,\theta)$
and $(M',\theta')$.
This means that we break up $(M,\theta_s)$ into a finite number
of strong deformations
$(M,\theta^1_s),\cdots,(M,\theta^k_s)$.
Let $V^i_s$ be the associated Liouville vector fields.
For a given strong deformation $(M,\theta^i_s)$ there is an unbounded
set $A_i \subset [0,\infty)$ and a function $g_i : M \rightarrow \R$
such that $g_i(V^i_s) > 0$ in $g_i^{-1}(A)$.
By the parameterized version of Lemma \ref{lemma:makingthingscomplete}
we can replace $\theta^i_s$ with $\theta^i_{s,t}$
such that:
\begin{enumerate}
\item $\theta^i_{s,0} = \theta^i_s$
\item $\omega^i_{s,t} := d\theta^i_{s,t}$ is a symplectic form.
\item If $V^i_{s,t}$ is the $\omega^i_{s,t}$-dual of $\theta^i_{s,t}$, then there is
a smooth family of functions $f^i_{s,t}$ such that $V^i_{s,t} = f^i_{s,t} V^i_{s,0}$
and $0 < f^i_{s,t} \leq 2$.
\item The vector field $V^i_{s,1}$ is integrable for each $s \in [0,1]$.
\end{enumerate}
We now replace the deformation $(M,\theta^1_s)$
by the concatenation $(M,{\theta'}^1_s)$ of the homotopies:
$(M,\theta^1_{0,t})$ and
$(M,\theta^1_{t,1})$.
This is still a strong deformation because
$g_1(V^1_{s,t})>0$ in the region $g_1^{-1}(A_1)$.
For $1 < i < k$ we replace $(M,\theta^i_s)$
with the concatenation $(M,{\theta'}^i_s)$ of:
$(M,\theta^{i-1}_{1,1-t})$,
$(M,\theta^i_{0,t})$ and
$(M,\theta^i_{t,1})$.
This again is a strong deformation using the function
$g_i$ and unbounded set $A_i \subset [0,\infty)$.
The reason why $g_i(V^{i-1}_{1,1-t}) > 0$ in $g_i^{-1}(A_i)$
is because $V^{i-1}_{1,1-t})$ is 
$f^{i-1}_{1,1-t} V^i_t$ where $f^{i-1}_{1,1-t}>0$.
Finally we replace the deformation
$(M,\theta^k_s)$ with the concatenation
$(M,{\theta'}^k_s)$ of
$(M,\theta^{k-1}_{1,1-t})$,
$(M,\theta^k_{0,t})$,
$(M,\theta^k_{t,1})$ and
$(M,\theta^k_{1,1-t})$.
Hence we have a new convex deformation equivalence
$(M,{\theta'}^1_t),\cdots,(M,{\theta'}^k_t)$ with the property
that ${\theta'}^i_0,{\theta'}^i_1$ is complete for all $i$.
By Lemma \ref{lemma:strongohomotopycompleting}, we can replace this with
a new convex deformation equivalence:
$(M,{\theta''}^1_t),\cdots,(M,{\theta''}^k_t)$ with the property
that ${\theta''}^i_t$ is complete for all $1 \leq i \leq k, t \in [0,1]$.
This is a complete convex deformation equivalence.
\qed

\begin{corollary} \label{cor:exactsymplectomorphicconvexdeformation}
Suppose that $(M,\theta),(M,\theta')$ are complete convex symplectic
manifolds that are convex deformation equivalent.
Then they are exact symplectomorphic to each other.
\end{corollary}
\proof of \ref{cor:exactsymplectomorphicconvexdeformation}.
By Lemma \ref{lemma:completingconvexdeformation} there is a
complete convex deformation equivalence between
$(M,\theta)$ and $(M,\theta')$.
Hence by \cite[Proposition 12.2]{CieliebakEliashberg:symplecticgeomofsteinmflds}
we get that they are exact symplectomorphic.
\qed

We also need another lemma similar to corollary \ref{cor:exactsymplectomorphicconvexdeformation}
except that we will be dealing with convex deformation equivalence at infinity
as in definition \ref{defn:convexdeformationequivalentatinfinity}.
\begin{lemma} \label{lemma:convexdeformationequivalentatinfinitycomplete}
Suppose that $(W,\theta_W),(W',\theta_{W'})$ are
complete finite type convex symplectic manifolds that
are convex deformation equivalent at infinity.
Then they are exact symplectomorphic at infinity.
\end{lemma}
\proof of Lemma \ref{lemma:convexdeformationequivalentatinfinitycomplete}.
By the definition of convex deformation equivalent at infinity,
we have a sequence of convex symplectic manifolds:
$(Q_i,\theta_{Q_i})$
for $i = 1,\cdots k$ such that:
\begin{enumerate}
\item  $Q_1$ is convex deformation equivalent to $W_1$.
\item $Q_k$ is convex deformation equivalent to $W_2$.
\item $Q_i$ is convex deformation equivalent to $Q_{i+1}$
or they are exact symplectomorphic at infinity.
\end{enumerate}
First of all (by using the fact that the identity
map is an exact symplectomorphism and also a convex deformation
equivalence), we can assume that $k$ is even and
that for odd $i$, that $Q_i,Q_{i+1}$ are convex deformation equivalent
and for even $i$ that $Q_i,Q_{i+1}$ are exact symplectomorphic at infinity.

Suppose we have two convex symplectic manifolds $(A,\theta_A),(B,\theta_B)$
that are exact symplectomorphic at infinity.
We wish to find two complete convex symplectic manifolds $\tilde{A}$
and $\tilde{B}$ that are exact symplectomorphic at infinity
and such that $\tilde{A}$ is convex deformation equivalent to $A$
and $\tilde{B}$ is convex deformation equivalent to $B$.
Let $\phi : A \setminus K_A \rightarrow B \setminus K_B$
be the exact symplectomorphism at infinity where
$K_A,K_B$ are relatively compact sets.
By possibly enlarging $K_A$ a little bit, we can ensure that there
is a function $f: A \rightarrow \R$ such that $\phi^* \theta_B = \theta_A + df$.
Let ${\mathcal N}$ be a small neighbourhood of the closure of $K_A$.
By Lemma \ref{lemma:makingthingscomplete},
there is a family of $1$-forms $\theta_A^s$
such that
\begin{enumerate}
\item $d\theta_A^s$ is symplectic.
\item $\theta_A+df = \theta_A^0$.
\item $\theta_A^s|_{ {\mathcal N} } = \theta_A^0|_{  {\mathcal N} }$
\item $\theta_A^1$ is complete.
\item The $d\theta_A^s$-dual $X_{\theta_A^s}$ of $\theta_A^s$ 
satisfies $X_{\theta_A^s} = gX_{\theta_A^0}$ for some positive function $g$.
\end{enumerate}
We can also define $\theta_B^s$ to be equal to $\theta_B$
near $K_B$ and equal to $\phi_*(\theta_A^s)$ outside $K_B$.

These are all convex symplectic manifolds because
their Liouville vector fields are all of the form $g_s X_{\theta_A + df}$
(resp. $g_s X_{\theta_B}$) for some family of functions $g_s$
and the same reason ensures that this is a convex deformation.
Also by \cite[Lemma 8.3]{McLean:symhomlef}, we have that
$(A,\theta_A)$ is convex deformation equivalent to $(A,\theta_A + df)$.
This means that $(A,\theta_A)$ (resp. $(B,\theta_B)$) is convex deformation equivalent
to the complete convex symplectic manifold $\tilde{A} := (A,\theta_A^1)$
(resp. $\tilde{B} := (B,\theta_B^1)$).
Also $\tilde{A}$,$\tilde{B}$ are exact symplectomorphic at infinity.

The previous discussion ensures (by changing the convex
deformation equivalences from $Q_i$ to $Q_{i+1}$)
that we can assume that the convex symplectic manifolds $Q_i$
are all complete convex symplectic manifolds.
Note we can also assume that $Q_1$ and $Q_k$ are complete
because $W$ and $W'$ are complete.
By Corollary \ref{cor:exactsymplectomorphicconvexdeformation}
this implies that for all odd $i$, $Q_i$ is exact symplectomorphic
to $Q_{i+1}$ and hence in particular they are exact symplectomorphic at infinity.
Because the property of being exact symplectomorphic at infinity is transitive,
we have that $W,W'$ are exact symplectomorphic at infinity.
This proves the Lemma.
\qed

\section{Appendix B: A Maximum Principle}

Let $N$ be a manifold and $\theta$ a $1$-form on $N$ so that
\begin{enumerate}
\item $\omega := d\theta$ is a symplectic form.
\item The $\omega$-dual $X_\theta$ of $\theta$
is transverse to the boundary of $N$ and pointing inwards.
\end{enumerate}
Let $S$ be a compact Riemann surface with boundary and complex structure $j$
and $\gamma$ a $1$-form on $S$.
Let $H : S \times N \rightarrow \R$ be a family of Hamiltonians parameterized by $S$.
We sometimes write this as a family of functions $H_\sigma : N \rightarrow \R$
parameterized by $\sigma \in S$.
Let $J_\sigma$ be a family $\omega$ compatible almost complex structures parameterized by $\sigma \in S$.
A small neighbourhood of $\partial N$ is diffeomorphic to $\partial N \times [1,1+\epsilon_N)$
where $\theta = r \alpha$.
Here $r$ parameterizes the interval $[1,1+\epsilon_N)$ and $\alpha = \theta|_{\partial N}$.
We require that $\theta \circ J_\sigma = dr$ and that $H_\sigma = f(r)$ for some function $f$ with $f(1)=1$ and $f'(1)=1$
near $\partial N$.
The differential $dH$ uniquely splits up as $d_S H + d_N H$
where vectors tangent to $N$ are contained in the kernel of $d_S H$
and vectors tangent to $S$ are in the kernel of $d_N H$.
We can view $d_S H$ as a family of $1$-forms on $S$
parameterized by $N$, so for each $p \in N$ we define
$d_S H_{(\cdot,p)}$ to be $d_S H$ restricted to $S \times \{p\}$.
We require that  $d_S H_{(\cdot,p)} \wedge \gamma + H(\cdot,p)d(\gamma) \geq 0$ for each $p \in N$.

Let $u : S \rightarrow N$ satisfy:
$(du - X_{H_\sigma} \otimes \gamma)^{0,1}_{J_\sigma} = 0$ at each point $\sigma \in S$.
In other words:
\begin{equation} \label{equation:parameterizedcauchyriemann}
du - X_{H_\sigma} \otimes \gamma + J_\sigma \circ (du - X_{H_\sigma} \otimes \gamma) \circ j = 0.
\end{equation}

The aim of this section is to prove:
\begin{lemma} \label{lemma:maximumprinciple}
If $u(\partial S) \subset \partial N$ then $u(S) \subset \partial N$.
\end{lemma}
This lemma is similar to 
\cite[Lemma 7.2]{SeidelAbouzaid:viterbo}.
Before we prove Lemma \ref{lemma:maximumprinciple} we need to define the {\it geometric energy} and {\it topological energy}
of $u$. The geometric energy is defined as:
\[ E^{\text{geom}} := \int_S \left\|du - X_{H_\sigma} \otimes \gamma \right\|^2_{J_\sigma}\]
where $\|\cdot\|_{J_\sigma}$ is the norm coming from the metric $\omega(\cdot,J_\sigma(\cdot))$
and some compatible metric on $(S,j)$.
The topological energy is defined as:
\[ E^{\text{top}} := \int_S u^* \omega + d(u^* H_\sigma \gamma).\]
Here $u^* H_\sigma$ is the function sending $\sigma \in S$ to
$H(\sigma,u(\sigma))$.
We define $u^* d_N H$ to be the $1$-form on $S$ such that for each vector $V$ on $S$,
$(u^* d_N H)(V) = d_N H (V,u_*(V))$.
Similarly we define $(u^* d_S H)(V)$ to be $(d_S H)(V,u_*(V))$.

Let $\sigma$ be a point on $S$ and $s+it$ a local holomorphic chart
around $\sigma$ where $\partial_s,\partial_t$ have magnitude $1$ at $(s,t)=(0,0)$
and where $\sigma$ is the point $(0,0)$.
We have:
\[\| du - X_{H_\sigma} \otimes \gamma \|^2_{J_\sigma} =
\omega(\partial_s u - X_{H_\sigma} \gamma(\partial_s),\partial_t u - X_{H_\sigma} \gamma(\partial_t))\]
\[= u^* \omega(\partial_s,\partial_t) - \gamma(\partial_s) d_N H(\partial_t u)
+ \gamma(\partial_t) d_N H (\partial_s u) \]
\[ = u^* \omega(\partial_s,\partial_t) + (u^*(d_N H) \wedge \gamma)( \partial_s,\partial_t).\]
Also $d_N H(\partial_t u)$ means that we consider $\partial_t u$ as a vector
inside $\{\sigma\} \times N$ and then contract it with $dH$.
The expression $d_N H(\partial_s u)$ has a similar meaning.
Hence $E^{\text{geom}}(u) = \int_S u^* \omega + \int_S u^* d_N H \wedge \gamma$.
We have that
\[d( u^* H_\sigma \gamma ) = u^* d_N H \wedge \gamma + u^* d_S H  \wedge \gamma
+ u^* H_\sigma d\gamma.\]
Hence we have:
\[ E^{\text{geom}}(u) = E^{\text{top}}(u) - \int_S \left(u^* d_S H \wedge \gamma + u^*H_\sigma d\gamma\right)\]
which implies that $E^{\text{geom}}(u) \leq E^{\text{top}}(u)$.

\proof of Lemma \ref{lemma:maximumprinciple}
We suppose that $u(\partial S) \subset \partial N$.
We just need to show that $E^{\text{geom}}(u) = 0$ as this will force our surface
$S$ to map to a Reeb orbit of $\partial N$.
We have that $E^{\text{top}}(u) \geq E^{\text{geom}}(u)$,
so we now need to show that $E^{\text{top}}(u) = 0$.
By Stokes' theorem we have:
\[ E^{\text{top}}(u) = \int_{\partial S} u^* \theta + u^* H_\sigma \gamma .\]
Because $H_\sigma = f(r)$ with $f'(1) = 1$, we have that $-X_{H_\sigma}$ is equal
to the Reeb vector field on $\partial N$.
Hence
$\theta(X_{H_\sigma}) = -1 = -H_\sigma$ along $\partial N$ as $f(1)=1$.
So
\[ E^{\text{top}}(u) = \int_{\partial S} \theta \circ (du - X_{H_\sigma} \otimes \gamma) \]
\[ = \int_{\partial S} \theta \circ J_\sigma \circ (du - X_{H_\sigma} \otimes \gamma) \circ (-j)
 = \int_{\partial S} dr \circ (du - X_{H_\sigma} \otimes \gamma) \circ (-j).\]
Because $dr(X_{H_\sigma}) = 0$ our integral becomes:
\[ E^{\text{top}}(u) = -\int_{\partial S} dr \circ du \circ j .\]
If a vector $V$ is tangent to $\partial S$ and pointing in the direction
in which $\partial S$ is oriented then $j(V)$ points inwards.
This implies that $dr \circ du \circ j(V) \geq 0$ because $r$ increases
as we move towards the interior of $N$.
Hence \[ E^{\text{top}}(u) \leq 0.\]
Hence $E^{\text{geom}}(u)$ vanishes which gives us our result.
\qed

Here are two applications of this lemma:
Let $M$ be a Liouville domain with $1$-form $\theta_M$.
Then its completion $\widehat{M}$ has a cylindrical end
$\partial M \times [1,\infty)$ with cylindrical coordinate $r_M$.
We define
$K : \widehat{M} \rightarrow \R$ to be an autonomous Hamiltonian on $\widehat{M}$
which is equal $k(r_M)$ near $\partial M$ where $k' > 0$.
Let $H_{s,t}$ be a family of Hamiltonians parameterized
by $(s,t) \in \R \times S^1$ so that $H_{s,t} = K+a_s$ near $\partial M$
where $a_s$ is a smooth family of constants.
We require that $\frac{\partial H}{\partial s} \geq \frac{\partial a_s}{\partial s}$.
We define $J_{s,t}$ to be a smooth family of almost complex structures
which are cylindrical near $\partial M$.
\begin{corollary} \label{corollary:parameterizedfloermax}
Suppose that $u_1 : \R \times S^1 \rightarrow \widehat{M}$ satisfies the perturbed
Cauchy-Riemann equations:
\[ \partial_s u_1 + J_{s,t} \partial_t u_1 = J_{s,t} X_{H_{s,t}}\]
and that $u_1(s,t) \in M$ for $|s| \gg 1$.
Then $u_1(s,t) \in M$ for all $(s,t) \in \R \times S^1$.
\end{corollary}
\proof of Corollary \ref{corollary:parameterizedfloermax}.
Let $S$ be equal to $u_1^{-1}(\partial M \times [1,\infty))$.
We perturb $\partial M$ slightly so that $S$ is a codimension
$0$ submanifold with boundary.
We define $H_1 := \frac{H-k(1)-a_s}{k'(1)} +1$.
We let $\gamma = k'(1) dt$ and so $(d_S H_1)_{(\cdot,p)} \wedge \gamma + (H_1)_{(\cdot,p)} d\gamma \geq 0$.
We have that $u_1$ satisfies the perturbed Cauchy Riemann equations
with respect to $H_1$ and $J_{s,t}$.
We also have that $H_1$ is equal to $h(r_M) = \frac{k(r_M)-k(1)}{k'(1)} +1$ near $\partial M$
and so $h(1)=h'(1)=1$.
By using Lemma \ref{lemma:maximumprinciple} with
$N = \partial M \times [1,\infty)$ we have that $u_1(S)$ must be
contained inside $\partial M$.
Hence the image of $u_1$ is contained in $M$.
\qed

\begin{corollary} \label{corollary:maximumprinciplerescaling}
Let $g : \R \rightarrow \R$ be a function satisfying $g'(s) \geq 0$.
Suppose in addition that $H \geq k(1) - k'(1) + a_s$  inside $\partial M \times [1,\infty)$,
then any solution:
$u_2 : \R \times S^1 \rightarrow \widehat{M}$
of
\[ \partial_s u_2 + J_{s,t} \partial_t u_2 = J_{s,t} X_{g(s) H_{s,t}}\]
with $u_2(s,t) \in M$ for $|s| \gg 1$ has image contained in $M$.
\end{corollary}
\proof of \ref{corollary:maximumprinciplerescaling}.
Again we define $H_1 := \frac{H-k(1)-a_s}{k'(1)} +1$.
We let $\gamma = k'(1) g(s) dt$.
Because $H_1, \frac{\partial H_1}{\partial s} \geq 0$,
we have $(d_S H_1)_{(\cdot,p)} \wedge \gamma + (H_1)_{(\cdot,p)} d\gamma \geq 0$.
We then apply Lemma \ref{lemma:maximumprinciple} to give us our result.
\qed

\bibliography{references}

\newcommand{\etalchar}[1]{$^{#1}$}
\begin{thebibliography}{CFHW96}

\bibitem[{Abo}10]{Abouzaid:contangentgenerate}
{Abouzaid, Mohammed}.
\newblock {A cotangent fibre generates the Fukaya category}.
\newblock pages 1--40, 2010, arXiv:1003.4449.

\bibitem[Ano81]{Anosov:homotopiescurves}
D.~V. Anosov.
\newblock Some homotopies in a space of closed curves.
\newblock {\em Math. USSR-Izv}, 17:423--453, 1981.

\bibitem[AS06]{AbbondandoloSchwartz:cotangentloop}
Alberto Abbondandolo and Matthias Schwarz.
\newblock On the {F}loer homology of cotangent bundles.
\newblock {\em Comm. Pure Appl. Math.}, 59(2):254--316, 2006.

\bibitem[AS10]{SeidelAbouzaid:viterbo}
Mohammed Abouzaid and Paul Seidel.
\newblock An open string analogue of {V}iterbo functoriality.
\newblock {\em Geom. Topol.}, 14(2):627--718, 2010.

\bibitem[BEH{\etalchar{+}}03]{BEHWZ:compactnessfieldtheory}
F.~Bourgeois, Y.~Eliashberg, H.~Hofer, K.~Wysocki, and E.~Zehnder.
\newblock Compactness results in symplectic field theory.
\newblock {\em Geom.Topol.}, 7:799--888, 2003, arXiv:SG/0308183.

\bibitem[BTE09]{eliashberg:symplectichomology}
F.~Bourgeois and Yakov~Eliashberg Tobias~Ekholm.
\newblock Effect of legendrian surgery.
\newblock pages 1--79, 2009, arXiv:SG/0911.0026.

\bibitem[CE]{CieliebakEliashberg:symplecticgeomofsteinmflds}
K.~Cieliebak and Y.~Eliashberg.
\newblock Symplectic geometry of {S}tein manifolds.
\newblock {\em In preparation}.

\bibitem[CFHW96]{CieliebakFloerHoferWysocki:SymhomIIApplications}
K.~Cieliebak, A.~Floer, H.~Hofer, and K.~Wysocki.
\newblock Applications of symplectic homology {II}:stability of the action
  spectrum.
\newblock {\em Math. Z}, 223:27--45, 1996.

\bibitem[DS94]{DostoglouSalamon:instantonscurves}
S.~Dostoglou and D.~Salamon.
\newblock Self-dual instantons and holomorphic curves.
\newblock {\em Ann. of Math. (2)}, 139:581--640, 1994.

\bibitem[EG91]{EliahbergGromov:convexsymplecticmanifolds}
Y.~Eliashberg and M.~Gromov.
\newblock Convex symplectic manifolds.
\newblock In {\em Several complex variables and complex geometry Part 2},
  volume~52, pages 135--162. Amer. Math. Soc., Providence, RI., 1991.

\bibitem[FHS95]{FHS:transversalitysymplectic}
A.~Floer, H.~Hofer, and D.~Salamon.
\newblock Transversality in elliptic {M}orse theory for the symplectic action.
\newblock {\em Duke Math.J.}, 80:251--292, 1995.

\bibitem[FHT82]{FelixHalperinThoamas:homotopylie}
Yves F{\'e}lix, Stephen Halperin, and Jean-Claude Thomas.
\newblock The homotopy {L}ie algebra for finite complexes.
\newblock {\em Inst. Hautes \'Etudes Sci. Publ. Math.}, (56):179--202 (1983),
  1982.

\bibitem[FHT01]{FelixHalperinThomas:rationalhomotopytheory}
Yves F{\'e}lix, Stephen Halperin, and Jean-Claude Thomas.
\newblock {\em Rational homotopy theory}, volume 205 of {\em Graduate Texts in
  Mathematics}.
\newblock Springer-Verlag, New York, 2001.

\bibitem[FT82]{FelixThomasRadius}
Y.~F{\'e}lix and J.-C. Thomas.
\newblock The radius of convergence of {P}oincar\'e series of loop spaces.
\newblock {\em Invent. Math.}, 68(2):257--274, 1982.

\bibitem[Gro78]{Gromov:homotopical}
Mikhael Gromov.
\newblock Homotopical effects of dilatation.
\newblock {\em J. Differential Geom.}, 13(3):303--310, 1978.

\bibitem[Hir64]{hironaka:resolution}
Heisuke Hironaka.
\newblock Resolution of singularities of an algebraic variety over a field of
  characteristic zero. {I}, {II}.
\newblock {\em Ann. of Math. (2) 79 (1964), 109--203; ibid. (2)}, 79:205--326,
  1964.

\bibitem[Lam01]{Lambrechts:Betti}
Pascal Lambrechts.
\newblock The {B}etti numbers of the free loop space of a connected sum.
\newblock {\em J. London Math. Soc. (2)}, 64(1):205--228, 2001.

\bibitem[McL09]{McLean:symhomlef}
M.~McLean.
\newblock Lefschetz fibrations and symplectic homology.
\newblock {\em Geom. Topol.}, 13(4):1877--1944, 2009.

\bibitem[MS98]{McduffSalamon:sympbook}
Dusa McDuff and Dietmar Salamon.
\newblock {\em Introduction to symplectic topology}.
\newblock Oxford Mathematical Monographs. The Clarendon Press Oxford University
  Press, New York, second edition, 1998.

\bibitem[Oan04]{Oancea:survey}
A.~Oancea.
\newblock A survey of {F}loer homology for manifolds with contact type boundary
  or symplectic homology.
\newblock {\em Ensaios Mat.}, 7, 2004, arXiv:SG/0403377.

\bibitem[Rit98]{Ritter:transfer}
Alexander Ritter.
\newblock Topological quantum field theory structure on symplectic cohomology.
\newblock pages 1--60, 1998, arXiv:SG/1003.1781.

\bibitem[Sei08]{Seidel:biasedview}
P.~Seidel.
\newblock A biased view of symplectic cohomology.
\newblock {\em Current Developments in Mathematics}, 2006:211--253, 2008.

\bibitem[SW06]{salamonweber:loop}
D.~A. Salamon and J.~Weber.
\newblock Floer homology and the heat flow.
\newblock {\em Geom. Funct. Anal.}, 16(5):1050--1138, 2006.

\bibitem[Tot03]{Totaro:complexifications}
Burt Totaro.
\newblock Complexifications of nonnegatively curved manifolds.
\newblock {\em J. Eur. Math. Soc. (JEMS)}, 5(1):69--94, 2003.

\bibitem[Vit96]{Viterbo:functorsandcomputations2}
C.~Viterbo.
\newblock Functors and computations in {F}loer homology with applications, part
  {II}.
\newblock {\em Preprint.}, 1996.

\bibitem[VP84]{ViguePoirrier:homotopie}
Micheline Vigu{\'e}-Poirrier.
\newblock Homotopie rationnelle et croissance du nombre de g\'eod\'esiques
  ferm\'ees.
\newblock {\em Ann. Sci. \'Ecole Norm. Sup. (4)}, 17(3):413--431, 1984.

\end{thebibliography}

\end{document}